\documentstyle[12pt,makeidx,amssymb]{article}
\def\eqn#1{\def\theequation{#1}}

\input psfig

\begin{document}

\begin{center} 
{\Large \bf  The ``spectral'' decomposition for one-dimensional maps} 
\end{center} 
\vspace{.2in}

\begin{center} 
Alexander M. Blokh\\
Department of Mathematics, Wesleyan University\\
Middletown, CT 06459-0128, USA
\end{center}
\vspace{.2in}

\noindent 
{\bf Abstract.}
{\it We construct the ``spectral'' decomposition of the sets $\overline{Per\,f},$ $\omega(f)=\cup\omega(x)$ and $\Omega(f)$
for a continuous map $f:[0,1]\rightarrow [0,1]$. Several corollaries are obtained; the main ones describe the generic
properties of $f$-invariant measures, the structure of the set $\Omega(f)\setminus \overline{Per\,f}$ and the
generic limit behavior of an orbit for maps without wandering intervals.
The ``spectral'' decomposition for piecewise-monotone maps is deduced from the Decomposition Theorem. Finally 
we explain how to extend the results of the present paper for a continuous map of a one-dimensional branched manifold into itself.}
\vspace{.2in} 

\baselineskip=20pt

\noindent
{\large \bf 1.  Introduction and main results}
\vspace{.2in}

\noindent
{\bf 1.0. Preliminaries }

	Let $T:X \rightarrow X$ be a continuous map of a compact space into itself (in what follows we consider continuous maps
only). For $x\in X$ the set $orb\,x \equiv \{T^ix:i\ge 0\}$\index{$orb\,x$} is called {\it the orbit of $x$} or 
{\it the $x$-orbit}\index{orbit of a point}.
The set $\omega (x)$\index{$\omega (x)$} 
of all limit points of the $x$-orbit is called {\it the $\omega$-limit set of $x$} or {\it the limit set
of $x$}\index{$\omega$-limit set of a point}. Topological dynamics studies the properties of limit sets.  
Let
us define some objects playing an important role here. A point $x\in X$ is called {\it non-wandering}\index{non-wandering point} 
if for any open $U\ni x$ there exists $n>0$ such that $T^nU\cap U\neq\emptyset$.  
The set $\Omega (T)$\index{$\Omega(T)$} of all non-wandering points is called 
{\it the non-wandering set}\index{non-wandering set}; clearly, $\Omega (T)$ is closed. 

	Let us give an important example. A point $p\in X$ is called {\it periodic}\index{periodic point} 
if $T^np=p$ for some $n\in \Bbb N$. 
Such an $n$ is called {\it a period of $p$}\index{period of a periodic point} 
and the set $orb\,p=\bigcup_{i\ge 0} T^ip$ is called {\it a cycle}\index{cycle}. The set of all
periodic points of $T$ is denoted by $Per(T)$\index{Per\,T}. Clearly, {\it periodic points are non-wandering}.

	We denote the set $\bigcup_{x\in X} \omega (x)$ by $\omega (T)$\index{$\omega(T)$}. 
The following assertion explains the role of the set
$\Omega(T)$.

\noindent
{\bf Assertion 1.1.}{\it For any open set $U\supset \Omega(T)$ and a point $x\in X$ there exists $N$ such that $T^nx \in U$ for all $n>N$,
 and so $\omega(T) \subset \Omega(T)$.}

	Sometimes it is important to know where a point $x\in X$ spends not all the time but almost all the time. The following
definition is useful for considering this problem: a point $x\in X$ is called {\it recurrent}\index{recurrent point} 
if $x\in \omega(x)$. The set of all
recurrent points is denoted by $R(T)$\index{$R(T)$}. 
The set $\overline {R(T)}\equiv C(T)$\index{$C(T)$} is called {\it the center}\index{center of a map} 
of $T$ (here $\overline Z$\index{$\overline Z$} is
{\it the closure} of the set $Z$).

\noindent 
{\bf Assertion 1.2} (see, e.g., [Ma]). {\it For any open $U\supset C(T)$ and $x\in X$ the following property holds:
$\displaystyle \lim_{n \rightarrow \infty} card\{i\le n:T^ix \in U\} \cdot n^{-1}=1$.}

	Let us summarize the connection between the sets $Per\,T,R(T),C(T),\omega(T)$ and $\Omega(T)$ as follows:
\eqn{1.1}
\begin{equation}
Per\,T \subset R(T) \subset \omega(T) \subset \Omega(T)
\end{equation}
\eqn{1.2}
\begin{equation}
\overline {Per\,T}\subset \overline {R(T)}=C(T) \subset \omega(T) \subset \Omega(T)
\end{equation}

	It is useful to split the sets $\Omega(T)$ and $\omega(T)$ into components such that for any $x\in X$ the set $\omega (x)$
belongs to one of them. The remarkable example of such a splitting is the famous Smale spectral decomposition theorem [S] 
(see also [B4]). The aim
of this paper is to show that {\it in the one-dimensional case for any continuous map there exists the decomposition which is in a sense
analogous to that of Smale}.

\noindent 
{\bf 1.1. Historical remarks }

	We start with the history of the subject. From now on fix an arbitrary continuous map $f:[0,1]\rightarrow [0,1]$. 
Speaking of maximality, minimality and ordering among sets we mean that sets are ordered by inclusion. The following definitions are
due to A.N.Sharkovskii [Sh3-6]. Let $\omega(x)$ be an infinite limit set maximal among all limit sets. The set $\omega(x)$
is called {\it a set of genus 1}\index{set of genus 1} 
if it contains no cycles; otherwise it is called {\it a set of genus 2}\index{set of genus 2}. A maximal among
limit sets cycle is called {\it a set of genus 0}\index{set of genus 0} (see [Bl4]; periodic attractors and isolated periodic repellers
are the most important and well-known examples of sets of genus 0).   

	In [Sh3-6] A.N.Sharkovskii has in fact constructed the decomposition of the set $\omega(f)$ into sets of genus 0,1 and 2.
He studied mostly properties of partially ordered family of limit sets belonging to a maximal limit set. Furthermore, he obtained
a number of fundamental results on properties of the sets $\Omega(f),\omega(f),C(f)$ and $\overline {Per\,f}$. Here we formulate
some Sharkovskii's results we need.

\noindent
{\bf Theorem Sh1 [Sh2].} {\it $C(f)\equiv\overline {R(f)}=\overline {Per\,f}=\Omega(f|\Omega(f))$.}

\noindent
{\bf Theorem Sh2 [Sh5].}{\it A point $x$ belongs to $\omega(f)$ if and only if at least 
one of the following properties holds:

	1) for any $\varepsilon>0$ there exists $n>0$ such that $(x-\varepsilon,x)\cap f^n(x-\varepsilon,x)\neq\emptyset$;

	2) for any $\varepsilon>0$ there exists $n>0$ such that $(x+\varepsilon,x)\cap f^n(x+\varepsilon,x)\neq\emptyset$;

	3) $x\in Per\,f$.

	In particular, $\omega(f)$ is closed and so $\overline {Per\,f} \subset \omega(f)$.}

	The main idea of the proofs here is to consider a special kind of recurrence which may occur for maps
of the interval and also to use the following

\noindent
{\bf Property C.} {\it If $I=[a,b]$ is an interval and either $fI\supset I$
or $fI\subset I$ or points $a$ and $b$ move under the first
iteration of $f$ in different directions then there is $y\in I$ such that $fy=y$. }

	 We illustrate this approach considering Theorem Sh1.
Indeed, let $U$ be a complementary to $\overline {Per\,f}$ interval. Then by Property C for any $n$ either $f^nx>x\, (\forall x\in U)$
or $f^nx<x\, (\forall x\in U)$. Suppose that for some $n$ and $x\in U$ we have $f^nx>x, f^nx\in U$. 
Then $f^n(f^nx)>f^nx>x$, i.e. $f^{2n}x>x$; moreover,
if $f^{kn}x>x$ then $f^{(k+1)n}x=f^{kn}(f^nx)>f^nx>x$ which proves that $f^{in}x\ge f^nx>x$ for all $i$. Now suppose that there exists
$y\in U$ and $m$ such that $f^my\in U$ and $f^my<y$. Then by the same arguments $f^{jm}y<y$ for any $j$. It implies that
$f^{mn}x>x$ and $f^{mn}y<y$; so by Property~C there is a periodic point in the interval $(x,y)\subset U$, which is 
a contradiction. 

	Hence if $z\in U,\,k$ is the minimal number such that $f^kz \in U$ and, say, $z<f^kz$ then 
for any $l>k$ the fact that $f^lz\in U$ implies that $f^lz=f^{l-k}(f^kz)>f^kz>z$. So by the definition there is no recurrent point of $f$
inside $U$. Moreover, if we consider a sufficiently small neighborhood $V\subset U$ of a point $\zeta\in \Omega(f)\cap U$
and repeat the arguments from above we will see that $f^n\zeta \not \in U$ for any $n$. 
So we have proved that $R(f)\cap U=\emptyset$ and thus  $\overline {R(f)}=\overline {Per\,f} $; we have proved also that
every point from $\Omega(f)\cap U$ never returns to $U$ and hence is not a non-wandering point in $\Omega(f)$, i.e.
$\overline {Per\,f}=\Omega(f|\Omega(f))$. It completes the sketch of the proof of Theorem Sh2.

	One of the most well-known and surprised results about one-dimensional dynamics is, perhaps, the famous Sharkovskii theorem.  
To state it let us consider the set $\Bbb N$ of positive integers with the following 
{\it Sharkovskii ordering}\index{Sharkovskii ordering}:

\eqn{*}
\begin{equation}
3\prec 5\prec 7\prec \ldots \prec 2\cdot 3\prec 2\cdot 5\prec 2\cdot 7\ldots \prec 2^3 \prec 2^2 \prec 2\prec 1
\end{equation}

\noindent
{\bf Theorem Sh3 [Sh1].} {\it Let $m\prec n$ and $f$ have a cycle of minimal period $m$. Then $f$ has a cycle of minimal period $n$.}

	We say that {\it $m$ is $\prec$-stronger than $n$} if $m\prec n$ and $m\neq n$. 
We say that {\it $f$ is a map of type $m$}\index{interval map of type $m$} [Bl7] 
if the $\prec$-strongest period of cycles of $f$ is $m$; in
other words, $m$ is the largest period which appears in terms of the Sharkovskii ordering. Such a period does not
exist if the periods of cycles of $f$ are exactly $1,2,2^2,2^3,\ldots$; then say that $f$ is {\it of type} $2^\infty$
\index{interval map of type $2^\infty$}[Bl1].

	For piecewise-monotone continuous maps the splittings of the sets $\Omega(f),\omega(f)$ and $C(f)$ in fact analogous
to that of Sharkovskii were constructed later by using the different technique (see the articles of Jonker and Rand [JR1-JR2],
Nitecki [N2] and the books of Preston [P1-P2]). For piecewise-monotone maps with a finite number of discontinuities the construction
of the splitting is due to Hofbauer [H1-H2].

\noindent
{\bf 1.2. A short description of the approach presented}

	The approach in this paper is different from the one of Sharkovskii and the ``piecewise-monotone'' approach; it is
based on the author's articles [Bl1-Bl12].
 First we need more definitions. Let $T:X\rightarrow X$ and $F:Y\rightarrow Y$ be maps of compact spaces. If there exists 
a {\it surjective} map $\phi:X\rightarrow Y$ such that $\phi \circ T=F \circ \phi$ then it is said that $\phi$ {\it semiconjugates} 
$T$ to $F$ and $\phi$ is called {\it a semiconjugation}\index{semiconjugation} between $T$ and $F$; 
if $\phi$ is a homeomorphism then it is said that $\phi$ {\it conjugates} $T$ to $F$ and $\phi$ is called
{\it a conjugation}\index{conjugation} between $T$ and $F$.

	Roughly speaking, our approach to one-dimensional maps is the following: we propose models for different kinds
of limit sets, study the properties of models, extend the properties to the limit sets and the map itself, and also obtain 
some corollaries. In the rest of this Section we formulate the main results of the present paper. The proofs will be given in 
Sections 2-11. 
At the end of this Section we apply the results to piecewise-monotone maps and explain how to extend the decomposition to
continuous maps of one-dimensional branched manifolds.

	An interval $I$ is called {\it periodic (of period $k$)}\index{periodic interval} 
or {\it $k$-periodic} if $J,\ldots,f^{k-1}J$ are pairwise
disjoint and $f^kJ=J$ (if it is known only that $f^kJ\subset J$ then $J$ is called 
{\it a weakly periodic interval}\index{weakly periodic interval}). The set
$\bigcup^{k-1}_{i=0}f^iJ\equiv orb\,J$\index{$orb\,J$} is called {\it a cycle of intervals}\index{cycle of intervals} 
if $J$ is periodic and {\it a weak cycle of intervals}\index{weak cycle of intervals}
if $J$ is weakly periodic (here $k$ is a period of $J$).

	Let us explain briefly how we will classify limit sets. Fix an infinite set $\omega(x)$ and consider a family $\cal A$
of all cycles of intervals $orb\,I$ such that $\omega(x)\subset orb\,I$. There are two possibilities.

	1) {\it Periods of sets $orb\,I \in \cal A$ are not bounded.} Then there exist ordered cycles of intervals
containing $\omega(x)$ with periods tending to infinity. It allows us to semiconjugate $f|\omega(x)$ to a transitive translation
in a compact group and implies many properties of $f|\omega(x)$. The set $\omega(x)$ corresponds to a Sharkovskii's set of
genus 1.

	2) {\it Periods of sets $orb\,I \in \cal A$ are bounded.} Then there exists a minimal cycle of intervals
$orb\,J \in \cal A$. It is easy to see that all points $y\in \omega(x)$ have the following property: if $U$ is a neighborhood
of $y$ in $orb\,J$ then $\overline {orb\,U}=orb\,J$ (otherwise $\overline {orb\,U}$ generates a cycle of intervals $orb\,K$
such that $\omega(x) \subset orb\,K \subset orb\,J,\,orb\,K\neq orb\,J$ which is a contradiction).
The idea is to consider all the points $z\in orb\,J$ with this property. They form a set $B$ which is another example of a maximal
limit set. The set $B$ is a Sharkovskii's set of genus 2. 

\noindent
{\bf 1.3. Solenoidal sets}

	Let us proceed more precisely. Let $T:X\rightarrow X$ be a map of a compact metric space $(X,d)$ into itself. The map $T$
is said to be {\it transitive}\index{transitive map} if there exists an $x$ such that $\omega(x)=X$, to be 
{\it minimal}\index{minimal} if for any
$x\in X$ we have $\omega(x)=X$, to be {\it topologically mixing }\index{topologically mixing, mixing map} 
or simply {\it mixing} if for any open $U,V$ there exists 
an $N$ such that $T^nU\cap V\neq \emptyset$ for any $n>N$.

	We will also need the definition of the topological entropy; the notion was introduced in
[AKMcA] but we give the definition following Bowen[B1].
A set $E\subset X$ is said to be $(n,\varepsilon)$-separated\index{$(n,\varepsilon)$-separated set} 
if for any two distinct points $x,y \in E$ there exists $k, 0\le k<n$
such that $d(T^kx,T^ky)>\varepsilon$. By $S_n(\varepsilon)$ 
we denote the largest cardinality of an $(n,\varepsilon)$-separated subset of $X$.
Let $S(\varepsilon)\equiv \limsup n^{-1}\cdot \ln S_n(\varepsilon)$. 
Then the limit $h(T)=\lim_{\varepsilon\rightarrow 0} S(\varepsilon)$ exists
and is called {\it the topological entropy}\index{topological entropy} of $T$ (see [B1]). Now let us turn back to interval maps.

	Let $I_0\supset I_1\supset \ldots$ be periodic intervals with periods $m_0,m_1,\ldots$. Obviously $m_{i+1}$ is a multiple
of $m_i$ for all $i$. If $m_i \rightarrow \infty$ then the intervals $\{I_j\}^\infty_{j=0}$ are said to be 
{\it generating}\index{generating intervals} and
any invariant closed set $S\subset Q=\bigcap_{j\ge 0}orb\,I_j$ is called {\it a solenoidal set}\index{solenoidal set}; 
if $Q$ is nowhere dense then we
call $Q$ {\it a solenoid}\index{solenoid}. In the sequel we use the following notation:

	$\bigcap_{j\ge 0}orb\,I_j \equiv Q(\{I_j\}^\infty_{j=0}) \equiv Q$\index{$Q(\{I_j\}^\infty_{j=0}) \equiv Q$};

	$Q\cap \overline {Per\,f} \equiv S_p(Q) \equiv S_p$\index{$S_p(Q) \equiv S_p$};

	$Q\cap \omega(f) \equiv S_\omega(Q) \equiv S_\omega$\index{$S_\omega(Q) \equiv S_\omega$};

	$Q\cap \Omega(f) \equiv S_\Omega(Q) \equiv S_\Omega$\index{$S_\Omega(Q) \equiv S_\Omega$}.

	Observe that $S_p \subset S_\omega \subset S_\Omega$ and all these sets are invariant and closed ( for $S_\omega$ it follows
from Theorem Sh2).

	One can use a transitive translation in an Abelian zero-dimensional infinite group as a model for the map on a solenoidal set.
Namely, let $D=\{n_i\}^\infty_{i=0}$ be a sequence of integers, $n_{i+1}$ be a multiple of $n_i$ for all $i$ and 
$n_i \rightarrow \infty$. Let us consider a subgroup 
$H(D) \subset \Bbb Z$$_{n_{0}} \times \Bbb Z$$_{n_{1}} \times \cdots$, defined
by $H(D)\equiv \{(r_0,r_1,\ldots):r_{i+1} \equiv r_i \pmod {m_i} (\forall i)\}$\index{$H(D)$}. 
Denote by $\tau$ the minimal translation in $H(D)$ by the
element $(1,1,\ldots)$. 

\noindent
{\bf Theorem 3.1[Bl4,Bl7].}
{\it Let $\{I_j\}^\infty_{j=0}$ be generating intervals with periods $\{m_i\}^\infty_{i=0}=D,\, Q=\bigcap_{j\ge 0}orb\,I_j$. Then
there exists a continuous map $\phi:Q \rightarrow H(D)$ with the following properties:

	1) $\phi$ semiconjugates $f|Q$ to $\tau$ (i.e. $\tau \circ \phi=\phi \circ f$ and $\phi$ is surjective);

	2) there exists a unique set $S\subset S_p$ such that $\omega(x)=S$ for any $x\in Q$ and, moreover, $S$ is the set of all limit 
points of $S_\Omega$ and $f|S$ is minimal);

	3) if $\omega(z)\cap Q \neq \emptyset$ then $S \subset \omega(z) \subset S_\omega$;

	4) for any {\bf r}$\in H(D)$ the set $J=\phi^{-1}(${\bf r}$)$ is a connected component of $Q$ and:

	a) if $J=\{a\}$ then $a\in S$;

	b) if $J=[a,b], a\neq b$ then $\emptyset \neq S\cap J\subset S_\Omega\cap J \subset \{a,b\}$;

	5) $S_\Omega \setminus S$ is at most countable and consists of isolated points;

	6) $h(f|Q)=0$.}

	It should be noted that the best known example of a solenoid is the Feigenbaum attractor ([CE],[F]) for which generating
intervals have periods $\{2^i\}^\infty_{i=0}$. If for a solenoid or a solenoidal set generating intervals
have periods $\{2^i\}^\infty_{i=0}$ then we call it {\it $2$-adic}\index{$2$-adic solenoidal set, $2$-adic solenoid}.

\noindent
{\bf 1.4. Basic sets}

	Let us turn to another type of maximal infinite limit set. Let $\{J_i\}^l_{i=1}$ be an ordered collection
of intervals (one can imagine these intervals lying on the real line in such a way that $J_1<J_2<\ldots<J_l,\, J_i\cap J_r=\emptyset$
for $i \neq r$); set $K=\bigcup^l_{i=1}J_i$. A continuous map $\psi: K \rightarrow K$ which permutes the intervals $\{J_i\}^l_{i=1}$ 
cyclically is called {\it non-strictly periodic}\index{non-strictly periodic map} 
(or {\it $l$-periodic}). Note that this term concerns a map, not an interval;
we speak of non-strictly periodic maps to distinguish them from {\it periodic maps} which are traditionally those with all points
periodic. An example of a non-strictly periodic map is a map of the interval restricted on a weak cycle of intervals.

	Now let $\psi:K \rightarrow K$ and   $\psi':K' \rightarrow K'$ be non-strictly $l$-periodic maps (so that $K$ and $K'$ are
unions of $l$ intervals). Let 
$\phi:K \rightarrow K'$ be a monotone semiconjugation between $\psi$ and $\psi'$ and $F\subset K$ be a $\psi$-invariant
closed set such that $\phi (F)=K', $ for any $x\in K'$ we have $int\,\phi^{-1}(x) \cap F=\emptyset$ and so 
$\phi^{-1}(x)\cap F \subset \partial \phi^{-1}(x),\,1\le card\{\phi^{-1}(x)\cap F\} \le 2$. Then we say that $\phi$
{\it almost conjugates} $\psi |F$ to $\psi'$ or $\phi$ is an
{\it almost conjugation}\index{almost conjugation} between $\psi |F$ and $\psi'$. Remark that here 
$int\, Z$\index{$int\, Z$} is {\it an interior} of a set $Z$ and
$\partial Z$\index{$\partial Z$} is {\it a boundary } of a set $Z$.

	Finally let $I$ be an $n$-periodic interval, $orb\,I=M$. Consider a set 
$\{x\in M:$ for any relative neighborhood $U$ of $x$ in $M$ 
we have $\overline {orb\,U}=M\}$; it is easy to see that this is a closed invariant set. It 
is called {\it a basic set}\index{basic set} and denoted by $B(M,f)$\index{$B(M,f)$} 
provided it is infinite. Now we can formulate

\noindent
{\bf Theorem 4.1[Bl4,Bl7].}
{\it Let $I$ be an $n$-periodic interval, $M=orb\,I$ and $B=B(M,f)$ be a basic set. Then there exist 
a transitive non-strictly $n$-periodic map $g:M' \rightarrow M'$ and a monotone map $\phi: M \rightarrow M'$
such that $\phi$ almost conjugates $f|B$ to $g$. Furthermore, $B$ has the following properties:

	a) $B$ is a perfect set;

	b) $f|B$ is transitive;

	c) if $\omega(z) \supset B$ then $\omega(z)=B$ (i.e. $B$ is a maximal limit set);

	d) $h(f|B) \ge \ln 2 \cdot (2n)^{-1}$;

	e) $B\subset \overline {Per\,f}$;

	f) there exist an interval $J\subset I$, an integer $k=n$ or $k=2n$ and a set $\widetilde B=\overline {int\,J\cap B}$ such
that $f^kJ=J,\:f^k \widetilde B=\widetilde B,\:f^i \widetilde B \cap f^j \widetilde B$ contains no more than $1$ point
($0\le i<j<k$), $\bigcup^{k-1}_{i=0}f^i \widetilde B=B$ and $f^k|\widetilde{B}$ is almost conjugate to a mixing interval map
(one can assume that if $k=n$ then $I=J$).}

	So we use transitive non-strictly periodic maps as models for the map $f$ on basic sets. Note that Theorems 3.1
and 4.1 allow to establish the connection between sets of genus 1 and solenoidal sets, sets of genus 2 and basic sets 
(see Assertion 4.2 in Section 4). Moreover, Theorems 3.1 and 4.1 easily imply that sets of genus 0 and limit solenoidal sets may be
characterized as those $\omega(x)$ for which the inclusion $\omega(y)\supset \omega(x)$ implies that $h(f|\omega(y))=0$ for any $y$ 
(see Assertion 4.3 in Section 4).

%
%
%
%

	Now we can construct the ``spectral'' decomposition for the sets $\overline {Per\,f}$ and $\omega(f)$. However, 
to extend the decomposition to the set $\Omega(f)$ we need the following definition. Let $B=B(orb\,I,f)$ be a basic set and 
$A$ be the set of all endpoints $x$ of the intervals of $orb\,I$ with the following properties:

	1) $x\in \Omega(f)$;

	2) there exists an integer $n$ such that $f^nx \in B$ and if $m$ is the least such integer then
$x,fx,\ldots,f^{m-1}x \not \in int(orb\,I)$.

	We denote the set $B\cup A$ by $B'(orb\,I,f)$\index{$B'(M,f)$} and call it {\it an $\Omega$-basic set}\index{$\Omega$-basic set}.

	Let us consider an example of an $\Omega$-basic set (cf. [Sh2], [Y]). Construct a map $f:[0,1] \rightarrow [0,1]$ in the
following way: fix $6$ points $c_0=0<c_1<\ldots<c_5=1$, define $f|\bigcup^5_{i=0}c_i$ and then extend $f$ on each interval 
$[c_i,c_{i+1}]$ as a linear function. Namely:

	1) $c_0=0,fc_0=2/3$;

	2) $c_1=1/3,fc_1=1$;

	3) $c_2=1/2, fc_2=5/6$;

	4) $c_3=2/3, fc_3=1/6$;

	5) $c_4=5/6, fc_4=5/6$;

	6) $c_5=1,fc_5=1$.

	It is easy to see that the interval $I=[1/6,1]$ is $f$-invariant and the point $5/6$ is fixed. Let us show
that there exists a basic set $B=B(orb\,I,f)$ and $1/2,5/6 \in B$.
For the moment we know nothing about the cardinality of the set $\{x\in orb\,I:$ for any relative neighborhood $U$ of $x$ in $orb\,I$ 
we have $\overline {orb\,U}=orb\,I\}$, so let us denote this set by $L$; by the definition of a basic set
we need to prove that $L$ is infinite. Indeed, any left semi-neighborhood of $5/6$ covers
the whole interval $I$ after some iterations of $f$, so $5/6 \in L$. On the other hand 
$f(1/2)=5/6$ and $f$-image of any right semi-neighborhood of $1/2$ covers some left semi-neighborhood
of $5/6$. So $1/2 \in L$ as well. But it is easy to see that there are infinitely many points $z\in (1/2,5/6)$ such that $f^nz=1/2$
for some $n$ and $f^n$-image of any neighborhood of $z$ covers some right semi-neighborhood of $1/2$  which implies
that $z \in L$; so $L$ is infinite and by the definition $L=B(orb\,I,f)=B$ is a basic set. Furthermore, the map $f$ coincides with
the identity on $[5/6,1]$, at the same time $f[1/6,1/2]=[5/6,1]=f[5/6,1]$ and $f$-image of any right semi-neighborhood of $1/6$
is some right semi-neighborhood of $5/6$. So by the definition we see that
there are no points of $B$ in $[1/6,1/2)$; in particular, $1/6 \not \in B$.

	Moreover, it is easy to see that there are no periodic points of $f$ in $[0,1/2)$. Indeed, there are no periodic points in 
$[1/6,1/2]$ because $f[1/6,1/2]=[5/6,1]=f[5/6,1]$. On the other hand there are no periodic points in $[0,1/6)$ because 
$f[0,1/6)\subset [1/6,1]=f[1/6,1]$. So $Per\,f\cap [0,1/2]=\emptyset$. 

	Now let us show that $1/6 \in B' \setminus B$ where $B'=B'(orb\,I,f)$. Indeed, we have already seen that $f(1/6)=5/6\in B$.
So by the properties of the point $5/6$ established above we see that for any open $U \ni 1/6$ there exists 
$m$ such that $f^mU \supset I=[1/6,1]$. It proves that $1/6 \in \Omega(f)$. Now the definition implies that $1/6 \in B' \setminus B$. 
It remains to note that by the definition $B' \setminus B$ consists only of some endpoints of intervals from $orb\,I$, i.e.
in our case of some  of the points $1/6,1$. Clearly, $1 \not \in B'$ and so $\{1/6\}=B' \setminus B$.

\noindent
{\bf 1.5. The decomposition and main corollaries}

	Now we can formulate the Decomposition Theorem. Let us denote by $X_f$\index{$X_f$} the union of all limit sets of genus $0$.

\noindent
{\bf Theorem 5.4 (Decomposition Theorem)[Bl4,Bl7].}
{\it Let $f:[0,1] \rightarrow [0,1]$ be a continuous map. Then there exist an at most countable family of pairs of basic and
$\Omega$-basic sets $\{B_i \subset B'_i\}$ and a family of collections of solenoidal sets 
$\{S^{(\alpha)}\subset S^{(\alpha)}_p\subset S^{(\alpha)}_\omega\subset S^{(\alpha)}_\Omega\subset Q^{(\alpha)}\}_{\alpha\in A}$ 
with the following properties:

	1) $\Omega(f)=X_f\cup (\bigcup_\alpha S^{(\alpha)}_\Omega) \cup (\bigcup_i B'_i)$;

	2) $\omega(f)= X_f\cup (\bigcup_\alpha S^{(\alpha)}_\omega) \cup (\bigcup_i B_i)$;

	3) $\overline {Per\,f}= X_f\cup (\bigcup_\alpha S^{(\alpha)}_p) \cup (\bigcup_i B_i)$;

	4) the set $S^{(\alpha)}_\Omega \setminus S^{(\alpha)}$ is at most countable set of isolated points, the set \linebreak 
$\{\alpha: int\,Q^{(\alpha)} \neq \emptyset\}$ is at most countable and $S^{(\alpha)}=Q^{(\alpha)}$ for all
other $\alpha \in A$;

	5) intersections in this decomposition are possible only between different basic or $\Omega$-basic sets, each three
of them have an empty intersection and the intersection of two basic or two $\Omega$-basic sets is finite.}

	Note that in statement 5) of Decomposition Theorem we do not take into account intersections between a basic set
and an $\Omega$-basic set with the same subscript and also between different solenoidal sets with the same superscript.

	The Decomposition Theorem in the full formulation is quite cumbersome but the idea is fairly clear and may be expressed in the
following rather naive version of the Decomposition Theorem: 

	{\it for any continuous map $f:[0,1]\rightarrow[0,1]$ the non-wandering set $\Omega(f)$ and related sets 
(like $\omega(f)$ and $\overline {Per\,f}$) are unions of the set $X_f$, solenoidal sets and basic sets.} 
 
	The main corollaries of this picture of dynamics are connected with the following problems.

	1) What is the generic limit behavior of orbits for maps without wandering intervals (we call an interval 
$I$ {\it wandering}\index{wandering interval} if $f^n\cap f^mI=\emptyset$ for $n>m\ge 0$ and $I$ does not tend to a cycle)(Section 6)?

	2) What is the related structure of the set $\Omega(f),\; \omega(f)$ and $\overline {Per\,f}$(Section 7)?

	3) How does dynamics of a map depend on its set of periods of cycles (Section 9)?

	4) What are the generic properties of invariant measures (Section 10)?

	Note that in order to study the generic properties of invariant measures we establish in Section 8 some important properties
of transitive and mixing interval maps. In Section 11 we also investigate the connection between the results of the present paper
and some recent results of Block and Coven [BC] and Xiong Jincheng [X]. 

	In the following subsections 1.6-1.10 we outline the way we are going to obtain the aforementioned corollaries of the 
Decomposition Theorem. In subsections 1.11-1.12 we describe the decomposition for piecewise-monotone interval maps and
in subsection 1.13 we discuss further generalizations. 

\noindent
{\bf 1.6. The limit behavior and generic limit sets for maps without wandering intervals}

	In this subsection we describe the results of Section 6. We start 
with the reformulation of the Decomposition Theorem for maps without wandering intervals. Namely, Theorem 3.1 
implies that if a map $f$ does not have wandering intervals then in the notation from the Decomposition Theorem for any
$\alpha \in A$ we have  $\{S^{(\alpha)}=S^{(\alpha)}_p=S^{(\alpha)}_\omega=S^{(\alpha)}_\Omega=Q^{(\alpha)}\}_{\alpha\in A}$;
in other words all solenoidal sets are in fact solenoids (recall that solenoids are nowhere dense intersections of
cycles of generating intervals which in particular implies that the map on a solenoid is topologically conjugate with the 
translation in the corresponding group). This makes the formulation of the Decomposition Theorem easier, so let
us reformulate it in this case. 

\noindent
{\bf Decomposition Theorem for maps without wandering intervals.}
{\it Let $f$ be a continuous interval map without wandering intervals. Then there exist an at most countable 
family of pairs of basic and
$\Omega$-basic sets $\{B_i \subset B'_i\}$ and a family of solenoids  
$\{Q^{(\alpha)}\}_{\alpha\in A}$ 
with the following properties:

	1) $\Omega(f)=X_f\cup (\bigcup_\alpha Q^{(\alpha)}) \cup (\bigcup_i B'_i)$;

	2) $\omega(f)= X_f\cup (\bigcup_\alpha Q^{(\alpha)}) \cup (\bigcup_i B_i)$;

	3) $\overline {Per\,f}= X_f\cup (\bigcup_\alpha Q^{(\alpha)}) \cup (\bigcup_i B_i)$;

	4) intersections in this decomposition are possible only between different basic or $\Omega$-basic sets, each three
of them have an empty intersection and the intersection of two basic or two $\Omega$-basic sets is finite.}

	A set $A$ which is a countable 
intersection of open subsets of a compact metric space $X$ is said to be a $G_\delta$-set\index{$G_\delta$-set}. 
A set $G$ containing a dense $G_\delta$-set is said to be {\it residual}\index{residual set}.
A property which holds for a residual subset of a compact metric space is said to be 
{\it topologically generic}\index{topologically generic}. One of the corollaries of the aforementioned version of the 
Decomposition Theorem is the description of generic limit sets for maps without wandering intervals. First let us explain why   
the concept of wandering interval appears naturally while studying the problem in question. Indeed,  
consider a pm-map $g$ without {\it flat spots} (i.e. intervals $I$ such that $fI$ is a point). Take any point $x\in [0,1]$ with
an infinite orbit not tending to a cycle. 
Then instead of points $z$ from the set $\bigcup^\infty_{i,j}g^{-i}(g^jx)$ we can ``paste in'' intervals $I(z)$ in
such a way that a new map will have a wandering interval $I(x)$ and that $orb_fI(x)$ will have essentially the same structure as
$orb_gx$. Therefore to consider the problem in question one should forbid the existence of wandering intervals. This remark makes the
following Theorem 6.2 quite natural.

\noindent
{\bf Theorem 6.2}(cf.[Bl1],[Bl8]).
{\it Let $f:[0,1]\rightarrow [0,1]$ be a continuous map without wandering intervals. Then there exists a residual subset
$G\subset [0,1]$ such that for any $x\in G$ one of the following possibilities holds:

	1)$\omega(x)$ is a cycle;

	2) $\omega(x)$ is a solenoid;

	3) $\omega(x)=orb\,I$ is a cycle of intervals}.

\noindent
{\bf Remark.} Note, that possibility 3) of Theorem 6.2 will be essentially specified in Section 10 where we show that in fact generic 
points $x$ for which $\omega(x)=orb\,I$ is a cycle of intervals may be chosen in such a way that the set of all limit measures
of time averages of iterates of $\delta$-measure $\delta_x$ coincides with the set of all invariant measures of $f|orb\,I$ 
(precise definitions will be given in subsection 1.10). 

\noindent
{\bf 1.7. Topological properties of sets $\overline {Per\,f},\: \omega(f)$ and $\Omega(f)$ } 

	In this subsection we summarize the results of Section 7. The following Theorem 7.6 which is the main theorem of Section 7 
describes the structure of the set $\Omega(f)\setminus \overline {Per\,f}$.

\noindent
{\bf Theorem 7.6.}
{\it Let $U=(a,b)$ be an interval complementary to $\overline {Per\,f}$. Then up to the orientation one of the following four possibilities
holds.

	1) $\Omega(f)\cap U=\emptyset$.

	2) $\Omega(f)\cap U=\{x_1<x_2<\ldots<x_n\}$ is a finite set, 
$card(orb\,x_1)<\infty,\ldots ,\linebreak 
card(orb\,x_{n-1})<\infty,\;(\bigcup^{n-1}_{i=1}x_i)\cap \omega(f)=\emptyset$ and there exist periodic
intervals $J_i=[x_i,y_i]$ such that $x_i\in B'(orb\,J_i,f)$ for $1\le i\le n-1$ and $J_i\supset J_{i+1}$ for $1\le i\le n-2$.
Moreover, for $x_n$ there exist two possibilities:
a) $x_n$ belongs to a solenoidal set;
b) $x_n$ belongs to an $\Omega$-basic set $B'(orb\,J_n,f)$ where $J_n=[x_n,y_n]\subset J_{n-1}$.

	3) $\Omega(f)\cap U=(\bigcup^\infty_{i=1}x_i)\cap x,\;x_1<x_2<\ldots,\;x=\lim\,x_i,\;$ and there exist generating intervals
$J_i=[x_i,y_i]$ such that:

	a) $x_i\in B'(orb\,J_i,f),\;card(orb\,x_i)<\infty\:(\forall i)\;$ and $(\bigcup^\infty_{i=1}x_i)\cap \omega(f)=\emptyset$;

	b) $x\in S_\omega(\{orb\,J_i\}^\infty_{i=1})=\omega(f)\cap (\bigcap^\infty_{i=1}orb\,J_i)$.

	4) $\Omega(f)\cap U=\bigcup^\infty_{i=1}x_i,\:x_1<x_2<\ldots,\:\lim\,x_i=b,\:card(orb\,x_i)<\infty\:(\forall i),\linebreak 
(\bigcup^\infty_{i=1}x_i)\cap \omega(f)=\emptyset$ and there exist periodic intervals $J_i=[x_i,y_i]$ such that
$x_i\in B'(orb\,J_i,f),J_i\supset J_{i+1}\,(\forall i)$ and $\bigcap^\infty_{i=1} J_i=b$. Moreover, either periods of $J_i$
tend to infinity, $\{J_i\}$ are generating intervals and $b$ belongs to the corresponding solenoidal set, or periods of $J_i$ do not
tend to infinity and $b$ is a periodic point.

	In any case $card\{\omega(f)\cap U\}\le 1$.}

	We also prove in Section 7 that $\displaystyle \omega(f)=\bigcap_{n\ge 0}f^n\Omega(f)$. Finally, we extend for any 
continuous interval map a result from [Y] where it is proved that if $f$ is a pm-map and 
$x\in \Omega(f)\setminus \overline{Per\,f}$ then there exists $n>0$ and a turning point $c$ such that $f^nc=x$ (the similar
result was obtained in the recent paper [Li], see Theorem 2 there).

\noindent
{\bf 1.8. Properties of transitive and mixing maps}

	Theorem 4.1 implies that properties of a map on basic sets are closely related to properties of 
transitive and mixing interval maps. We investigate these properties in Section 8 and give here a summary of the 
corresponding results.

	The following lemma shows the connection between transitive and mixing maps of the interval.

\noindent
{\bf Lemma 8.3[Bl7].} {\it Let $f:[0,1] \rightarrow [0,1]$ be a transitive map. Then one of the following possibilities holds:

	1) the map $f$ is mixing and, moreover, for any $\varepsilon>0$ and any non-degenerate interval $U$ there exists
$m$ such that $f^nU \supset [\varepsilon, 1-\varepsilon]$ for any $n>m$;

	2) the map $f$ is not mixing and there exists a fixed point $a \in (0,1)$ such that $f[0,a]=[a,1],f[a,1]=[0,a],f^2|[0,a]$ 
and $f^2|[a,1]$ are mixing.

	In any case $\overline {Per\,f}=[0,1].$}

	It turns out that mixing interval maps have quite strong expanding properties: any open interval under iterations of a mixing map
$f:[0,1]\rightarrow [0,1]$ eventually covers any compact subset of $(0,1)$. More precisely,  
let $A(f)\equiv A$  be the set of those from points $0,1$ which have no preimages in $(0,1)$. 

\noindent
{\bf Lemma 8.5[Bl7].}
{\it If $f:[0,1]\rightarrow [0,1]$ is mixing then there are the following possibilities for $A$:

	1) $A=\emptyset$;

	2) $A=\{0\},f(0)=0$;

	3) $A=\{1\},f(1)=1$;

	4) $A=\{0,1\},f(0)=0,f(1)=1$;

	5) $A=\{0,1\},f(0)=1,f(1)=0$.

	Moreover, if $I$ is a closed interval, $I\cap A=\emptyset$, then for any open $U$ there exists $n$ such that $f^mU\supset I$
for $m>n$ (in particular, if $A=\emptyset$ then for any open $U$ there exists $n$ such that $f^nU=[0,1]$.}

	In fact this lemma is one of the basical tools in the  proof of Theorem 8.7 where we show that mixing interval maps 
have the specification property. It is well-known ([Si1-Si2], [DGS]) that this implies a lot of generic properties of invariant
measures of a map, and we will rely on them in the further studying of interval dynamics.

	Let us give the exact definition. Let $T:X \rightarrow X$ be a map of a compact infinite metric space $(X,d)$ into itself.
A dynamical system $(X,T)$ is said to have {\it the specification property}\index{specification property} 
or simply {\it the specification} [B2] if for any
$\varepsilon>0$ there exists an integer $M=M(\varepsilon)$ such that for any $k>1$, for any $k$ points $x_1,x_2,\ldots,x_k\in X$,
for any integers $a_1\le b_1<a_2\le b_2<\ldots <a_k\le b_k$ with $a_i-b_{i-1} \ge M,\,2\le i\le k$ and for any integer $p$ with 
$p\ge M+b_k-a_1$ there exists a point $x\in X$ with $T^px=x$ such that $d(T^nx,T^nx_i)\le \varepsilon$ for $a_i\le n \le b_i,1\le i\le k$.

\noindent
{\bf Theorem 8.7[Bl4,Bl7].}
{\it If $f:[0,1] \rightarrow [0,1]$ is mixing then $f$ has the specification property.}

\noindent
{\bf Remark.} In Section 8 we in fact introduce slightly stronger version of the specification property (i-specification property) 
related to the properties of interval maps and prove that mixing maps of the interval have the i-specification.

\noindent
{\bf 1.9. Corollaries concerning periods of cycles for interval maps}

	Here we formulate two results concerning periods of cycles 
for interval maps which are proved in Section 9. 
We explain also how the famous Misiurewicz theorem on maps with zero entropy is connected with our results.

	Well-known properties of the topological entropy and Theorem Sh1 imply that $h(f)=h(f|\overline {Per\,f})$.   
However, it is possible to get a set $D$ such that
$h(f)=h(f|D)$ using essentially fewer periodic points of $f$. Indeed, let $A\subset \Bbb N$ , $K_f(A)\equiv \{y\in Per\,f:$ minimal period 
of $y$ belongs to $A\}$\index{$K_f(A)$}. 

\noindent
{\bf Theorem 9.1[Bl4,Bl7].} {\it The following two properties of $A\subset \Bbb N$ are equivalent:

	1) $h(f)=h(f|\overline {K_f(A)})$ for any $f$;

	2) for any $k$ there exists $n\in A$ which is a multiple of $k$. }

	In Theorem 9.4 we study how the sets $\Omega(f),\Omega(f^2),\ldots$ vary for maps with a fixed set of periods of cycles.
In [CN] this problem was investigated for arbitrary continuous map of the interval and it was proved that
$\Omega(f)=\Omega(f^n)$ for any odd $n$ and any continuous interval map. The following theorem is related to the results of [CN].

%
%
%
%
%
%
\noindent
{\bf Theorem 9.4[Bl4,Bl8].}
{\it Let $n\ge 0,k\ge 1$ be fixed and $f$ have no cycles of minimal period $2^n(2k+1)$. Then the following statements hold:

	1) if $B=B(orb\,I,f)$ is a basic set and $I$ has a period $m$ then $2^n(2k+1)\prec m\prec 2^{n-1}$;

	2) $\Omega(f)=\Omega(f^{2^n})$;

	3) if $f$ is of type $2^m, 0\le m\le \infty$, then $\Omega(f)=\Omega(f^l)\;(\forall l)$.}

	Note that we use here the Sharkovskii ordering ``$\prec$'' (see ($*$) in the beginning of this Section) in the strict sense
(i.e. $m\prec n$ implies $m\neq n$). The assertion close to statement 1) of Theorem 9.4 was proved in [N2], Theorem 1.10.

	Now we explain the connection between the Misiurewicz theorem on maps with zero entropy and our results. Bowen and Franks 
[BF] have proved that if $f$ is a map of type $m,m\neq 2^n (0\le n\le \infty)$ then $h(f)>0$. The converse was proved first for
pm-maps [MiS] and then for arbitrary continuous maps of the interval into itself [Mi2]. Let us show how to deduce the converse
assertion from our results. 

	If $T:X\rightarrow X$ is a map of a compact metric space $X$ with the specification then there exists $N$
such that for any $n>N$ there exists a periodic point $y\in X$ of a minimal period $n$. Let $f$ be of type $2^n,n\le \infty$.
Suppose that $f$ has a basic set $B=B(orb\,I,f)$. By Theorem 4.1 properties of the map restricted on $B$ are close to 
properties of the corresponding mixing interval map. Furthermore, by Theorem 8.7 this mixing interval map has the
specification, and so by the property of maps
with the specification mentioned above one can find integers $l,k$ such that for any  $m>l$ there exists an $f$-cycle of minimal period
$km$ , which contradicts to the fact that $f$ is of type $2^n$. Thus $f$ has no basic sets and $h(f)=0$ by well-known properties
of the topological entropy  
and Theorem 3.1. Note  that now the Decomposition Theorem implies that infinite limit sets of a map with zero entropy are $2$-adic
solenoidal sets (another proof of this assertion follows from Misiurewicz papers [Mi2, Mi3]).

\noindent
{\bf 1.10. Invariant measures for interval maps}

	We describe here the results from Section 10. 
To investigate the properties of invariant measures it is natural to consider the restriction of $f$ on a component
of the decomposition. We start with studying of $f|B$ for a basic set $B$. By Theorems 4.1 and 8.7 we may apply the results
of [Si1-Si2],[DGS] where a lot of generic properties of maps with the specification are established. 
To formulate the theorem which summarizes the results from [Si1-Si2], [DGS] we need some definitions. 

	Let $T:X \rightarrow X$ be a map of a compact metric space $(X,d)$ into itself. By $M(X)$\index{$M(X)$} we denote the set of all
Borel normalized measures on $X$ (i.e. $\mu (X)=1$ for any $\mu \in M(X)$) with {\it weak topology} (see, e.g., [DGS] for the 
definition). 
If $\mu =T_\ast\mu$ then $\mu$ is said to be {\it invariant}\index{invariant measure}. 
The set of all $T$-invariant measures $\mu \in M(X)$
with the weak topology is denoted by $M_T(X) \equiv M_T$\index{$M_T(X) \equiv M_T$}.
A measure $\mu \in M(X)$ is said to be {\it non-atomic}\index{non-atomic measure} if $\mu (x)=0$ for any point $x\in X$. 
{\it The support of $\mu$}\index{support of $\mu$} 
is the minimal closed set $\:S\equiv supp\,\mu\;$\index{$supp\,\mu$} such that $\mu(S)=1$. A measure $\mu \in M_T$
whose $supp\,\mu$ coincides with one closed periodic orbit is said to be {\it a CO-measure}\index{CO-measure} ([DGS], Section 21);
if $a\in Per\,T$ then the corresponding CO-measure is denoted by $\nu(a)$\index{$\nu(a)$}. 
The set of all CO-measures which are concentrated
on cycles with minimal period $p$ is denoted by $P_T(p)$\index{$P_T(p)$}.

	For any $\:x\in X\;$ let $\:\delta_x \in M(X)\;$\index{$\delta_x$} be the corresponding 
$\delta$-measure\index{$\delta$-measure} (i.e. $\delta_x(x)=1$).
Let $V_T(x)$\index{$V_T(x)$} be the set of limit measures of time averages $N^{-1}\cdot \sum^{N-1}_{j=0} T^j_\ast\delta_x$; 
it is well-known
([DGS], Section 3) that $V_T(x)$ is a non-empty closed and connected subset of $M_T$. A point $x\in X$ is said to have
{\it maximal oscillation}\index{maximal oscillation} 
if $V_T(x)=M_T$. If $V_T(x)=\{\mu\}$ then a point $x$ is said to be {\it generic}\index{generic point for a measure} for $\mu$.

	A measure $\mu \in M_T$ is said to be {\it strongly mixing}\index{strongly mixing measure} if 
$\lim_{n\rightarrow\infty}\mu(A\cap T^{-n}B)=\linebreak\mu(A)\cdot\mu(B)$ for all measurable sets $A,B$. A measure $\mu$ is said to be
{\it ergodic}\index{ergodic measure} if there is no set $B$ such that $TB=B=T^{-1}B,0<\mu(B)<1$.

	We summarize some of the results from [Si1-Si2] and [DGS] in the following
 
\noindent
{\bf Theorem DGS[Si1-Si2],[DGS].} 
{\it Let $T:X \rightarrow X$ be a continuous map of a compact metric space $X$ into itself with the specification property.
Then the following statements are true.

	1) For any $l \in \Bbb N$ the set $\bigcup_{p\ge l}P_T(p)$ is dense in $M_T$.

	2) The set of ergodic non-atomic invariant measures $\mu$ with $supp\,\mu=X$ is residual in $M_T$.

	3) The set of all invariant measures which are not strongly mixing is a residual subset of $M_T$.

	4) Let $V\subset M_T$ be a non-empty closed connected set. Then the set of all points $x$ such that $V_T(x)=V$ is
dense in $X$ (in particular, every measure $\mu\in M_T$ has generic points).

	5) The set of points with maximal oscillation is residual in $X$.}

	Let us return to interval maps. In fact an interval map on its basic set does not necessarily have the specification property.
However, 
applying Theorem DGS and some of the preceding results
(Theorem 4.1, Lemma 8.3 and Theorem 8.7) we prove in Theorem 10.3 that a restriction of an interval map on its basic set
has all the properties 1)-5) stated in Theorem DGS. The fact that statement 5) of Theorem DGS holds for mixing maps 
(since mixing maps has the specification by Theorem 8.7) allows us to specify the third possible type of generic behavior of an orbit  
for maps without wandering intervals (as it was explained after the formulation of Theorem 6.2).  
Moreover, we prove also the following Theorem 10.4 and Corollary 10.5.
 
\noindent
{\bf Theorem 10.4.} 
{\it Let $\mu$ be an invariant measure. Then the following properties of $\mu$ are equivalent.

	1) There exists $x\in [0,1]$ such that $supp\,\mu\subset\omega(x)$.

	2) The measure $\mu$ has a generic point.

	3) The measure $\mu$ can be approximated by CO-measures.}

\noindent
{\bf Remark.} In fact one can deduce Theorem 10.4 for a non-atomic invariant measure directly from Theorem DGS and the above
mentioned Theorem 4.1, Lemma 8.3, Theorem 8.7; this version of Theorem 10.4 was obtained in [Bl4],[Bl7].
Note that even this preliminary version implies the following 

\noindent
{\bf Corollary 10.5[Bl4],[Bl7].}
{\it CO-measures are dense in all ergodic invariant measures of an interval map.}

%

	In what follows we need the definition of a piecewise-monotone continuous map. A continuous map $f:[0,1] \rightarrow [0,1]$
is said to be {\it piecewise-monotone (a pm-map)}\index{piecewise-monotone map, pm-map} 
if there exist $n\ge 0$ and points $0=c_0<c_1<\ldots<c_{n+1}=1$ such that
for any $0\le k\le n\;\; f$ is monotone on $[c_k,c_{k+1}]$ and $f[c_k,c_{k+1}]$ is not a point (by monotone we mean non-strictly  
monotone).  Each $c_k,\,1\le k\le n$ is called 
{\it a turning point of $f$}\index{turning point}; we 
denote the set $\{c_k\}^n_{k=1}$ by $C(f)$\index{C(f)}. 

%

	Let $T:X \rightarrow X$ be a map of a compact metric space $(X,d)$ into itself, 
$\mu \in M_T$ be a $T$-invariant measure. Denote by $h_\mu(T)$ {\it the measure-theoretic entropy of $T$ (with respect to $\mu$)} 
[K]; Theorem DG1 
proved in [Di] and [Go]  
plays an important role in the theory of dynamical systems.

\noindent
{\bf Theorem DG1.}
$\displaystyle h(T)=\sup_{\mu \in M_T}h_\mu(T)=\sup\{h_\mu(T):\mu\in M_T\: ${\it is ergodic} $\}$.
%

	Now let us return to one-dimensional maps. In his recent paper [H3] F.Hofbauer has proved statements 1)-3) of 
Theorem 10.3 for pm-maps. He used the technique which seems to be essentially piecewise-monotone. Moreover, he has proved
that for a pm-map $f$ the set of all $f$-invariant measures $\mu$ such that $h_\mu(f)=0$ is a residual subset of $M_{f|B}$.
The last result can be deduced also from Theorem 10.3 and the following theorem of Misiurewicz and Szlenk.

\noindent
{\bf Theorem MiS[MiS].}
{\it If $f:[0,1] \rightarrow [0,1]$ is a pm-map then the entropy function $h:M_f \rightarrow \Bbb R$ defined by
$h(\mu)=h_\mu(f)$ is upper-semicontinuous}.

	Indeed, by Theorem 10.3.1) and Theorem MiS the set $h^{-1}(0)\cap M_{f|B}$ is a dense $G_\delta$-subset of $M_{f|B}$.
However, the corresponding problem for an arbitrary continuous map of the interval (not necessarily a pm-map) has not been solved yet.
By Theorem 4.1 and Lemma 8.3 it is sufficient to consider a mixing map of the interval. Thus the natural question is whether entropy
$h:M_f\rightarrow \Bbb R$ is upper-semicontinuous provided that $f:[0,1]\rightarrow [0,1]$ is mixing.

	Suppose that the answer is affirmative. Then by Theorem DG1 for any mixing $f$ there exists a measure $\mu$ of maximal entropy.
However, [GZ] contains an example of a mixing map without any such measure. So we get to the following

\noindent
{\bf Problem.}
{\it Let $g:[0,1]\rightarrow [0,1]$ be a mixing map. 

	1) Do measures $\mu\in M_g$ with zero entropy form a residual subset of $M_g$?

	2) What are the conditions on $g$ that would imply the upper-semicontinuity of the entropy $h:M_g\rightarrow \Bbb R$
or at least the existence of a measure of maximal entropy for $g$?}

\noindent
{\bf Conjecture}(cf. Lemma 8.3).
{\it Let $g:[0,1]\rightarrow [0,1]$ be a mixing map and let for any open $U$ there exists an integer $n\in \Bbb N$ such that
$g^nU=[0,1]$. Then the entropy $h:M_g\rightarrow \Bbb R$ is upper-semicontinuous and $g$ possesses a measure of maximal entropy.}

\noindent
{\bf 1.11. The decomposition for piecewise-monotone maps}

	For pm-maps we can make our results more precise. In fact we are going to illustrate how our technique and the above
formulated results work applying them to pm-maps. It will be shown that the Decomposition Theorem for pm-maps is an easy
consequence of the Decomposition Theorem for arbitrary continuous maps and the properties of basic and solenoidal sets.
The properties of pm-maps will not be discussed in Sections 2-11.

	Let $f$ be a pm-map with the set of turning points $C(f)=\{c_i\}^m_{i=1}$. 

\noindent
{\bf Lemma PM1.}
{\it Let $A=\overline {\bigcup_{c\in C(f)} orb\,c}, t'=\inf A,t''=\sup A$. Then $f[t',t'']\subset [t',t'']$.}

\noindent
{\bf Proof.} Left to the reader. $\Box$

\noindent
{\bf Lemma PM2.}
{\it Let $\{I_j\}^\infty_{j=1}$ be generating intervals with periods $\{m_j\}^\infty_{j=1}$. Then in the notation of Theorem 3.1 
$S=S_p=S_\omega$ (i.e. for any $x\in Q=\cap orb\,I_j$ we have $\omega(x)=S=Q\cap \overline {Per\,f}=S_p=Q\cap\omega(f)=S_\omega)$.}

\noindent
{\bf Proof.} By Theorem 3.1 $S\subset S_p=Q\cap \overline{Per\,f}\subset S_\omega=Q\cap \omega(f)$. Suppose that 
$S_\omega\setminus S \neq\emptyset$ and $x\in S_\omega\setminus S$. Then we can make the following assumptions.

	1) Replacing if necessary $x$ by an appropriate preimage of $x$ we can assume that $x\not \in \{orb\,c:c\in C(f)\}$.
Indeed, the fact that $x\in S_\omega$ implies that there exists $y$ such that $x\in \omega(y)$. But $f|\omega(y)$ is surjective,
so for any $n>0$ there exists $x_{-n}\in \omega(y)$ such that $f^nx_{-n}=x$. Moreover, $x\not \in S$, so $x_{-n}\not \in S$
too.

	Now suppose that there are some points $c\in C(f)$ such that for some $m=m(c)$ we have $f^mc=x$. Clearly
$x\not \in Per\,f$; thus the number $m(c)$ is well defined. Take the maximum $M$ of the numbers $m(c)$ over all $c\in C(f)$
which are preimages of $x$ under iterations of $f$ and then replace $x$ by $x_{M+1}$. Obviously $x_{M+1}\not \in \{orb\,c:c\in C(f)\}$.

	2) We can assume $x$ to be an endpoint of non-degenerate component $[x,y], x<y$ of $Q$.

	3) We can assume $i$ to be so large that $orb\,I_i\cap C(f)=Q\cap C(f)\equiv C'$, where $I_i=[x_i,y_i] \ni x$.
Indeed, if $c\in C(f)\setminus Q$ then there exists $j=j(c)$ such that $c\not \in orb\,I_j$. So if $N$ is the maximum
of such $j(c)$ taken over all $c\in C(f)\setminus Q$ then for any $i>N$ we have $orb\,I_i\cap C(f)=Q\cap C(f)$.

	4) By Theorem 3.1 $\omega(c')=S$ for any $c'\in C'$. So $x\not \in \omega(c')$ for any $c'\in C'$. Now the first assumption
implies that $x\not \in \{\overline {\cup orb\,c'}:c'\in C'\}=D\supset S$. Thus we can assume $i$ to be so large that 
$[x_i,x+\varepsilon]\cap D=\emptyset$ for some $\varepsilon>0$.

	Let $A=I_i\cap D, t'=\inf A,t''=\sup A$. Note that then by the assumption 4) $x\not \in [t',t'']$. 
The fact that $x\in S_\omega$ implies that there exists a point $z$ such that $x\in \omega(z)$; by Theorem 3.1 $S \subset \omega(z)$.
At the same time $S\subset D$ and so $S\cap I_i \subset D\cap I_i \subset [t',t'']$. 
But by Lemma PM1 the interval $[t',t'']$ is $f^{m_i}$-invariant and so by the properties of solenoidal sets 
the fact that $S\subset \omega(z)$ implies that $\omega(z) \cap I_i \subset [t',t'']$. Thus we see that
$x\in \omega(z) \cap I_i \subset [t',t'']$ which contradicts the assumption 4). $\Box$


\noindent
{\bf Lemma PM3.}
{\it Let $g:[0,1]\rightarrow [0,1]$ be a continuous map, $I$ be an $n$-periodic interval, $B=B(orb\,I,g)$ be a basic set,
$J\subset I$ be another $n$-periodic interval. Furthermore, let $I=L\cup J\cup R$ where $L$ and $R$ are the components of $I\setminus J$.
Then at least one of the functions $g^n|L,g^n|R$ is not monotone. In particular, if $g$ is a pm-map then 
$C(g)\cap orb\,I \supset C(g)\cap orb\,J$ and  $C(g)\cap orb\,I \neq C(g)\cap orb\,J$.}

\noindent
{\bf Remark.} It is easy to give an example of a fixed interval $I$ containing a fixed interval $J$ and such that
a basic set $B=B(orb\,I,g)$ exists. Indeed, consider a mixing pm-map $f$ with a fixed point $a$ and
then to ``glue in'' intervals instead of $a$ and all preimages of $a$ under iterations of $f$. It is quite easy to see 
that this may be done
in such a way that we will get a new map $g$ with the required property; $J$ will be an interval which replaces $a$ itself.

\noindent
{\bf Proof.} We may assume $n=1$ and $J$ to be a complementary interval to $B$. Suppose  that $g|L$ and $g|R$  are monotone. By 
Theorem 4.1 there exists a transitive map $\psi:[0,1]\rightarrow [0,1]$ which is a monotone factor of the map $g$; in other words
$g$ is semiconjugate to $\psi$ by a monotone map $\phi$. By the definition of a basic set $\phi(J)=a$ is a point $a$.  
The monotonicity of $\phi$ and the fact that $g|L$ and $g|R$ are monotone
imply that $\psi|[0,a]$ and $\psi|[a,1]$ are monotone; moreover, $\psi(a)=a$. Clearly, it contradicts 
the transitivity of $\psi$.$\Box$

	We need the following definition: if $B=B(orb\,I,f)$ is a basic set then the period of $I$ is called 
{\it the period of $B$}\index{period of a basic set}
and is denoted by $p(B)$\index{$p(B)$}. 
To investigate the decomposition for a pm-map let us introduce the following ordering in the family of all
basic sets of 
$f:\: B(orb\,I_1,f)=B_1 \succ B(orb\,I_2,f)=B_2$ if and only if $orb\,I_1\supset orb\,I_2$\index{$\succ$-ordering among basic sets}. 
The definition is correct
for a continuous map of the interval. So it is possible to analyze the structure of the decomposition via $\succ$-ordering 
in the continuous case. We do not follow this way to avoid unnecessary complexity.

	For any set $D\subset C(f)$ consider the family $G(D)$ of all basic sets $B(orb\,I,f)$ such that $D=orb\,I\cap C(f)$.
Let us investigate the properties of the family $G(D)$ with the $\succ$-ordering. Fix a subset $D\subset C(f)$ and suppose that
$B_1=B(orb\,I_1,f)\in G(D), B_2=B(orb\,I_2,f)\in G(D)$. Then either $orb\,I_1\supset orb\,I_2$ or $orb\,I_1\subset orb\,I_2$.
Indeed, otherwise let $J=I_1\cap I_2\neq \emptyset$ and let for instance $p(B_1)\le p(B_2)$. It is easy to see that the period of $J$
is equal to $p(B_2)$. Now Lemma PM3 implies $C(f)\cap orb\,I_2\neq C(f)\cap orb\,J$ which is a contradiction.

	Thus we may assume $orb\,I_1\supset orb\,I_2$; by Lemma PM3 it implies that $p(B_1)<p(B_2)$. So if $D\subset C(f)$
and $G(D)$ is infinite then $G(D)=\{B_1\succ B_2\succ \ldots\}$. Moreover, assume that $B_i=B(orb\,I_i,f)$; then 
$Q(D)\equiv \cap orb\,I_i$ is a solenoidal set, $D\subset Q(D)$ and the corresponding group is $H(\{p(B_i)\}^\infty_{i=1})$. Let us
show that there is no basic set $B(orb\,J,f)=B \not \in G(D)$ such that $orb\,J\subset orb\,I_1$.

	Indeed, let $B=B(orb\,J,f)$
be such a basic set. Let $E=orb\,J\cap C(f)$; then $\emptyset \neq E \subset D,E\neq D$. 
At the same time 
it is easy to see that $Q(D)\subset orb\,J$. Indeed, let $z\in E\subset D$. Then $\omega(z)$ is a solenoidal set belonging to $Q(D)$
and to $orb\,J$ as well; in other words, $\omega(z)\subset Q(D)\cap orb\,J$,      
which implies $Q(D)\subset orb\,J$ and contradicts the fact that $E\neq D$. 

	Note that if $G(D)=\{B_1\succ B_2\succ\ldots\}$ then the well-known methods of one-dimensional symbolic dynamics 
easily yield that $f|B_i$ is semiconjugate by a map $\phi$ to a one-sided shift $\sigma:M\rightarrow M$ and 
$1\le card\,\phi^{-1}(\xi)\le 2$ for any $\xi\in M$. Indeed, let $B_i=B(orb\,I_i,f), B_{i+1}=B(orb\,I_{i+1},f),
orb\,I_{i+1} \subset orb\,I_i$ and $\cal R$ be a collection of components of the set $orb\,I_i\setminus orb\,I_{i+1}$. Then for each
interval $J\in \cal R$ a map $f|J$ is monotone and for some finite subset $\cal F$ of $\cal R$ we have $fJ\supset J'$ if
$J'\in \cal F$ and $fJ\cap J'=\emptyset$ if $J'\not \in \cal F$. Construct an oriented graph $X$ with vertices which are elements
of $\cal R$ and oriented edges connecting $J\in \cal R$ with $J'\in \cal R$ if and only if $fJ\supset J'$. 
This graph generates a one-sided shift $\sigma:M\rightarrow M$ in the corresponding topological Markov chain. Let 
$K=\{x:f^nx\in orb\,I_i\setminus orb\,I_{i+1}\}$. Then $f|K$ is monotonically semiconjugate to $\sigma:M\rightarrow M$
({\it monotonically } means here that a preimage of any point is an interval, probably degenerate) and $B$ coincides with $\partial K$;
in other words, to get a set $B$ from $K$ we need to exclude from $K$ interiors of all non-degenerate intervals which are
components of $K$.

	Let us return to the properties of the family of all basic sets. If $D\subset C(f)$ then $G(D)$ is either infinite or
finite. Let $\{D^i\}^k_{i=1}$ be all $D^i$ such that $G(D)$ is infinite and $\{\widetilde {B_r}\}^R_{r=1}$ be all basic sets belonging to
finite sets $G(D)$. The family of all possible sets $D\subset C(f)$ is finite so $R<\infty,k<\infty$ and basic sets from 
$\{G(D^i)\}^k_{i=1}$ together with the collection  $\{\widetilde {B_r}\}^R_{r=1}$ form the family of all basic sets.
Note that  
$D^i\cap D^j=\emptyset (i\neq j)$. Indeed, otherwise $\emptyset\neq D^i\cap D^j \subset Q(D^i)\cap Q(D^j)$; by the
Decomposition Theorem this is only possible if $Q(D^i)=Q(D^j)$. But $D^r=C(f)\cap Q(D^r)\:(r=i,j)$ and thus $D^i=D^j$ which is
a contradiction. Moreover, if $E\subset C(f)$ is such that $G(E)$ is finite and $E\cap D^i\neq\emptyset$ then $E\supset D^i$. 
Indeed, considering points from $E\cap D^i$ it is easy to see that for any $B=B(orb\,J,f)\in G(E)$ we have $Q(D^i)\subset orb\,J$
and hence $E\supset D^i$. 

	Clearly, we have already described all basic and some solenoidal sets via $\succ$-ordering.
However, there may exist generating intervals $\{I_j\}$ with periods $\{m_j\}$ and the corresponding solenoidal set
$Q=\cap orb\,I_j$ such that $Q\cap C(f)=F$ and $F\neq D^i\:(1\le i\le k)$. Then by the Decomposition Theorem 
$F\cap D^i=\emptyset\;(1\le i\le k)$ and $f|orb\,I_N$ has no basic sets for sufficiently large $N$. Applying the analysis of
maps with zero entropy to $f|orb\,I_N$ we finally obtain the Decomposition Theorem for pm-maps.
\pagebreak

\noindent
{\bf Theorem PM4 (Decomposition Theorem for pm-maps).}
{\it Let $f$ be a \mbox{pm-map.} Then there exist an at most countable family of pairs of basic and $\Omega$-basic sets
$\{B_i\subset B'_i\}$ and a family of triples of solenoidal sets 
$\{S^{(\alpha)}\subset S^{(\alpha)}_\Omega\subset Q^{(\alpha)}\}_{\alpha\in A}$ such that:

	1) $\Omega(f)=X_f\cup (\bigcup_\alpha S^{(\alpha)}_\Omega)\cup (\bigcup_i B'_i)$;

	2) $\omega(f)=\overline {Per\,f}=X_f\cup (\bigcup_\alpha S^{(\alpha)})\cup (\bigcup_i B_i)$;

	3) $card\,A\le card\,C(f)$;

	4) $S^{(\alpha)}=S^{(\alpha)}_p=S^{(\alpha)}_\omega$ for any $\alpha\in A$;

	5) intersections in this decomposition are only possible between different basic or $\Omega$-basic sets, intersection
of any three sets is empty and intersection of any two sets is finite;

	6) there exist a finite number of pairwise disjoint subsets $\{D^i\}^k_{i=1},\{F_j\}^l_{j=1}$ of $C(f)$,
a finite collection of basic sets $\{\widetilde B_r\}^R_{r=1}$ and a finite collection of cycles of intervals $\{orb\,K_j\}^l_{j=1}$
with the following properties:

	a) for any $i,1\le i\le k$ the family $G(D^i)$ is an infinite chain $B^i_1\succ B^i_2\succ\ldots$ of basic sets with periods
$p^i_1<p^i_2<\ldots$ and $Q(D^i)=\bigcap_n orb\,I^i_n$ is a solenoidal set with the corresponding group $H(p^i_1,p^i_2,\ldots)$;

	b) $f|B^i_n$ is semiconjugate to a one-sided shift in a topological Markov chain and the semiconjugation is at most
2-to-1;

	c) if $i\neq j,B\in G(D^i),\widehat B\in G(D^j)$ then neither $B\succ \widehat B$ nor $\widehat B\succ B$;

	d) all basic sets of $f$ are $\{\widetilde B_r\}^R_{r=1}\cup \bigcup^k_{i=1}\{B^i_n\}^\infty_{n=1}$;

	e) for any $j,1\le j\le l$ the cycle of intervals $orb\,K_j$ has period $N_j$, there exists a unique solenoidal set
$Q_j\subset orb\,K_j,\;h(f|orb\,K_j)=0,\; orb\,K_j\cap C(f)=F_j\subset Q_j\subset orb\,K_j$ and the group corresponding to
$Q_j$ is $H(N_j,2N_j,4N_j,\ldots)$;

	f) $\{Q(D^i)\}^k_{i=1}\cup \{Q_j\}^l_{j=1}=\{Q^{(\alpha)}\}_{\alpha\in A}$;

	g) there exists a countable set of pairwise disjoint cycles of intervals $\{orb\,L_j\}$
(perhaps some of them are degenerate) such that $C(f)\cap int\,(orb\,L_j)=\emptyset\;(\forall j)$ and $X_f\subset \bigcup_j orb\,L_j$.}

	Remark that we have not included the proofs of statements 3) and 6.g) which are left to the reader.

	Let us make several historical remarks. Jonker and Rand [JR1,JR2] constructed the ``spectral'' decomposition of 
$\Omega(f)$ for a map with a unique turning point ({\it a unimodal map}\index{unimodal map}); 
they used the kneading theory of Milnor and Thurston [MilT].
The unimodal case was studied also in [Str]. The decomposition was extended to an arbitrary pm-map by Nitecki [N2] and 
Preston [P1-P2].

	Our Decomposition Theorem for a pm-map is related to those of Nitecki and Preston. However, we would like to note some differences:
1) we deduce the Decomposition Theorem for a a pm-map from the general Decomposition Theorem for a continuous map of the interval;
2) we investigate the properties of basic sets using the approach which seems to be new.

\noindent
{\bf 1.12. Properties of specific kinds of piecewise-monotone maps}

	To conclude the part of Introduction concerning pm-maps we discuss some specific kinds of pm-maps. First
we need some definitions. A pm-map $f$ is said to be
{\it topologically expanding} or simply {\it expanding}\index{topologically expanding, expanding} 
if there exists $\gamma>1$ such that $\lambda(fI)\ge \gamma\cdot\lambda(I)$
for any interval $I$ provided $f|I$ is monotone (here $\lambda(\cdot)$ is Lebesgue measure on $\Bbb R$). 
Let $g$ be a continuous interval map, $J$ be a non-degenerate interval such that 
$g^n|J$ is monotone $(n\ge 0)$; following Misiurewicz we call $J$ {\it a homterval}\index{homterval}. 
Remark also that one can define the topological
entropy of $f|K$ without assuming $K$ to be an invariant or even compact set [B3]. Now we are able to formulate

\noindent
{\bf Lemma PM5[Bl3].}
{\it The following properties of $f$ are equivalent:

	1) $f$ is topologically conjugate to an expanding map;

	2) if $d<b$ then $f|[d,b]$ is non-degenerate and if $\{c_1,\ldots,c_k\}=C(f)\cap int\,(\overline{Per\,f})$ then
$\bigcup^k_{i=1}(\bigcup_{n\ge 0}f^{-n}c_i)$ is a dense subset of $[0,1]$;

	3) there exists $\delta>0$ such that $h(f|J)\ge \delta$ for any non-degenerate interval $J$;

	4) $f$ has no homtervals and solenoidal sets.}

\noindent
{\bf Proof.} We give here only a sketch of the proof. 
It is based on the Decomposition Theorem and the following important theorem of Milnor and Thurston,
proved in [MilT].

\noindent
{\bf Theorem MT.} {\it Let $f$ be a pm-map with $h(f)>0$. Then there exists an expanding map $g$ with two properties:

	1) $\lambda (g[d,b])=e^{h(f)}\cdot\lambda([d,b])$ for any $d<b$ provided that $g|[d,b]$ is monotone;

	2) $f$ is topologically semiconjugate to $g$ by a monotone map.}

	An expanding map $g$ with the properties from Theorem MT is called {\it a map of a constant slope}\index{map of a constant slope}.

	Now suppose that statement 1) from Lemma PM5 holds for a map $f$. 
Then the properties of solenoidal sets and the definition of a homterval imply that statement 4) holds. 

	Now suppose statement 4) holds. Then by  
the Decomposition Theorem we see that because of the non-existence of solenoidal sets there are only finitely many basic sets.
Besides it follows from the non-existence of solenoidal sets and homtervals that there are no periodic intervals on which $f$
has zero entropy; it implies  
that  all $\succ$-minimal basic sets are cycles of intervals on which the map $f$ is transitive.
 
	Let us show that it implies statement 2) of Lemma PM5. Indeed, the non-existence of homtervals implies that
the map $f$ is non-degenerate on every  non-degenerate interval, so 2a) holds. Now let us prove that 2b) holds too. 
Let $J$ be an interval; consider the
orbit of $J$ under iterations of $f$. The non-existence of homtervals implies that there are numbers
$n<m$ such that $f^nJ\cap f^mJ\neq\emptyset$. 
It is easy to see now that
there is a weak cycle of intervals $I,fI,\ldots,f^{k-1}I,f^kI\subset I$ and a number $n$ such that 
$\bigcup_{i\ge n}f^iJ=\bigcup^{k-1}_{r=0}f^rI$. But by what has already been proved $\bigcup^{k-1}_{r=0}f^rI$ should contain
a cycle of intervals $M$ on which $f$ is transitive. On the other hand by the properties of basic sets 
$M\subset \overline{Per\,f}$ and clearly, there is $c\in C(f)$ such that $c\in int\,(M)$. Finally we see that there exists
$c\in int\,(\overline{Per\,f})$ with preimages in the interval $J$. It proves statement 2) of Lemma PM5.

	Now suppose that statement 2) holds. Let us show that there exist cycles of intervals
on which the map is transitive. Indeed, otherwise every basic or solenoidal set has an empty
interior. Hence by the Decomposition Theorem for pm-maps we see that the set 
$\omega(f)=\overline {Per\,f}=X_f\cup (\bigcup_\alpha S^{(\alpha)})\cup (\bigcup_i B_i)$ may have a non-empty
interior only if $X_f$ has a non-empty interior and interiors of $\omega(f)=\overline {Per\,f}$ and $X_f$ coincide. 
But $X_f\subset Per\,f$; so if $c\in C(f)\cap int(\overline {Per\,f})$ then $c\in C(f)\cap int({Per\,f})$. Now that
$Per\,f=\bigcup_{n\ge 0}\{x: f^nx=x\}$ we see that $c\in  int(\{x: f^nx=x\})$ for some $n$ which contradicts the fact that
$c$ is a turning point of $f$. So  there exist cycles of intervals
on which the map is transitive. 

	Now take all cycles of intervals on which $f$ is transitive. On each cycle there exists
a semiconjugation with a constant slope map (Theorem MT); indeed, the semiconjugation exists because
transitive interval maps have positive entropy (see e.g. Lemma 9.3; this fact also may be easily deduced from the
Decomposition Theorem). Moreover, in fact it must be a conjugation because otherwise ``expanding'' properties of transitive
maps (see Lemma 8.3) imply that the constant slope map in question is degenerate. Using some technical arguments and
statement 2) itself one can now construct a conjugation between the map $f$ and some expanding map, i.e. statement 1) holds.
The equivalence of all these statements and statement 3) may be proved by similar methods. It completes the sketch of the proof of
Lemma 8.3. $\Box$

	For a map with constant slope the Decomposition Theorem may be refined. Namely, in [Bl10] the following theorem is proved.

\noindent
{\bf Theorem PM6[Bl10].} 
{\it Let $f$ be a map of constant slope and $\{B_i\}^N_{i=1}$ be the family of all basic sets of $f$. Then $N\le card\,C(f)$, the
family of limit sets of genus $0$ is finite and there is no solenoidal sets.}

	Let us apply Theorem PM6 and investigate the continuity of topological entropy for pm-maps. Let $M_n$ be the class
of pm-maps $f$ such that $card\,C(f)\le n$. For any $c\in C(f)$ let $q(c,f)$ be the number of basic sets $B=B(orb\,I,f)$ such
that $c\in orb\,I$ if $c\in Per\,f$ or $\infty$ otherwise.

	The machinery of discontinuity of the entropy as a function $h:M_n\rightarrow\Bbb R$ was investigated in [Mi5]
and in different way in [Bl12]; in [Mi\'Sl]
the analogous result was obtained for piecewise-monotone maps with discontinuities. Roughly speaking if $h$ is not continuous
at $f\in M_n$ (where $card\,C(f)=n$) then there exists $c\in C(f)\cap Per\,f$ which can ``blow up'' turning into a periodic
interval $J$ such that for a new map $g\in M_n$ we have $h(g|orb_gJ)>h(f)$. However, this is impossible if $q(c,f)\ge n$. Namely,
the following theorem holds.

\noindent	
{\bf Theorem PM7[Bl12].}
{\it Let $f\in M_n, card\,C(f)=n$ and $q(c,f)\ge n$ for any $c\in C(f)$. Then the entropy function $h:M_n\rightarrow\Bbb R$
is continuous at $f$.}

	As a corollary we obtain in [Bl12] a new proof of the following result of \linebreak M.Misiurewicz [Mi5].

\noindent
{\bf Corollary PM8[Mi5],[Bl12].}
{\it Let $f\in M_1, C(f)=\{c\}$ and either $h(f)=0$ and $c\not \in Per\,f$ or $h(f)>0$. Then the entropy function 
$h:M_1\rightarrow \Bbb R$ is continuous at $f$.}
  
	The most important example of a pm-map is perhaps {\it a smooth map of an interval}\index{smooth map of an interval}, 
by which we mean a $C^\infty$-map $f:[0,1]\rightarrow [0,1]$
with a finite number of non-flat critical points. We denote the set of all smooth maps with $n$ critical points by $Sm_n$\index{$Sm_n$};
$Sm\equiv \cup Sm_n$\index{$Sm$}. Let us define {\it the Schwarzian derivative}\index{Schwarzian derivative} as 
$Sf\equiv {f'''}/{f'}-{3/2}\cdot{(f''/f')}^2$\index{$Sf$}. If for $f\in Sm_n$ we have $Sf<0$ outside the critical points of $f$ 
then we say that $f$ is a map with {\it negative Schwarzian}\index{map with negative Schwarzian}. 
The family of all such $f$ is denoted by $NS_n$\index{$NS_n$};$NS\equiv\bigcup_{n\ge 0}NS_n$\index{$NS$}.

	{\it Does there exist a smooth map with a wandering interval? } Since Denjoy theorem[D] the question has been attracting 
great attention. The main conjecture was that the answer is negative. Let us describe the history of the verification of this
conjecture.

	0) [D] for a circle diffeomorphism $f$ with the irrational rotation number and $\log f'$ of bounded variation; 

	1) [Mi1] for a map $f\in NS_1$ with a $2$-adic solenoid;

	2) [Gu] for a map $f\in NS_1$;

	3) [MSt] for a map $f\in Sm_1$;

	4) [Yo] for a smooth homeomorphism of the circle with a finite number of non-flat critical points;

	5) [L] for a map $f\in NS$ with critical points which are turning points (the principal step towards the polymodal case);

	6) [BL] for a map $f\in Sm$ with critical points which are turning points;

	7) [MMSt] for a map $f\in Sm$.

	Remark also that in [MMSt] the following nice theorem was proved.

\noindent
{\bf Theorem MMS.}
{\it Let $f\in Sm$. Then there exist $N$ and $\xi>0$ such that $|Df^n(p)|\ge 1+\xi$ for any periodic point $p$ of minimal period $n>N$.}

\noindent
{\bf Remark.} G.Hall constructed an example of a $C^\infty$-piecewise-monotone map with finitely many critical points
(among them there are flat critical points) which has a homterval. It shows that $C^\infty$-property alone is not
sufficient for the conjecture in question to be true.  

	Together with Theorem 6.2 and the Decomposition Theorem for pm-maps these results imply the following

\noindent
{\bf Corollary PM9.}
{\it Let $f\in Sm$. Then there exist $k$ cycles of intervals $\{orb\,I_j\}^i_{j=1}, q$ solenoids $\{Q_j\}^q_{j=1}$ 
and $l$ cycles of intervals $\{L_j\}^l_{j=1}$ such that $i+q\le C(f)$ and the following statements are true:

	1) $f|orb\,I_j$ is transitive $(1\le j\le i)$;

	2) $int\,(orb\,L_j)\cap C(f)=\emptyset \,(1\le j\le l)$;

	3) there exists a residual subset $G\subset [0,1]$ such that for $x\in G$ either $\omega(x)\subset orb\,L_j$ is a cycle for
some $1\le j\le l$, or $\omega(x)=Q_j$ for some $1\le j\le q$, or $\omega(x)=orb\,I_j$ and $V_f(x)=M_{f|orb\,I_j}$ for some $1\le j\le i$.}

\noindent
{\bf Remark.} In [Bl1] we describe generic limit sets for pm-maps without intervals having pairwise disjoint forward iterates.

\noindent
{\bf 1.13. Further generalizations}

	Now we would like to discuss possible generalizations of these results. First note that we consider 
a pm-map as a particular case of a continuous map of the interval; at the same time one can consider a continuous map as
a generalization of a pm-map. It is natural to ask whether there are other generalizations and here a pm-map with finite number 
of discontinuities is another important example.

	This class of maps was investigated by F.Hofbauer in his papers [H1-H3] where he constructed and studied the corresponding
``spectral'' decomposition. It is necessary to mention also the paper [HR] where components of Hofbauer's decomposition with zero
entropy are investigated and the paper [W] where topologically generic limit behavior of pm-maps with finite number
of discontinuities is studied.

	However, we are mostly interested in continuous maps; this leads to the generalization of our results to continuous maps
$f:M\rightarrow M$ of a one-dimensional branched manifold (``graph'') into itself. 
It turns out that the ``spectral'' decomposition and the classification
of its components can be generalized for a continuous map of a ``graph'' with slight modifications.

	More precisely, 
let $f:M\rightarrow M$ be a continuous map of a ``graph''. Let $K=\bigcup^n_{i=1}K_i$ be a submanifold with connected
components $K_1,\ldots,K_n$; we call $K$ a {\it cyclical submanifold}\index{cyclical submanifold} 
if $K$ is invariant and $f$ cyclically permutes the components 
$K_1,\ldots,K_n$. A cyclical submanifold $R$ can generate a maximal limit set; the definition is analogous to that for the interval. 
Namely, let $L=\{x\in R:$ for any open neighborhood $U$ of $x$ in $R$ we have $\overline {orb\,U}=R\}$ be an infinite set. There are two
possibilities.

	1) {\it $f|R$ has no cycles}. Then $f|L$ acts essentially as an irrational rotation of the circle. In this case we denote $L$
by $Ci(R,f)$ and call $Ci(R,f)$ {\it a circle-like} set. For instance, if $g:S^1\rightarrow S^1$ is the Denjoy map of the circle
(i.e. the example of the circle homeomorphism 
with a wandering interval) then $R=S^1$ and $Ci(S^1,g)$ is the unique minimal set of $g$.
The existence of a monotone map which semiconjugates $g$ with the irrational rotation is in this case a well-known fact;
moreover, this semiconjugation is at most 2-to-1 on $Ci(S^1,g)$, i.e. essentially $g|Ci(S^1,g)$ is similar to the corresponding
irrational rotation.
Actually monotone semiconjugation exists in general case as well and shows that in general case a map on its circle-like set
is similar to some irrational rotation.

	2) {\it $f|R$ has cycles}. Then we denote $L$ by $B(R,f)$ and call $B(R,f)$ {\it a basic set}. The properties of a basic set
of a map of a ``graph'' are analogous to those of a basic set of a map of the interval.

	The definitions of a solenoidal set and of a limit set of genus 0 are similar to those for the interval. Limit sets of genus 0,
solenoidal sets, circle-like sets and basic sets are the components of the ``spectral'' decomposition for a map of a ``graph''.

	The Decomposition Theorem for a map of a ``graph'' and its several corollaries are proved in [Bl5,Bl9,Bl11]. For
example, the generic properties of invariant measures are analogous to those for a map of the interval (clear modifications
are connected with the existence of circle-like sets). It should be mentioned also that the famous Sharkovskii theorem on the 
co-existence of periods of cycles (Theorem Sh3) was generalized for continuous maps of the circle[Mi4], of the letter $Y$[ALM]
and of any $n$-od[Ba]. There are also some recent results concerning the description of sets of periods of cycles
for continuous maps of an arbitrary finite ``graph'' into itself [Bl13, Bl14] and for continuous maps of an arbitrary finite ``tree''
into itself [Bl15] (here ``tree'' is a finite ``graph'' which does not contain subsets homeomorphic to the circle).

	Almost all the results of this paper are contained in the author's Ph.D. Thesis (Kharkov,1985).

\vspace{.2in}

\noindent
{\bf Acknowledgments.} I would like to express my gratitude to the Institute for Mathematical Sciences in SUNY at Stony Brook 
for the kind hospitality which made possible the preparation of this paper. I am also grateful to M. Lyubich and
J. Milnor for providing useful comments. 
\vspace{.2in}

\noindent
{\large \bf 2. Technical lemmas}
\vspace{.2in}

	From now on we will use all notions introduced in Section 1 without repeating definitions. Also we will only seldom repeat
formulations of those theorems and lemmas which have already been stated in Section 1. Fix a continuous map $f:[0,1]\rightarrow [0,1]$.
We will prove in this Section some elementary preliminary lemmas which nevertheless seem quite important.
Let us start with the following easy

\noindent
{\bf Lemma 2.1.}
{\it 1) Let $U$ be an interval ,$f^mU\cap U\neq\emptyset$ for some $m$. Then there exists a weakly periodic closed interval
$I$ of period $n$ such that $\overline {orb\,U}=\bigcup^{n-1}_{i=0}f^iI=orb\,I$ and $\{orb\,I\setminus orb\,U\}$ is a finite set.

	2) Let $J$ be a weakly $l$-periodic closed interval. Then $L=\bigcap_{i\ge 0}f^{il}J$ is an $l$-periodic interval.}

\noindent
{\bf Proof.} 1) Clearly, $\bigcup^\infty_{i=0}f^{mi+k}U=J_k$ is an interval for $0\le k<m$. Thus the set  
$\overline {orb\,U}=\bigcup^{m-1}_{k=0}\overline {J_k}$ consists of a finite number of its components and 
$card\,(\overline{orb\,U}\setminus orb\,U)<\infty$. Let $I\supset U$ be a component of $\;\overline {orb\,U}$ and $n$ be the minimal
integer such that $f^nI\cap I\neq\emptyset$. Then $f^nI\subset I$ and the first statement is proved.

	2) The proof is left to the reader. $\Box$

	Denote by $L$ the left side and by $R$ the right side of any point $x\in [0,1]$. 
Now if $T=L$ or $T=R$ is a side of $x\in [0,1]$ then denote by $W_T(x)$ a one-sided semi-neighborhood 
of $x$. Let $U=[\alpha,\beta]$ be an interval, $\alpha<\beta,x\in (\alpha,\beta)$. By $Si_U(x)\equiv\{L,R\}$\index{$Si_U(x)$}
we denote the set
of the sides of $x$; also let $Si_U(\alpha)\equiv \{R\},Si_U(\beta)\equiv \{L\}$.
We will consider a pair $(x,T)_U$\index{$(x,T)_U$} where $T\in Si_U(x)$ and call $(x,T)_U$ 
{\it a $U$-pair}\index{$U$-pair, pair in $U$}
or {\it a pair in $U$}. A set of all $U$-pairs is denoted by $\widehat U$\index{$\widehat U$}. 
If $U=[0,1]$ then we write simply $Si(x)$\index{$Si(x)$} , $(x,T)$\index{$(x,T)$} and call $(x,T)$ {\it a pair}\index{pair}.
If $(x,T)$ is a pair in $U$ then we also say that $T$ is {\it a side of $x$ in $U$}\index{side of a point in an interval}.
Finally, if $f|W_T(x)$ is not degenerate for any $W_T(x)$ then we say that 
$f$ {\it is not degenerate on the side $T$ of $x$}\index{map which is not degenerate on a side of a point}.

	Let us define the way $f$ acts on pairs. Namely, say that 
{\it $(y,S)$ belongs to $f(x,T)$}\index{pair belonging to the image of another pair} if $y=fx$ and 
for any $W_T(x)$ there exists $W_S(y)$ such that $fW_T(x)\supset W_S(y)$.

	Let us formulate without proof some properties of a continuous map of the interval.

\noindent
{\bf Property C1.}
{\it Let $U$ be an interval, $x\in fU=V$ and $T\in Si_V(x)$. Then there exists $y\in U$ and $S\in Si_U(y)$ such that $(x,T)\in f(y,S)$.
In particular:

	1) if $x\in int\,V$ then for any side $T\in Si(x)$ there exists $y\in int\,U$ and a side $S$ of $y$ in $U$ such that 
$(x,T)\in f(y,S)$;

	2) if $x$ is an endpoint of $V$ and there exists $y\in int\,U$ such that $fy=x$ then there exists $z\in int\,U$  and
$S\in Si_U(z)$ such that $f(z,S) \ni (x,T)$.}

\noindent
{\bf Property C2.}
{\it Let $f$ be non-degenerate on the side $T$ of $x$. Then $f(x,T)$ is non-empty.}

\noindent
{\bf Property C3.}
{\it If $I,J$ are closed intervals and $I\subset fJ$ then there exists a closed interval $K\subset J$ such that $fK=I$.}

\noindent
{\bf Property C4.}
{\it Let $U$ be an interval, $x\in U$ be a point, $\lambda(U)\ge \varepsilon>0,\, n>0$. Then there exists an interval $V$ such that
$x\in V\subset U,\,\lambda(f^iV)\le \varepsilon\:(0\le \varepsilon \le n)$ and $\lambda(f^jV)=\varepsilon$ for some $j\le n$.}

	Let us consider some examples.

\noindent
{\bf Example 2.1.}
Let $f(x)\equiv x$. Then $f(x,L)=(x,L)$ and $f(x,R)=(x,R)$ for any $x\in [0,1]$.

\noindent
{\bf Example 2.2.}
Let $f(x)=4x(1-x)$; then $f(1/2,L)=f(1/2,R)=(1,L)$.

\noindent
{\bf Example 2.3.}
Let $f$ be continuous and $x$ be a point of local strict maximum of $f$. Then $f(x,L)=f(x,R)=(fx,L)$.

	Let $I$ be a $k$-periodic interval, $M=orb\,I=\bigcup^{k-1}_{i=0}f^iI$. For every $x\in M$ we consider
three sets which are similar to the well-known prolongation set. Let $\cal U$ be either the family $L$ of all left 
semi-neighborhoods of $x$ in $M$ or the family $R$ of all right semi-neighborhoods of $x$ in $M$ or the family $A$ of all
neighborhoods of $x$ in $M$.
For any $W\in \cal U$ and $n\ge 0$ let us consider the invariant closed set $\overline {\bigcup_{i\ge n}f^iW}$.
Set $P^{\cal U}_M(x,f)\equiv P^{\cal U}_M\equiv \bigcap_{W\in {\cal U}}\bigcap_{n\ge 0}(\overline {\bigcup_{i\ge n}f^iW})$
\index{$P^{\cal U}_M(x,f)\equiv P^{\cal U}_M$}. Let us formulate (without proof) some properties of these sets (we will
write $P_M(x)$\index{$P_M(x)$} instead of $P^A_M(x)$ and $P^{\cal U}(x)$\index{$P^{\cal U}(x)$} instead of $P^{\cal U}_{[0,1]}$).

\noindent
{\bf Property P1.}
{\it $P^{\cal U}_M(x)$ is an invariant closed set and $P_M(x)=P^L_M(x)\cup P^R_M(x)$.}

\noindent
{\bf Property P2.}
{\it Let $y\in \overline {orb\,x}$. Then $P_M(x)\subset P_M(y)$.}

\noindent
{\bf Property P3.}
{\it If $y=f^nx$ and $f^n(x,T)=\{(y,S_i)\}^t_{i=1}$ then $P^T_M(x)=\bigcup^t_{i=1}P^{S_i}_M(y)$.}

	We say that a point $y$ is {\it a limit point of $orb\,x$ from the side $T$}\index{limit point of an orbit 
from the side} or that {\it a side $T$ is a limit side of $y\in \omega(x)$}\index{limit side of a point} 
if for any open semi-neighborhood $W_T(y)$ we have $W_T(y)\cap orb\,x\neq\emptyset$. 

\noindent
{\bf Property P4.}
{\it If $y$ is a limit point of $orb\,x$ from the side $T$ then $P^T_M(y)\supset P_M(x)$ and $P^T_M(y)\supset \omega (x)$.} 

\noindent
{\bf Property P5.}
{\it $f|P^{\cal U}_M(x)$ is surjective.}

\noindent
{\bf Property P6.}
{\it $P^{\cal U}_M(x)=\bigcup^{m-1}_{i=0}f^iP^{\cal U}_M(x,f^m)$.}

	Moreover, the following lemma is true (note that by the definition if $W\in \cal U$ 
then either $x$ is an endpoint of $W$ or
$x\in W$).

\noindent
{\bf Lemma 2.2.}
{\it Let $I$ be a periodic interval, $M=orb\,I,x\in M$. Then one of the following possibilities holds for the set $P^{\cal U}_M(x)$. 

	1) There exists an interval $W\in \cal U$ with pairwise disjoint forward iterates 
and $P^{\cal U}_M(x)=\omega(x)$ is a $0$-dimensional set.

	2) There exists a periodic point $p$ such that $P^{\cal U}_M(x)=orb\,p$.

	3) There exists a solenoidal set $Q$ such that $P^{\cal U}_M(x)=Q$.

	4) There exists a periodic interval $J$ such that $P^{\cal U}_M(x)=orb\,J$.

	If additionally $x\in \Omega(f)$ then $x\in P(x)$.}

\noindent
{\bf Proof.}  
The possibility 1) is trivial. Suppose this possibility does not hold. 
Clearly, it means that if $W\in \cal U$ then for some $l<n$ we have $f^lW\cap f^mW\neq\emptyset$. 
By Lemma 2.1 there exists a periodic interval
$J_W$ such that $\bigcap_{k\ge l} (\overline {\bigcup_{i\ge k}f^iW})=orb\,J_W$. Let us choose a family of intervals 
$\{W_m\}$ so that $W_m\in {\cal U},\;W_m\supset W_{m+1}$ and $\lambda(W_m)\rightarrow 0$. Denote $J_{W_m}$ by $J_m$. Then
$orb\,J_m\supset orb\,J_{m+1}\;(\forall m)$ and $P^{\cal U}_M(x)=\bigcap_{m\ge 0}orb\,J_m$. If periods of $J_m$ tend to infinity
then we get to the case 3) of the lemma. Otherwise $orb\,J_m$ tend either to a cycle (the case 2)) or to a cycle of intervals 
(the case 4)). $\Box$

	Let us consider some examples.

\noindent
{\bf Example 2.4. } Let $f:[0,1]\rightarrow [0,1]$ be a transitive map. Then for any pair $(x,T)$ we have $P^T(x)=[0,1]$.

\noindent
{\bf Example 2.5.} Let $f:[0,1]\rightarrow [0,1], f(0)=0,f(1)=1$ and $fx>x$ for any $x\in (0,1)$. Then for the pair $(0,R)$
we have $P^R(0)=[0,1]$ and for any other pair $(x,T)$ we have $P^T(x)=\{1\}$.
\vspace{.2in}

\noindent
{\large \bf 3.Solenoidal sets}
\vspace{.2in}

	The main theorem concerning solenoidal sets is Theorem 3.1; its formulation may be found in subsection 1.3 of Section 1.

\noindent
{\bf Proof of Theorem 3.1.} If $y\in Q$ then there exists a well-defined element {\bf r}$=(r_0,r_1,\ldots)\in H(D)$ such that 
$y\in f^{r_i}I_i\;(\forall i)$. Let us define $\phi:Q\rightarrow H(D)$ as follows: $\phi(y)\equiv {\bf r}(y)$. Then $\phi$ is continuous,
surjective and $\phi^{-1}({\bf s})=\bigcap_{i\ge 0}f^{s_i}I_i$ is a component of $Q$ for any ${\bf s}=(s_0,s_1,\ldots)\in H(D)$.
Clearly, $\tau\circ\phi=\phi\circ f$ and all the components of $Q$ are wandering.

	Now we are going to prove statement 2). Let us denote by $J_z$ the component of $Q$ containing $z$. 
Besides let  $S$ be the set of all limit points of $S_\Omega$ and also $x\in Q$.
We show that $\omega(x)=S$. First observe that $J_x\cap S_\Omega\neq \emptyset$;
this easily implies that $\omega(x)\subset S$.

	On the other hand let $y\in S$. By the definition there exists a sequence $\{U_i\}$ of intervals, where every $U_i$
is a component of $orb\,I_i$, with the following property: 
$U_i\rightarrow y,\; y\not \in U_i\;(\forall i)$.
Since $U_i\cap Per\,f \neq \emptyset$ we have $y\in \overline {Per\,f}$. Moreover, we can choose a sequence $\{n_i\}$ such that
$f^{n_i}J_x\subset U_i\;(\forall i)$. Therefore $y\in \omega(x)$ and $\omega(x)=S \subset \overline {Per\,f}$. Statement 2)
is proved.

	Statements 3) and 6) easily follow from what has been proved and are left to the reader 
(statement 3) follows from the construction and statement 6) may be deduced from statement 3) and the well-known properties
of the topological entropy). Statement 4) follows from statements 1)-2) and the observation that $J_z$
is wandering for any $z\in Q$ (indeed, $\phi(J_z)$ as a point of $H(D)$ has an infinite
$\tau$-orbit and infinite $\omega$-limit set which together with $\tau\circ\phi=\phi\circ f$ implies that $J_z$ itself is a  
wandering interval). Statement 5) follows from statements 2) and 4).$\Box$

	In the sequel it is convenient to use the following

\noindent
{\bf Corollary 3.2.}
{\it Let $\{I_j\}$ be a family of generating intervals, $Q=\bigcup_{j\ge 0}orb\,I_j$. Then the following statements hold:

	1) $Q\cap Per\,f=\emptyset$;

	2) if $J\subset int\,Q$ is an interval then $J$ is wandering;

	3) if $int\,Q=\emptyset$ (i.e. $Q$ is a solenoid) then $f|Q$ is conjugate to the minimal translation $\tau$ in $H(D)$.}

\noindent
{\bf Proof.} Left to the reader. $\Box$
\vspace{.2in}

\noindent
{\large \bf 4. Basic sets}
\vspace{.2in}

	Now we pass to the properties of basic sets. The main role here plays Theorem 4.1 (see subsection 1.4 of Section 1 
for the formulation).  
Before we prove it let us formulate some assertions which easily follow from Theorems 3.1 and 4.1 and 
show the connection between basic sets and sets of genus 1 and 2 introduced by Sharkovskii in [Sh3-Sh6].

\noindent
{\bf Assertion 4.2[Bl4,Bl7].}
{\it 1) Limit sets of genus $1$ are solenoidal sets which are maximal among all limit sets, and vice versa;

	2) limit sets of genus $2$ are basic sets, and vice versa.}

\noindent
{\bf Assertion 4.3[Bl4,Bl7].} {\it Two following properties of a set $\omega(x)$ are equivalent:

	1) for any $y$ the inclusion $\omega(y)\supset \omega(x)$ implies that $h(f|\omega(y))=0$;

	2) $\omega(x)$ is either a solenoidal set or a set of genus $0$.}

	Now we pass to the proof of Theorem 4.1.

\noindent
{\bf Proof of Theorem 4.1.} We divide the proof by steps. The proofs of the first three ones are left to the reader. 

\noindent
{\bf Step B1.} {\it $f|M$ is surjective.} 

\noindent
{\bf Step B2.} {\it $B$ is an invariant closed set.} 

\noindent
{\bf Step B3.} {\it $B(M,f)=B(orb\,I,f)=\bigcup_{i=0}^{n-1}B(f^iI,f^n)$.} 

\noindent
{\bf Example. }
Let $f:[0,1]\rightarrow [0,1]$ be a transitive map. Then $B([0,1],f)=[0,1]$.

\noindent
{\bf Remark.} One can make the Steps B1-B3 without the assumption $card\,B=\infty$.

	In the rest of the proof we assume $I=M=[0,1]$.

\noindent
{\bf Step B4.} 
{\it For any $x\in B$ there exists a side $T$ of $x$ such that $P^T(x)=[0,1]$ } (we call such $T$ 
{\it a source side}\index{source side}).

\noindent
{\bf Remark.} In general case if $I$ is an $n$-periodic interval, $M=orb\,I, x\in I$ and $T$ is a side of $x$ in $I$ such that
$P^T_M(x)=M$ then we call $T$ {\it a source side of $x$ for $F|M$}.

	Suppose that for some $x\in B$ there is no such side. Then $x\neq 0,1$ 
(indeed, if, say, $x=0$ then the fact that $x\in B$ implies that 
$P^R(x)=[0,1]$ which proves Step B4). Furthermore, the assumption implies that $P^L(x)\neq [0,1]$ and 
$P^R(x)\neq [0,1]$. On the other hand $x\in B$, i.e. by the definition $P(x)=P^L(x)\cup P^R(x)=[0,1]$
(the fact that  $P(x)=P^L(x)\cup P^R(x)$ follows from Property P1 in Section 2). By Lemma 2.2 it implies that $P^L(x)$ and $P^R(x)$ are
cycles of intervals or orbits of periodic points. But the set $B$ is infinite; hence there exist a point $y\in B$ and a side $S$ such
that $y\in int\,P^S(x)$ and so necessarily $P^S(x)=[0,1]$ which is a contradiction.

\noindent
{\bf Step B5.} {\it Let $U$ be an interval and $x\in B\cap int\,(fU)$. Then there exists $y\in (int\,U)\cap B$.}

	Indeed, first let us choose the side $S$ of $x$ in $int\,U$ such that $P^S(x)=[0,1]$ (it is possible by Step B4 and because
$int\,(fU)$ is open). Then by Property C1.1) from Section 2 we can find a point $y\in int\,(U)$ such that $fy=x$ and, moreover, 
$(x,S)\in f(y,T)$ for some side $T$ of $y$ in $U$. Now by the definition of a basic set we see that $y\in int\,(U)\cap B$.

	Let us denote by $\cal B$ the set of all maximal intervals complementary to $B$. 

\noindent
{\bf Step B6.}
{\it If $U\in \cal B$ then $(int\,fU)\cap B=\emptyset$ and either $\overline U$ has pairwise disjoint
forward iterates or for some $m,n$ we have $f^{m+n} \overline U\subset f^m \overline U$.}

	Follows from Step B5.

\noindent
{\bf Step B7.} {\it Let $x\in B$ and $T$ be a source side of $x$. Then for any $V_T(x)$ we have $(int\,V_T(x))\cap B\neq\emptyset$ 
(and so $B$ is a perfect set).}

	Suppose that there exists $V_T(x)$ such that $(int\,V_T(x))\cap B=\emptyset$. We may assume that $V_T(x)\in \cal B$.
By Step B6 and the definition of a source side it is easy to see that $f^n \overline {V_T(x)}\subset \overline {V_T(x)}$ for some $n$ and
$\bigcup_{i=0}^{n-1}f^i \overline {V_T(x)}=[0,1]$. But $B$ is infinite which implies that  $(int\,f^iV_T(x))\cap B\neq\emptyset$
for some $i$. Clearly, it contradicts Step B6.

\noindent
{\bf Step B8.}
{\it Let $\phi:[0,1]\rightarrow [0,1]$ be the standard continuous monotone increasing surjective map such that for any interval 
$U$ the set $\phi(U)$ is degenerate if and only if $(int\,U)\cap B=\emptyset$. Then $\phi$ almost conjugates $f|B$ to a transitive
continuous map $g:[0,1]\rightarrow [0,1]$.}

	The existence of the needed map $\phi$ is a well-known fact. Moreover, by Steps B6 and B7 one can easily see that there exists
the continuous map $g$ with $g\circ \phi=\phi\circ f$. Now let us take any open interval $W\subset [0,1]$ and prove that its $g$-orbit
is dense in $[0,1]$. Indeed, by the construction $\phi^{-1}W$ is an open interval containing points from $B$, so the $f$-orbit
of $\phi^{-1}W$ is dense in $[0,1]$ which implies that $g$-orbit of $W$ is dense in $[0,1]$ as well. So $g$-orbit of any open
set is dense and $g$ is transitive.

\noindent
{\bf Step B9.} {\it $f|B$ is transitive.}

	Follows from Step B8.

	Statements a)-c) of Theorem 4.1 are proved. Statements d)-f) follow from the lemmas which will be proved later. 
Namely in Lemma 9.3 we will prove that $h(g)\ge {\ln 2}/2$ provided $g:[0,1]\rightarrow [0,1]$ is transitive. 
Clearly, it implies statement d). In Lemma 8.3 we establish the connection between transitive and mixing maps 
of the interval into itself and show that $\overline {Per\,g}=[0,1]$ provided $g:[0,1]\rightarrow [0,1]$ is transitive;
statements e) and f) will follow from Lemma 8.3. These remarks complete the proof of the theorem.$\Box$

\noindent
{\bf Corollary 4.4.}
{\it Let $B$ be a basic set. Then $B$ is either a cycle of intervals or a Cantor set.}

\noindent
{\bf Proof.} Left to the reader.$\Box$

	Now we may construct the ``spectral'' decomposition for the sets $\overline {Per\,f}$ and $\omega(f)$. However
to extend the decomposition to the set $\Omega(f)$ we need some more facts. Let $I$ be a $k$-periodic interval,
$M=orb\,I$. Set $E(M,f)\equiv \{x\in M$: there exists a side $T$ of $x$ in $M$ such that $P^T_M(x)=M\}$ (in the case of a basic set 
we call such side {\it a source side}). By Theorem 4.1 if there exists the set $B=B(orb\,I,f)$ then $E(M,f)=B$.
In particular, if $card\,E(M,f)=\infty$ then $E(M,f)=B(M,f)$. The other possibilities are described in the following

\noindent
{\bf Lemma 4.5.} {\it Let $N=[a,b]$ be an $s$-periodic interval, $M=orb\,N,\;E=E(M,f)$ is finite and non-empty. Then $E=orb\,x$
is a cycle of period $k$, $M\setminus E$ is an invariant set and one of the following possibilities holds:

	1) $k=s,f^s[a,x]=[x,b],f^s[x,b]=[a,x]$;

	2) $k=s$ and either $x=a$ or $x=b$;

	3) $k=2s$ and we may assume $x=a,f^s=b$.}

\noindent
{\bf Remark.} Note that by Theorem 4.1 and Lemma 4.5 $E(M,f)\subset \overline {Per\,f}$.

\noindent
{\bf Proof.} Let us assume $N=M=[0,1]$. Clearly, $E$ is closed and $f$ is surjective. Let $\cal B$ be the family of all intervals
complementary to $B$. As in Steps B5-B6 of the proof of Theorem 4.1 we have that

\noindent
{\bf (E1)} {\it for any  $U\in \cal B$ there exists $V\in \cal B$ such that $f \overline U \subset \overline V$.}

	Surjectivity of $f$ implies that

\noindent
{\bf (E2)} {\it $\cal B$ consists of several cycles of intervals; moreover, if $U,V\in \cal B$ and
$f \overline U \subset \overline V$ then $f \overline U=\overline V$.}

	Let us now consider some cases. 

\noindent
{\bf Case 1.}
{\it There are no fixed points $a\in (0,1)$.}

	Clearly, Case 1 corresponds to the possibility 2) of the lemma. 

\noindent
{\bf Case 2.}
{\it There is a fixed point $a\in (0,1)\setminus E$.}

	Let $a\in U=(\alpha,\beta)\in \cal B$. First assume that $U \not \supset (0,1)$. Then 
by E1 we see that $\overline U$ is $f$-invariant and by E2 we see that $\overline {[0,1]\setminus U}$ is $f$-invariant.
Clearly, it implies that neither $\alpha$ nor $\beta$ have a source side which is a contradiction. 

	So we may assume that $U\supset (0,1)$. First suppose that $0\in E$ and there exists $x\in (0,1)$ such that $fx=0$.
Then by Property C1 from Section 2 we see that $(0,1)\cap E \neq \emptyset$ which is a contradiction. The similar statement holds for 1.
We conclude that $E$ is invariant and $M\setminus E$ is invariant.

	It remains to show that the possibility ``$f(0)=0,f(1)=1$'' is excluded (the other possibilities correspond to the possibilities
2) and 3) of Lemma 4.5). Suppose that $E=\{0,1\}, f(0)=0,f(1)=1$. Then for any $b\in (0,1)$ neither $[0,b]$ nor $[b,1]$ are
invariant. Choose $\eta<1$ such that $|x-y|\le 1-\eta$ implies that $|fx-fy|\le \eta$ for any $x,y$. 

	Let us show that if $[c,d]\ne [0,1]$ is invariant then $d-c \le \eta$. Indeed, otherwise $[0,d]$
and $[c,1]$ are invariant which is a contradiction. Thus if $J=[c,d]$ is a maximal by inclusion invariant proper subinterval 
containing the fixed point $a$ then $\lambda (J)\le \eta$. Suppose that $c\neq 0$. Then by the maximality of $J$ for any
$\gamma \in [0,c)$ we get $[0,c]\supset \overline {\bigcup_{i\ge 0}f^i[\gamma,c]}$ 
and hence $[0,1]= \overline {\bigcup_{i\ge 0}f^i[\gamma,c]}$ which contradicts the fact that $c\not \in E$. 

\noindent
{\bf Case 3.}
{\it There is a fixed point $a\in (0,1)\cap E$ and there is no fixed point in $(0,1)\setminus E$.}

	Let $U=(c,a)$ and $V=(d,a)$ be the components of $\cal B$. At least one of them is not invariant (because of the fact that
$a\in E$). By E1-E2 we have $f[c,a]=[a,d]$ and $f[a,d]=[c,a]$; so $c=0$ and $d=1$. But by Case 2 we have $c\not \in E([c,a],f^2)$
and $d\not \in E([a,d],f^2)$. Hence $c,d\not \in E$, i.e. $E=\{a\}$. Now it is easy to see that $M\setminus E$ is invariant which completes the proof.$\Box$

	Now let us describe the properties of $\Omega$-basic sets and the set $\Omega(f)\setminus \omega(f)$ 
(more detailed investigation of the properties of this set one can find in Section 7). To this end we will need 
the results of Coven and Nitecki obtained in [CN]; we summarize them in the following 

\noindent
{\bf Theorem CN.} 
{\it Let $f:[0,1]\rightarrow [0,1]$ be an arbitrary continuous map of the interval $[0,1]$ into itself. Then the following
statements hold:

	1) $\Omega(f)=\Omega(f^n)$ for any odd $n$;

	2) if $x$ has an infinite orbit and $x\in \Omega(f)$ then $x\in \Omega(f^n) \;(\forall n)$;

	3) if $x\in \Omega(f)$ then $x\in \overline {\bigcup_{n>0}f^{-n}x}$;

	4) if $0\in \Omega(f)$ then $0\in \overline {Per\,f}$ and if $1\in \Omega(f)$ then $1\in \overline {Per\,f}$.}

	We will also need Theorem Sh2 which was formulated in subsection 1.1 of Introduction. 

\noindent
{\bf Lemma 4.6.}
{\it Let $x\in \Omega(f)\setminus \omega(f)$. Then there exist a number $m$ and an $m$-periodic interval $I$ such that 
the following statements are true:

	1) $x\in \Omega(f^m)$;

	2) $x$ is one of the endpoints of $I$;

	3) if $x$ does not belong to a solenoidal set then the following additional facts hold: a) $x\in B'(orb\,I,f)$; 
b) $f^kx\in B(orb\,I,f)$ provided $f^kx\in int\,(orb\,I)$; c) $f^{2m}x\in B(orb\,I,f)$.}

\noindent
{\bf Proof.} By Theorem CN.4) and Theorem Sh2 we have $x\neq 0,1$. By Theorem CN.3) we may assume that there exist sequences 
$n_i\nearrow \infty$ and $x_i\nearrow x$  such that $f^{n_i}x_i=x\;(\forall i)$. Finally by Theorem Sh2 we may assume 
that there exists $\eta>0$ such that the interval $(x-\eta,x)$ has pairwise disjoint forward iterates and 
the same is true for the interval $(x,x+\eta)$. 

	Fix $j$ such that $x_j\in (x-\eta,x)$ and consider the set $\overline {\bigcup_{n\ge 1}f^n[x_j,x]}$. Obviously, 
there exists an interval $J=[x,z]$ and an integer $u$ such that $J,fJ,\ldots,f^{u-1}J$ are pairwise disjoint, $f^uJ\subset J$
and $f^lJ\cap [x_j,x)=\emptyset\;(\forall l)$. Moreover, $\bigcap_{r\ge 0}f^{ru}J=N$ is a $u$-periodic interval such that $x\in N$ 
is its endpoint. In other words, we have proved the existence of a periodic interval having $x$ as its endpoint. 

	Remark that $f^k|[x-\delta,x]$ is not degenerate for any $\delta>0, k\in \Bbb N$ (otherwise $x\in Per\,f$). Moreover,
{\it for any $k\in \Bbb N$ and any side $T$ of $f^kx$ such that $T\in f^k(L,x)$ we have $x\in P^T(f^kx)$ and if $x$ does not belong
to a solenoidal set then $P^T(f^kx)$ is a cycle of intervals.} Now if $x$ belongs to a solenoidal set then $orb\,x$ is infinite
and by Theorem CN.2) $x$ belongs to $\Omega(f^n)$ for any $n$. So in case when $x$ 
belongs to a solenoidal set we are done and it remains to consider the case when $x$ does not 
belong to a solenoidal set. 

	Note that if $M=[x,\zeta]$ is a periodic interval then  $x\not \in E(orb\,M,f)$. Indeed, by Lemma 4.5 and Theorem 4.1
(see Remark after the formulation of Lemma 4.5) 
$E(M,f)\subset \overline {Per\,f}$ and at the same time $x\not \in \overline {Per\,f}$ so $x\in E(M,f)$ is impossible. 
Let us assume that
$I=[x,y]$ is the minimal by inclusion periodic interval among all periodic intervals having $x$ as an endpoint. Let $I$ have a period $m$.
Let us consider two possibilities. 

	A) {\it There exists $k\in \Bbb N$ and a side $T$ of $f^kx$ in $f^kI$ such that
$(f^kx,T)\in f^k(x,L)$ }(for example this holds provided that $f^kx\in int\,(f^kI)$).

	Choose the minimal integer $k$ among those existing by the supposition and prove that $f^kx\in E(orb\,I,f)$ and 
$E(orb\,I,f)=B(orb\,I,f)=B$ is infinite. Indeed,
by the minimality of the interval $I$ 
for any semi-neighborhood $V_T(f^kx)$ we easily have that \linebreak 
$\overline {orb\,V_T(f^kx)}=orb\,I$ and so $f^kx\in E(orb\,I,f)$. Now we see that $E(orb\,I,f)$ is not an $f^{-1}$-invariant set;
so by Lemma 4.5 the set  
$E(orb\,I,f)=B(orb\,I,f)=B$ is infinite. So $f^kx\in B$ and by the choice of $k$ we see that $f^vx\not \in int\,(orb\,I)$
for any $0\le v<k$. It proves that $x\in B'(orb\,I,f)$; moreover, we have also proved statement 3b) of Lemma 4.6.

	In the preceding paragraph we have shown that  $\overline {orb\,V_T(f^kx)}=orb\,I$ where $T$ is a 
side of $f^kx$ in $f^kI$ such that $(f^kx,T)\in f^k(x,L)$; clearly, it implies that $x\in \Omega(f^m)$. Furthermore, if $f^mx\in int\,I$
or $f^{2m}x\in int\,I$ then $f^{2m}x\in B$. Otherwise we may assume that $f^mx=f^{2m}x=y$; now the fact that $f^kx\in B$ 
and the choice of $k$ easily imply 
that $y=f^{2m}x\in B$ which completes the consideration of the possibility A).

	B) {\it There are no $k\in \Bbb N$ and side $T$ of $f^kx$ in $f^kI$ such that $T\in f^k(x,L)$}.

	Clearly, we see that   $f^mx=f^{2m}x=y$ and $f^{km}(x,L)=(y,R)$ for any $k\ge 1$. Let us consider the set $P^R(y)$.
By Lemma 2.2 $P^R(y)=orb\,K\ni y$ is a cycle of intervals; we may assume that $y\in K$. Clearly, 
the fact that $f^m(x,L)=(y,R)$ implies that $x\in orb\,K$
and $[x-\eta,x)\cap orb\,K=\emptyset$, so $x$ is an endpoint of $K$. Thus by the choice of $I$ we have $I\subset K$. 
Moreover, it is easy to see that
$I\neq K$ (otherwise the possibility B) is excluded) which implies that $y\in int\,K$. At the same time $y\in E(orb\,K,f)$
by the definition. Repeating now 
the arguments from the previously considered possibility A) we obtain the conclusion.$\Box$
\vspace{.2in}

\noindent
{\large \bf 5. The decomposition.}
\vspace{.2in}

	The aim of this section is to prove the Decomposition Theorem. First let us describe intersections between basic sets, 
solenoidal sets and sets of genus 0. 

\noindent
{\bf Lemma 5.1.}
{\it 1) Let $B_1=B(orb\,I_1,f)$ and $B_2=B(orb\,I_2,f)$ be basic sets, $B'_1$ and $B'_2$ be the corresponding $\Omega$-basic
sets. Let $B_1\neq B_2$ and $B_1\cap B_2\neq\emptyset$. Finally let $A$ be
the union of endpoints of intervals from $orb\,I_1$ and endpoints of intervals from $orb\,I_2$. Then 
$B_1\cap B_2\subset B'_1\cap B'_2\subset A$ and so
$B_1\cap B_2$ and $B'_1\cap B'_2$ are finite. Moreover, if $x\in B'_1\cap B'_2$ then $x$ is not a limit point for both 
$B_1$ and $B_2$ from the same side.

	2) Intersection of any three $\Omega$-basic sets is empty and intersection of any two basic sets is finite.}

\noindent
{\bf Proof.} 1) Obviously it is enough to consider the case when $x\in B_1\cap B_2$. 
It is easy to see that there is no side $T$ of $x$ such that $T$ is a source side for both $f|orb\,I_1$ and 
$f|orb\,I_2$. For the definiteness let $L$ be the only source side of $x$ for $f|orb\,I_1$ and $R$ be the only source side of $x$ for
$f|orb\,I_2$. Let us suppose that $x\in int\,(orb\,I_2)$ and prove that $x$ is an endpoint of one of the intervals 
from $orb\,I_1$. Indeed, otherwise for open $U$ such that $int\,(orb\,I_2) \cap int\,(orb\,I_1)\supset U\ni x$ we have
$\overline {orb\,U}=orb\,I_1=orb\,I_2$ which is a contradiction.

	2) Follows from 1).$\Box$

\noindent
{\bf Example. } Suppose that $g:[0,1]\rightarrow [0,1]$ has the following properties:

	1) $g[0,1/2]=[0,1/2],\:g|[0,1/2]$ is transitive;

	2)  $g[1/2,1]=[1/2,1],\:g|[1/2,1]$ is transitive.

	Then $B_1=[0,1/2]$ and $B_2=[1/2,1]$ are basic sets and $B_1\cap B_2=\{1/2\}$.

\noindent
{\bf Lemma 5.2.} 
{\it The family of all basic sets of $f$ is at most countable.}

\noindent
{\bf Proof.} First consider basic sets $B$ with non-empty interiors. Properties of basic sets
easily imply that these interiors are pairwise disjoint so the family of such sets is at most countable.

	Now let us consider a basic set $B=B(M,f)$ with an empty interior;
then by Corollary 4.4 $\;B$ is a Cantor set. We will show that there exists an interval $W\equiv W(B)$ in $M$ 
complementary to $B$ and such that its forward iterates are disjoint from it and its endpoints 
belong to $B$ and do not coincide with the endpoints 
of intervals from $M$. Indeed, denote by $\cal B$ the family of all complementary to $B$ in $M$
intervals; by Theorem 4.1 they are mapped one into another by the map $f$. Choose
two small intervals $I\in \cal B$ and $J\in \cal B$ belonging to the same interval $K\in M$. If one of them is not periodic
then it has the required properties. Otherwise we may suppose that $f^NI\subset I,\;f^NJ\subset J$ for some $N$; moreover, 
denoting by $L$ the interval lying between $I$ and $J$ we may assume that $L$ is non-degenerate and there are no
intervals from $orb\,I$ or $orb\,J$ in $L$. If for some $n$ we have $f^nL\cap I\neq\emptyset$ then one may
take as the required interval the subinterval of $L$ which is complementary to $B$ in $M$ and
is mapped by $f^n$ in $I$. On the other hand if for any $i$ 
we have $f^iL\cap (I\cup J)=\emptyset$ then we get to the contradiction with the fact that by the
definition of a basic set $\overline {orb\,L}=M\supset K\supset (I\cup J)$. 

	Now suppose that there are two basic sets $B_1\neq B_2$; then it is easy to see that $W(B_1)\cap W(B_2)=\emptyset$.
Indeed, $\overline{W(B_1)}$ and $\overline{W(B_2)}$ have no common endpoints (otherwise by Lemma 5.1 these points are endpoints
of intervals from generating $B_1$ and $B_2$ cycles of intervals which contradicts the choice of $W(B_1)$ and $W(B_2)$). 
On the other hand no enpoints of $\overline{W(B_1)}$ can belong to $W(B_2)$ because of the choice of $W(B_2)$ 
of the fact that endpoints of $\overline{W(B_1)}$ are non-wandering points. Similarly no endpoints of 
$\overline{W(B_2)}$ belong to $W(B_1)$. Hence $W(B_1)\cap W(B_2)=\emptyset$ which implies that the family
of intervals $W(B)$ and the family of all basic sets are at most countable.$\Box$

\noindent
{\bf Lemma 5.3.}
{\it 1) Let $I_0\supset I_1\supset\ldots$ be generating intervals and $Q=\bigcap_{j\ge 0}orb\,I_j$. Then $Q\cap B=\emptyset$ for
any basic set $B$ and if $J_0\supset J_1\supset\ldots$ are generating intervals and $Z=\bigcap_{i\ge 0}orb\,J_i$ then
either $Z\cap Q=\emptyset$ or $Z=Q$.

	2) There is at most countable family of those solenoidal sets $Q=\bigcap_{j\ge 0}orb\,I_j$ which have non-empty interiors.}

\noindent
{\bf Proof.} The proof easily follows from the properties of solenoidal sets (Theorem 3.1) and is left to the reader.$\Box$

	Now we can prove the Decomposition Theorem (Theorem 5.4); the formulation may be found in subsection 1.5 of Section 1. 
Recall that by $X_f$ we denote the union of all limit sets of genus 0 of a map $f$.

\noindent
{\bf Proof of the Decomposition Theorem (Theorem 5.4).} 
We start with statement 2). Let us consider some cases assuming that $x\in \omega(f)$.
If $x\in X_f$ then we have nothing to prove. If $x\in Q$ for some solenoidal set $Q$ then by Theorem 3.1 $x\in  S^{(\alpha)}_\omega$
for the corresponding solenoidal set $S^{(\alpha)}_\omega$. Thus we may assume that 
$x\not \in X_f\cup(\bigcup_\alpha Q^{(\alpha)}_\omega)$. Hence there exists $\omega(z)\ni x$ such that $\omega(z)$ 
is neither a cycle nor a solenoidal set. Clearly, we may assume that $\omega(z)$ is infinite.

	Let us construct a special cycle of intervals $orb\,I$ such that $x\in B(orb\,I,f)$. Recall that 
we say that a point $y$ is {\it a limit point of $orb\,\xi$ from the side $T$} or that 
{\it a side $T$ is a limit side of $y\in \omega(\xi)$}  
if for any open semi-neighborhood $W_T(y)$ we have $W_T(y)\cap orb\,\xi\neq\emptyset$. If $T$ is a limit side of $x\in \omega(z)$ then
by Property P4 $P^T(x)\supset \omega(z)$ and hence $P^T(x)=orb\,I$ is a cycle of intervals. Moreover, the fact that 
$\omega(z)\subset P^T(x)$ is infinite implies that if $\zeta\in \omega(z)$ and $N$ 
is a limit side of $\zeta$ then $N$ is a side of $\zeta$ in $P^T(x)$. 
Thus $P^N(\zeta)\subset P^T(x)$; the converse is also true
and thus $P^N(\zeta)=P^T(x)=orb\,I$ for any $\zeta\in \omega(z)$ and any limit side $N$ of $\zeta$. 
By the definition we have $\omega(z)\subset E(orb\,I,f)=B(orb\,I,f)$ which proves statement 2).

	It remains to note that now statement 1) follows from Lemma 4.6, statement 3) follows from Theorem 3.1 and Theorem 4.1,
statement 4) follows from Theorem 3.1 and Corollary 3.2 and statement 5) follows from Lemma 5.1 and Lemma 5.3. Moreover, the 
family of all basic sets is at most countable by Lemma 5.2. It completes the proof.$\Box$ 

\noindent
{\bf Corollary 5.5.} {\it For an arbitrary $x\in [0,1]$ one of the following  possibilities holds:

	1) $\omega(x)$ is a set of genus 0;

	2) $\omega(x)$ is a solenoidal set;

	3) $\omega(x)\subset orb\,I$ where $orb\,I$ is cycle of intervals and $f|orb\,I$ is transitive;

	4) $\omega(x)\subset B$ for some basic set $B$, $B$ is a Cantor set and if $x$ does not belong to a wandering interval
then $\omega(x)$ is a cycle or $f^nx\in B$ for some $n$.}

\noindent
{\bf Proof.} The proof is left to the reader.$\Box$
\vspace{.2in}

\noindent
{\large \bf 6. Limit behavior for maps without wandering intervals}
\vspace{.2in}

	In this section we describe topologically generic limit sets for maps without wandering intervals. 

We will need the following notions: $Z_f\equiv\{x:\omega(x)$ is a cycle $\}$, $Y_f\equiv int\,Z_f$. 

\noindent
{\bf Lemma 6.1.}
{\it A map $f:[0,1]\rightarrow [0,1]$ has no wandering intervals if and only if the set $Z_f$ is dense.}

\noindent
{\bf Proof.} By the definition if a map $f$ has a wandering interval $J$ then $Z_f$ is not dense because $int\,J\cap Z_f=\emptyset$. 
Now suppose that
$f$ has no wandering intervals and at the same time there is an interval $I$ such that $I\cap Z_f=\emptyset$ (and so  
$Z_f$ is not dense). Let us show that $I$ is a wandering interval. Suppose that 
there exist $n$ and $m$ such that $f^nI\cap f^{n+m}I\neq\emptyset$. Hence the set $\bigcup^\infty_{i=0}f^{n+im}I=K$ is an
interval of some type; moreover, $f^mK\subset K$ and on the other hand $K$ contains no cycle of $f$. It is 
easy to see now that all points from $int\,K$ tend under iterations of $f^m$ to one of the endpoints of
$\overline K$ which is a periodic point of $f$. In other words all points from $int\,K$ belong
to $Z_f$ which is a contradiction. So $I$ has pairwise disjoint forward iterates. But $I\cap Z_f=\emptyset$
and so $I$ is a wandering interval which is a contradiction. This completes the proof. $\Box$

	Note the following property of maps without wandering intervals: all solenoidal sets of 
such maps are in fact solenoids. Now let us pass to the proof of Theorem 6.2; its formulation may be found in subsection 1.6
of Section 1.

\noindent
{\bf Proof of Theorem 6.2.} Let us investigate the set $\Gamma_f=[0,1]\setminus Y_f$. It is easy to see that $\Gamma_f$ has the following
properties:

	1) $\Gamma_f$ is closed and invariant;

	2) $f|K$ is non-degenerate for any interval $K\subset \Gamma_f$;

	3) for any non-degenerate component $I$ of $\Gamma_f$ there exist a non-degenerate component $J$ of $\Gamma$ and
integers $m,n$ such that $J$ is a weakly $m$-periodic interval and $f^nI\subset J$.

	Clearly, property 2) easily implies that if there are two intervals $L, M$, $fL\subset M$ and, moreover, $W$ is a residual 
subset of $M$ then $f^{-1}W\cap L$ is a residual subset of $L$. Thus  
it remains to show that if an interval $J$ is a weakly periodic component of $\Gamma_f$  then Theorem 6.2 holds for
$f|orb\,J$. We may assume that $J=[0,1]$. Then $Y_f=\emptyset$ and $f|K$ is non-degenerate for any open interval $K$. Let $B$
be a nowhere dense basic set. Then $f^{-n}B$ is nowhere dense for any $n$. On the other hand 
by Lemma 5.2 the family of all basic sets is at most
countable. Let $D_f=\{x:$ there is no nowhere dense basic set $B$ such that $f^lx\in B$ for some $l\in {\Bbb N}\}$.
Clearly, it follows from what we have shown that $D_f$ is residual in $[0,1]$ and by Corollary 4.5 for $x\in D_f$ one of 
the following three possibilities holds:

	i) $\omega(x)$ is a cycle; 

	ii) $\omega(x)$ is a solenoid;

	iii) there is a cycle of intervals $orb\,I$ such that $f|orb\,I$ is transitive and $\omega(x)\subset orb\,I$.

	Denote by $\cal T$ the family of all cycles of intervals $orb\,I$ such that $f|orb\,I$ is transitive. Suppose that
there are chosen residual invariant subsets $\Pi_{orb\,I}$ of any cycle of intervals $orb\,I\in \cal T$. 
Now instead of condition iii) let us consider the following condition:

iii$^\ast$) there is a cycle of intervals $orb\,I\in \cal T$ such that $orb\,x$ eventually enters the set $\Pi_{orb\,I}$.

	Then it is easy to show that the set $G_{\Pi}$ of all the points for which one of the conditions i), ii) and iii$^\ast$)
is fulfilled is a residual subset of $[0,1]$. 
Indeed, since $D_f$ is residual in $[0,1]$ we may assume that $D_f=\bigcap^\infty_{i=0} H_i$ 
where $H_i$ is an open dense in $[0,1]$ set for
any $i$. Consider the set $R=\{x: orb\,x$ enters an interior of some cycle of intervals $orb\,I\in \cal T\}$. 
Then $R$ is an open subset of $D_f$. Now set $T_i\equiv int\,(H_i\setminus R)$ and replace every $H_i$ by $H'_i= R\cup T_i$.
Then $T_i$ is an open set and $T_i\cap R=\emptyset$ for any $i$. Moreover, $D'_f=\cap H'_i$ is a residual in $[0,1]$ set. 

	So by the choice of sets $\Pi_{orb\,I}$ and by the
previously mentioned consequence of property 2) we may conclude
that preimages of points from the set $\bigcup_{orb\,I\in \cal T}\Pi_{orb\,I}$ form a residual subset of $R$. 
Clearly, it implies that the set of the points for which one of the conditions  i), ii) and iii$^\ast$)
is fulfilled is a residual subset of $[0,1]$. 

	Now to prove Theorem 6.2 it is enough to observe that one can choose as $\Pi_{orb\,I}$ the set of all points in $orb\,I$
with dense in $orb\,I$ orbit. However, in Section 10 we will show that choosing sets $\Pi_{orb\,I}$ in a different way
one can further specify the limit behavior of those generic points whose orbits are dense in cycles of intervals. 
$\Box$
\vspace{.2in}

\noindent
{\large \bf 7. Topological properties of the sets $Per\,f,\;\omega(f)$ and $\Omega(f)$}
\vspace{.2in}

	In this section we are going mostly to investigate the properties of the set $\Omega(f)\setminus \overline {Per\,f}$. Set
$A(x)\equiv (\bigcup_{n\ge 0}f^{-n}x)\cap\Omega(f)$.

\noindent
{\bf Lemma 7.1.}
{\it 1) If $x\not \in \Omega(f)$ then $A(x)=\emptyset$.

	2) Let $x\in \Omega(f),\:I\ni x$ be a weakly periodic interval and $f^nx\in int\,(orb\,I)$ for some $n$. 
Then $A(x)\subset orb\,I$.

	3) Let $x\in \Omega(f)\setminus \overline {Per\,f},\:I$ be periodic interval such that $x$ is an endpoint of $I$.
Then $A(x)\cap int\,(orb\,I)=\emptyset$.}

\noindent
{\bf Proof.} The proof is left to the reader; note only that statement 3) follows from  Theorem CN.4) 
(Theorem CN was formulated in Section 4).$\Box$

\noindent
{\bf Corollary 7.2.}
{\it   Let $x\in \Omega(f)\setminus \overline {Per\,f},\:I$ be periodic interval such that $x$ is an endpoint of $I$ and 
$f^nx\in int\,(orb\,I)$ for some $n$. Then $A(x)\subset \partial(orb\,I)$.}

\noindent
{\bf Proof.} Follows immediately from Lemma 7.1, statements 2) and 3).$\Box$

\noindent
{\bf Lemma 7.3.}
{\it If $x\in \Omega(f)\setminus \omega(f)$ then there exists a periodic interval $J$ such that $x$ is an endpoint of $J$
and $A(x)\subset \partial(orb\,J)$; if $x$ does not belong to a solenoidal set then we may also assume that
$x\in B'(orb\,J,f)$.}

\noindent
{\bf Proof.} By Lemma 4.6 we may assume that there exists a periodic interval $I=[x,y]$ having $x$ as one of its endpoints. 
Let us consider two possibilities. 

	1) {\it There exists $k$ such that $f^kx\in int\,(orb\,I)$.}

	Then by Corollary 7.2 $A(x)\subset \partial(orb\,I)$ which together with 
Lemma 4.6 proves Lemma 7.3 in this case  
(recall that by Theorem Sh2 $\overline {Per\,f}\subset \omega(f)$ and so 
$\Omega(f)\setminus \omega(f)\subset \Omega(f)\setminus \overline {Per\,f}$). It implies Lemma 7.3 in 
the case when $x$ belongs to a solenoidal set 
(indeed, if $x$ belongs to a solenoidal set then obviously there exists an integer $k$ such that $f^kx\in int\,(orb\,J)$.). 

	2) {\it For any $k$ we have $f^kx\not \in int\,(orb\,I)$.}

	Then by Lemma 4.6 we may assume that $x\in B'(orb\,I,f),\:I=[x,y]$ has  a period $m$ and $f^mx=y=f^my$. 
By Lemma 7.1.3) $A(x)\cap int\,(orb\,I)=\emptyset$. Suppose $A(x)\not \subset \partial (orb\,I)$ and show that 
there exists a periodic interval $J$ such that $x\in B'(orb\,J,f)$ and $A(x)\subset \partial(orb\,J)$.

	Indeed, if $A(x)\not \subset \partial(orb\,I)$ then there exists $z\in A(x)\setminus orb\,I$. By Lemma 4.6 $z\in B'(orb\,J,f)$
for some $n$-periodic interval $J$. We may assume $x\in J$; then $f^{nm}x=y\in f^{nm}J=J$ and thus $I=[x,y]\subset J,\;I\neq J$.
Clearly, $x\in B'(orb\,J,f)\setminus B(orb\,J,f)$ because $z\in B'(orb\,J,f)$ is a preimage of $x$ under the corresponding
iteration of $f$ and at the same time $x\not \in \omega(f)$; so $x$ is an endpoint of $J=[x,\zeta]$. Hence $f^mx=y\in int\,J$ and 
as in case 1) we see that by 
Corollary 7.2 $A(x)\subset \partial(orb\,J)$. This completes the proof.$\Box$

	To formulate the next corollary connected with the results of [Y] and [N] we need some definitions.  
Let $c$ be a local extremum of $f$. It is said to be {\it an o-extremum}\index{o-extremum} 
in the following cases:

	1) $c$ is an endpoint of an interval $[c,b]$ such that i) $f|[c,b]$ is degenerate, ii) $f$ is not degenerate in any
neighborhood of each $c$ and $b$, iii) $c$ and $b$ are either both local minima or both local maxima;

	2) there is no open interval $(c,b)$ such that $f|(c,b)$ is degenerate (note that neither in case 1) nor in case 2)
we assume that $c<b$).

	In [Y] the following theorem was proved.
 
\noindent
{\bf Theorem Y.}
{\it Let $f:[0,1]\rightarrow [0,1]$ be a pm-map and $x\in \Omega(f)\setminus \overline{Per\,f}$. Then there exists $n>0$ and 
a turning point $c$ such that $f^nc=x$.}

	On the other hand [N1] contains the following 

\noindent
{\bf Theorem N.}
{\it If $f$ is a pm-map then $\overline {Per\,f}=\omega(f)$.}

\noindent
{\bf Remark.} Note, that Theorem N may be also deduced from Lemma PM2 (see subsection 1.11 of Introduction) and the Decomposition
Theorem. 

	So the following Corollary 7.4 generalizes Theorem Y.

\noindent
{\bf Corollary 7.4.}
{\it If $x\in \Omega(f)\setminus \omega(f)$ then there exist an o-extremum $c$ and $n>0$ such that $f^nc=x$.}

\noindent
{\bf Proof.} Take the interval $J$ existing for the point $x$ by Lemma 7.3. Then $f|orb\,J$ is a surjective map and at the same time 
$x$ is not a periodic point. 
Hence we may choose the largest $n$ such that there exists an endpoint $y$ of
an interval from $orb\,J$ with the following properties: $f^ny=x$ and $y,fy,\ldots,f^ny$ are endpoints of intervals
from $orb\,J$. Then by the choice of $y$ there exists a point $z$ and an interval $[a,b]$ from $orb\,J$ such that
$z\in (a,b), fz=y,fa\neq y,fb\neq y$. Now it is easy to see that we may assume $z$ to be an o-extremum.$\Box$

\noindent
{\bf Remark.} Corollary 7.4 was also proved in the recent paper [Li] (see Theorem 2 there).

	Theorem 7.5 describes another sort of connection between the sets $\omega(f)$ and $\Omega(f)$.

\noindent
{\bf Theorem 7.5[Bl1],[Bl7].}
$\displaystyle \omega(f)=\bigcap_{n\ge 0}f^n\Omega(f)$.

\noindent
{\bf Proof.} By the properties of limit sets for any $z$ we have $f\omega(z)=\omega(z)$. It implies that  
$\omega(f) =f\omega(f)\subset \bigcap_{n\ge 0}f^n\Omega(f)$. At the same time by Lemma 7.3 the set $A(x)$ is finite for any 
$x\in \Omega(f)\setminus \omega(f)$. So by the definition of $A(x)$ we see that if $x\in \Omega(f)\setminus \omega(f)$
then $x\not \in \bigcap_{n\ge 0}f^n\Omega(f)$ which implies the conclusion.$\Box$

	Finally in Theorem 7.6 we study the structure of the set $\Omega(f)\setminus\overline {Per\,f}$; the formulation of 
Theorem 7.6 may be found in subsection 1.7 of Section 1. 

\noindent
{\bf Proof of Theorem 7.6.} We divide the proof into steps.

\noindent
{\bf Step A.}
{\it $card\{\omega(f)\cap U\}\le 1$ and if $x\in \omega(f)\cap U$ then $x$ belongs to a solenoidal set.}

	If $x\in \omega(f)\cap U$ then by the Decomposition Theorem there exists a solenoidal set $Q\ni x$; so if
$x\in \omega(f)\cap U$ then $card(orb\,x)=\infty$. Now it follows from Theorem 3.1 and Corollary 3.2 that if $J$ is the component
of $Q$ containing $x$ then up to the orientation we may assume that $J=[x,b]$ and, moreover, $J$ is a wandering interval. 

	Suppose that there exists $y\in \omega(f)\cap U,\;y\neq x$. 
Then the fact that $J$ is a wandering interval implies that
$a<y<x$. Moreover, similarly to what we have seen in the previous paragraph it is easy to see now that there exists  
a solenoidal set $\widetilde Q$ such that $K=[a,y]$ is its component. By the
properties of solenoidal sets there exist intervals $M=[\tilde a, \tilde y]$ and $N=[\tilde x, \tilde b]$ such that
$y<\tilde y<\tilde x<x$ and $f^nM=M,f^nN=N$ for some $n$. Clearly, it implies that 
$f^n[\tilde y,\tilde x]\supset [\tilde y,\tilde x]$ and so there exists a point 
$z\in [\tilde y, \tilde x]$ such that $f^nz=z$ which is a contradiction. So $card\{\omega(f)\cap U\}\le 1$.

\noindent
{\bf Step B.}
{\it $\Omega(f)\cap U$ has in $U$ at most one limit point, which necessarily belongs to some solenoidal set $S_\omega$.}

	By Theorem Sh2 limit points of $\Omega(f)$ belong to $\omega(f)$. Thus Step B follows from Step A.

	Let $J$ be a periodic interval and suppose that 
one of the endpoints of $J$ belongs to $U$. Then the endpoint of $U$
belonging to $J$ is uniquely determined; we denote this endpoint of $U$ by $e=e(J)$.

\noindent
{\bf Step C.}
{\it The point $e$ is uniquely defined and does not depend on $J$.}

	Clearly, it is sufficient to show that there is no pair of periodic intervals $I=(a',y)$ and $J=(x,b')$ where 
$x,y\in U, a'<a, b'>b$.
To prove this fact observe that if these intervals existed then the interval $K$ with endpoints $x,y$ would have the property 
$f^nK\supset K$ for some $n$ which is impossible. 

	In the rest of the proof we assume that $e=b$. 

\noindent
{\bf Step D.} {\it If $z\in \Omega(f)\cap U$ and $orb\,z$ is infinite then $[z,b]$ has pairwise disjoint forward iterates.}

	If $z$ belongs to a solenoidal set then Step D is trivial by the properties of solenoidal sets
(see Theorem 3.1 and Corollary 3.2). So
we may assume that there exists a periodic interval $J=[z,c]$ such that $z\in B'(orb\,J,f)=B'$. Hence there exists an interval
$[z,d],\;d\ge b$, which is a complementary to $B(orb\,J,f)=B$ in $orb\,J$ interval. If $[z,d]$ does not have pairwise disjoint iterates 
then clearly, there exists a weakly periodic interval $K$ which is a complementary to $B$ in $orb\,J$ interval and, moreover,
$f^mz\in int(orb\,K)$ for large $m$.  
At the same time by the definition of an $\Omega$-basic set $f^mz\in B$ for a large $m$. Clearly, 
this is a contradiction.

\noindent
{\bf Step E.}
{\it If $x,y\in \Omega(f)\cap U,\,x<y,$ then $card(orb\,x)<\infty$ and $x$ belongs to an $\Omega$-basic set
$B'(orb\,[x,d], f)$ for some periodic interval $[x,d]$.}

	If $card\,(orb\, x)=\infty$ then $[x,b]$ has pairwise disjoint forward iterates which is impossible 
because $y\in \Omega(f)\cap (x,b)$.
Hence $card\,(orb\,x)<\infty$; by the Decomposition Theorem it implies that $x$ belongs to an $\Omega$-basic set 
$B'(orb\,[x,d], f)$ for some periodic interval $[x,d]$.	

\noindent
{\bf Step F.}
{\it Let $x,y\in \Omega(f)\cap U,\;x<y,\;x\in B'(orb\,J,f)$ and $y\in B'(orb\,I,f),$ where $J=[x,c]$ and $I=[y,d]$. Then $d<c$}.

	Suppose that $c\le d$. Then by Step C those iterations of $I$ which do not coincide with $I$ have empty intersections with $U$.
Thus by the definition of a basic set we have $B(orb\,J,f)\cap(y,c)=\emptyset$. Moreover, by the Decomposition Theorem 
$B(orb\,J,f)\subset \overline {Per\,f}$ and so $[x,y]\cap B(orb\,J,f)=\emptyset$.  Hence $B(orb\,J,f)\cap J\subset \{c\}$
which contradicts the definition of a basic set.

\noindent
{\bf Step G.}
{\it The point $a$ is not a limit point of $\Omega(f)\cap U$}.

	Suppose that $a$ is a limit point of $\Omega(f)\cap U$. We may assume that 
$x_{-i}\searrow a$ while $i\rightarrow \infty$ and (by Step E) that $\;card(orb\,x_{-i})<\infty\;(\forall i>0)$. 
By Step F we may assume also that for any $i>0$ there exists an $n_i$-periodic interval $J_i=[x_{-i},d_i]$ such that
$x_{-i}\in B'(orb\,J_i,f)$ and $J_{i+1}\supset J_i\;(\forall i>0)$. Clearly, we may assume that $n_i=1\;(\forall i>0)$.
Indeed, as we have just shown $J_{i+1}\supset J_i$, so periods of $J_i$ decrease and hence become equal to some constant;
we will consider the case when this constant is $1$, the arguments in the general case are similar.

	By the definition and Theorem 4.1 $B(orb\,J_{i+1},f)=B\subset [d_i,d_{i+1}]$ for any $i>0$.
Indeed, basic sets belong to $\overline{Per\,f}$, so $B\cap U=\emptyset$ and $B\subset [b,d_{i+1}]$.
But by the definition of a basic set and the fact that $[x_{-i},d_i]$ is invariant we see that there are no points of
$B$ in $[b,d_i)$ which implies that $B\subset [d_i,d_{i+1}]$.

	Let us choose $i>0$ such that for any $y,z$ we have $|fz-fy|<|d_1-x_{-1}|$ provided $|z-y|<d_{i+1}-d_i$; 
clearly, it is possible because $d_{i+1}-d_i \rightarrow 0$ while $i \rightarrow \infty$. 
We are going to show that the interval $[x_{-i},d_{i+1}]$ is invariant. Indeed, let $z\in [x_{-i},d_{i+1}]$. If in fact
$z\in [x_{-i},d_i]$ then $fz\in [x_{-i},d_i]\subset [x_{-i},d_{i+1}]$. If $z\in [d_i,d_{i+1}]$ then by the choice of $i$
we see that $|fz-f\zeta|<|d_1-x_{-1}|$ for any $\zeta\in [d_i,d_{i+1}]$. Choose any $\zeta\in B\subset [d_i,d_{i+1}]$; then
$f\zeta\in B\subset [d_i,d_{i+1}]$ as well and so $|d_1-x_{-1}|>|fz-f\zeta|>|fz-d_i|$ which implies that $fz\in [x_{-i},d_{i+1}]$.
Hence $[x_{-i},d_{i+1}]$ is invariant which contradicts the definition of a basic set and the existence of the basic set $B$.

	Recall that by Step B the set $\Omega(f)\cap U$ has at most one limit point which we denote by $x$. By Step G $x\neq a$. 
Now if $x=b$ then $\Omega(f)\cap U<b$. If $x\in U$ then by Step B $x$ belongs to some solenoidal set $S_\omega$ and the fact that
$e=b$ (see Step C) implies that $[x,b]$ has pairwise disjoint forward iterates 
and so all non-limit points of $\Omega(f)\cap U$ are less then
$x$. This observation shows that the formulation of Step H is correct.

\noindent
{\bf Step H.}
{\it Let $\Omega(f)\cap U\supset\{x_i\}^\infty_{i=0}$ where $\{x_i\}^\infty_{i=0}$ is the whole set of non-limit points
of $\Omega(f)\cap U$; moreover, let $x_0<x_1<\ldots,\;x_n\rightarrow x$. Then there exist periodic intervals
$J_i=[x_i,d_i],\;J_0\supset J_1\supset\ldots$ such that  $x_i\in B'(orb\,J_i,f)\;(\forall i)$ and $\cap J_i=[x,b]$. Moreover,
if periods of the intervals $J_i$ tend to infinity then $[x,b]$ belongs to a solenoidal set and either $x=b$ and
$\Omega(f)\cap U=\{x_i\}^\infty_{i=0}$ or $x<b$ and $\Omega(f)\cap U=\{x_i\}^\infty_{i=0}\cup\{x\}$. On the other hand,
if periods of $J_i$ do not tend to infinity then $x=b$ and so $\cap J_i=\{b\}$.}

	The existence of the intervals $J_i=[x_i,d_i]$ such that
$x_i\in B'(orb\,J_i,f)$ and $J_i\supset J_{i+1}\;(\forall i)$ follows from Steps E and F.  
If periods of $J_i$ tend
to infinity then the required property follows from the properties of solenoidal sets (Theorem 3.1, Corollary 3.2). Now suppose that
periods of $J_i$ do not tend to infinity; consider the case when all $J_i$ are invariant (i.e. have period $1$), the general
case may be considered in the similar way. 

	We are going to show that $\cap J_i=\{b\}$. Indeed, let $\cap J_i=[b',d'],\;b'<d'$; then clearly
$\lim x_i=b'\le b$. Choose $i$ such that for any $y,z$ we have
$|fz-fy|<d'-b'$ provided $|z-y|<|d_i-d'|$. Now repeating all the arguments from Step G we get the same contradiction. 
Indeed, for any $i$ the set $B(orb\,J_i, f)=B_i$ has an empty intersection with $[x_i, b')$ because 
$B_i\subset \overline {Per\,f}$ by Theorem 4.1 and at the same time there is no points of $\overline {Per\,f}$
in $[x_i,b')$. On the other hand the choice of $i$ and the fact that $B_i$ is invariant imply (as in Step G) that
$[b', d_i]$ is an invariant interval which contradicts the definition of a basic set and the existence of the set $B_i$.
This contradiction shows that $b'=d'=b$ which completes Step H.  

	Now let us consider different cases depending on the properties of the set $\Omega(f)\cap U$. First of all let us note that
the properties of points $x\in \Omega(f)\cap U$ such that $(x,b)\cap \Omega(f)\neq\emptyset$ are fully described in Steps E and F;
together with the definitions it completes the consideration of the case 2) and proves the corresponding statements from the other
cases. Furthermore, by Step B we see that 
$\Omega(f)\cap U$ has at most one limit point in $U$ and if so then by Steps E - H we get the case 3)
and also the first part of the case 4) of Theorem 7.4. The second part of the case 4) follows from Step H. This
completes the proof of Theorem 7.4.$\Box$
\vspace{.2in}

\noindent
{\large \bf 8. Transitive and mixing maps}
\vspace{.2in}

	In this section we will investigate the properties of transitive and mixing interval maps which are closely related 
to the properties of maps on their basic sets as it follows from theorem 4.1. 
Let us start with the following simple

\noindent
{\bf Lemma 8.1[Bl7].} {\it Let $f:[0,1]\rightarrow [0,1]$ be a transitive map, $x\in (0,1)$ be a fixed point,
$\eta>0$. Then there exists $y\in (x,x+\eta)$ such that $f^2y> y$ or $y\in(x-\eta,x)$ such that $f^2y< y$. }

\noindent
{\bf Proof.} First suppose there is a point $z\in (x,x+\eta)$ such that $fz>z$. Then choose the maximal fixed point $\zeta$ 
among fixed points which are smaller than $z$. Clearly, if we take $y>\zeta$ close enough to $\zeta$ we will see that $f^2y>y$.
Moreover, we can similarly consider the case when there is a point $z\in (x-\eta,x)$ such that $fz<z$. So we may assume that
for points from $(x-\eta,x+\eta)$ we have $fz<z$ if $x<z$ and $fz>z$ if $x>z$. 

	Now choose $\delta>0$ such that $\delta<\eta, f[x,x+\delta]\subset (x-\eta,x+\eta)$. The map $f$ is transitive so 
$f[x,x+\delta]=[a,b]$ where $a<x$ and $b\ge x$. Moreover, by the transitivity of $f$ one can easily see that there is a point $d\in [a,x]$
such that $fd>x+\delta$ (otherwise $[a,x+\delta]$ is an invariant interval). 
Take $y\in [x,x+\delta]$ such that $fy=d$; clearly, $y$ is the required point. $\Box$

\noindent
{\bf Lemma 8.2.} {\it Let $f:[0,1]\rightarrow [0,1]$ be a transitive map, $\eta>0$. Then there exist a fixed point
$x\in (0,1)$, a periodic point $y\in (0,1),\:y\neq x$ with minimal period $2$ and an interval $U\subset [x-\eta,x+\eta]$
such that $x\in U\subset fU$.}

\noindent
{\bf Proof.} The existence of a fixed point in $(0,1)$ easily follows from the transitivity of $f$. 
Let us show that there exists a point $y$ of minimal period $2$. We may assume that
$1$ is not a periodic point of minimal period $2$.   
Suppose that $x$ is a fixed point and there exists $\varepsilon>0$ such that
for points from $(x-\varepsilon,x+\varepsilon)$ we have $fz<z$ if $x<z$ and $fz>z$ if $x>z$. 
By Lemma 8.1 there exists, say, $\zeta\in (x,x+\varepsilon)$
such that $f^2\zeta>\zeta$. Now if there are no fixed points in $(x,1]$ then set $\xi=1$; otherwise
let $\xi$ be the nearest to $\zeta$ fixed point which is greater than $\zeta$. 
By the construction for any $\alpha\in (\zeta,\xi)$ we have $f\alpha<\alpha$; it easily implies that if  
a point $\beta\in (\zeta,\xi)$ is sufficiently close to $\xi$ then $f^2\beta<\beta$. Together with $f^2\zeta>\zeta$ it shows that
there is a periodic point $y\in (\zeta,\beta)$ such that $f^2y=y$; at the same time by the choice of $\xi$ we have
$fy\neq y$, so the minimal period of $y$ is $2$. 

	Now suppose that there is no fixed point $x$ for which there exists $\varepsilon>0$ such that
for points from $(x-\varepsilon,x+\varepsilon)$ we have $fz<z$ if $x<z$ and $fz>z$ if $x>z$. Then clearly, there are at least two
fixed points, say, $a$ and $b$, and we may assume that $a<b$ and $z<fz$ for $z\in (a,b)$. Let us show that
$a\in f[b,1]$. Indeed, otherwise $I=[b,1]\cup f[b,1]\neq [0,1]$ is an $f$-invariant interval which contradicts the
transitivity. Choose the smallest $c\in [b,1]$ such that $fc=a$; then $b<c$. It is easy to see that again by the
transitivity there exists $d\in (a,c)$ such that $fd=c$. Choose the fixed point $a'$ in such a way that 
the interval $(a',d)$ does not contain fixed points. Then for $z$ sufficiently close to $a'$ we have $f^2z>z$
which together with the fact that $f^2d=a<d$ implies that there is a periodic point $y\in (z,d)$ of minimal period $2$. 

	The proof of the existence of the interval
$U\subset (x-\eta,x+\eta)$ with $x\in U\subset fU$ uses arguments similar to those from Lemma 8.1. 
Indeed, if there is a point $z\in (x,x+\delta)$ such that $fz>z$ or $z\in (x-\delta,x)$ such that $fz<z$ then it is
sufficient to take $U=(x,z)$. So we may assume  that
for points from $(x-\eta,x+\eta)$ we have $fz<z$ if $x<z$ and $fz>z$ if $x>z$. Now take a point $y\in (x-\eta,x+\eta)$
which exists by Lemma 8.1; we may assume that $y\in (x,x+\delta, f(x,x+\delta)\subset (x-\eta,x+\eta)$ and $f^2y>y$.
Then it is easy to see that $U=[x,y]\cup f[x,y]$ is the required interval. $\Box$

	The following Lemma 8.3 establishes the close connection between mixing and transitive interval maps.  

\noindent
{\bf Lemma 8.3[Bl7].} 
{\it Let $f:[0,1]\rightarrow [0,1]$ be a transitive map. Then one of the following possibilities holds:

	1) the map $f$ is mixing  and, moreover, for any $\eta>0$ and any non-degenerate interval $U$ there exists 
$n_0$ such that $f^nU\supset [\eta,1-\eta]$ for any $n>n_0$;

	2) the map $f$ is not mixing and, moreover, there exists a fixed point $a\in (0,1)$ such that 
$f[0,a]=[a,1],\;f[a,1]=[0,a],\;f^2|[0,a]$ and $f^2|[a,1]$ are mixing.

	In any case $\overline {Per\,f}=[0,1]$.}

\noindent
{\bf Proof.} 1) First suppose there exists a fixed point $x\in (0,1)$ such that $x\in int\,f[0,x]$ or $x\in int\,f[x,1]$.
To be definite suppose that $x\in int\,f[0,x]$ and prove that $f$ is mixing and has all the properties from statement 1).
Clearly, we may assume that $x\in int\,f[b,x]$ for some $0<b<x$. By Lemma 8.2 there exists a closed interval 
$U\subset f[0,x]$ such that $x\in U\subset fU$. Let $V$ be any open interval. By Lemma 2.1 the set 
$[0,1]\setminus \bigcup_{m\ge 0}f^mV$ is finite. On the other hand, the set $\bigcup_{n\ge 0}f^{-n}x$ is infinite. So $x\in f^kV$
for some $k$. Now the transitivity implies that $f^lV\supset [b,x]$ for some $l$ and so $U\subset f^{l+1}V$. 

	At the same time the inclusion $U\subset fU$ and the transitivity imply that for any $\varepsilon>0$ there exists 
$N=N(\varepsilon)$ such that $f^nU\supset [\varepsilon,1-\varepsilon]$ for $n\ge N$. Thus $f^mV\supset [\varepsilon,1-\varepsilon]$
for $m>N+l$. It completes the consideration of the case 1).

	2) Suppose there exists a fixed point $a\in (0,1)$ such that $a\not \in int\,f[0,a]$ and $a\not \in int\,f[a,1]$.
By the transitivity $f[0,a]=[a,1],\;f[a,1]=[0,a]$; moreover, $f^2|[0,a]$ and $f^2|[a,1]$ are transitive and hence by the case 1)
$f^2|[0,a]$ and $f^2|[a,1]$ are mixing. The fact that $\overline {Per\,f}=[0,1]$ easily follows from what we have proved.$\Box$

	In the proof of Theorem 4.1 we announced that statements e) and f) of it would follow from Lemma 8.3. Let us prove
the statements now; for the sake of convenience we will recall their formulations. 

	{\it e) If $B=B(orb\,I,f)$ is a basic set then $B\subset \overline {Per\,f}$.}

\noindent
{\bf Proof.} Clearly, it is enough to consider the case when the period of the interval $I$ is $1$. In this case by the preceding 
statements of Theorem 4.1 $f|B$ is almost conjugate by a monotone map $\phi: I\rightarrow [0,1]$ to a transitive map 
$g:[0,1]\rightarrow [0,1]$. By Lemma 8.3 $\overline {Per\,g}=[0,1]$. Now the fact that $B$ is perfect (statement a) of Theorem 4.1) 
and monotonicity of $\phi$ easily imply that $B\subset \overline {Per\,f}$.

	{\it f) there exist an interval $J\subset I$, an integer $k=n$ or $k=2n$ and a set $\widetilde B=\overline {int\,J\cap B}$ such
that $f^kJ=J,\:f^k \widetilde B=\widetilde B,\:f^i \widetilde B \cap f^j \widetilde B$ contains no more than $1$ point
($0\le i<j<k$), $\bigcup^{k-1}_{i=0}f^i \widetilde B=B$ and $f^k|\widetilde{B}$ is almost conjugate to a mixing interval map
(one can assume that if $k=n$ then $I=J$).}

\noindent
{\bf Proof.} Again consider the sace when the period of the interval $I$ is $1$ and $f|B$ is   
almost conjugate by a monotone map $\phi: I\rightarrow [0,1]$ to a transitive map $g:[0,1]\rightarrow [0,1]$.
If $g$ is in fact mixing then set $k=n=1,\;J=I$; clearly then all the properties from statement b) hold. If $g$ is not mixing then
by Lemma 8.3 there exist such $a\in (0,1)$ that $g[0,a]=[a,1],\;g[a,1]=[0,a],\;g^2|[0,a]$ and $g^2|[a,1]$ are mixing. Set $k=2,
J=\phi^{-1}[0,a]$; again it is easy to see that all the properties from statement f) hold which completes the proof.  

\noindent
{\bf Corollary 8.4[Bl7].} 
{\it If $f:[0,1]\rightarrow [0,1]$ is mixing then there exist a fixed point $a\in (0,1)$ and a sequence of intervals
$\{U_i\}^\infty_{i=-\infty}$ with the following properties:

	1) $U_i\subset U_{i+1}=fU_i\;(\forall i)$;

	2) $\cap U_i=\{a\}$;

	3) for any open $V$ there exists $n=n(V)$ such that $f^nV\supset U_0$;

	4) $\bigcup^\infty_{i=-\infty}U_i\supset (0,1)$.}

\noindent
{\bf Proof.} Follows from Lemmas 8.2 and 8.3.$\Box$

	Let $A(f)\equiv A$  be the set of those from points $0,1$ which have no preimages in $(0,1)$. 

\noindent
{\bf Lemma 8.5[Bl7].}
{\it If $f:[0,1]\rightarrow [0,1]$ is mixing then there are the following possibilities for $A$:

	1) $A=\emptyset$;

	2) $A=\{0\},f(0)=0$;

	3) $A=\{1\},f(1)=1$;

	4) $A=\{0,1\},f(0)=0,f(1)=1$;

	5) $A=\{0,1\},f(0)=1,f(1)=0$.

	Moreover, if $I$ is a closed interval, $I\cap A=\emptyset$, then for any open $U$ there exists $n$ such that $f^mU\supset I$
for $m>n$ (in particular, if $A=\emptyset$ then for any open $U$ there exists $n$ such that $f^nU=[0,1]$.}

\noindent
{\bf Proof.} The map $f$ is surjective; thus $A$ is $f^{-1}$-invariant set which 
together with Lemma 8.3 implies the conclusion.$\Box$

\noindent
{\bf Lemma 8.6[Bl7].}
{\it 1) Let $A\neq \emptyset,\;a\in A,\;f(a)=a$. If $f$ is mixing then there exists a sequence $c_n\rightarrow a,\;c_n\neq a$
of fixed points.

	2)  Let $A=\{0,1\},f(0)=1,f(1)=0$. If $f$ is mixing then there exists a sequence of periodic points $\{c_n\}$ 
of minimal period $2$ such that   
$c_n\rightarrow 0, c_n\neq 0$.}

\noindent
{\bf Proof.} It is sufficient to consider the case $0\in A,f(0)=0$. Suppose that $0$ is an isolated fixed point. Then by the transitivity 
$fx>x$ for some $\eta>0$ and any $x\in (0,\eta)$. 
At the same time $0\in A$ and so $0\not \in f[\eta,1]$. Let $z=\inf f|[\eta,1]$; by the transitivity
$z<\eta$. Then because of the properties of 
$f|[0,\eta]$ we see that in fact $z=\inf_k f^k|[\eta,1]$ and so $[z,1]\subset (0,1]$ is an invariant interval which is
a contradiction.$\Box$

	Let us prove that a mixing map of the interval has the specification property. In fact we introduce a property
which is slightly stronger than the usual specification property (we call it {\it the i-specification}) 
and then prove that mixing maps of the interval have the i-specification. Actually, we need this variant of the specification property 
to make possible the consideration of interval maps on their basic sets which are closely related to mixing maps (see Theorem 4.1).

	We will not repeat the definition of the specification property (see Section 1); instead let us introduce the notion
of the i-specification. To this end we first need the following definition. Let $z\in Per\,f$ have a period $m$. Moreover,
let $f^m[z,z+\eta]$ lie to the left of $z$ and $f^m[z-\eta,z]$ lie to the right of $z$ for some $\eta>0$. Then we say that 
{\it the map $f^m$ at the point $z$ (of period $m$) is reversing}\index{map $f^m$ at the point $z$ (of period $m$) is reversing}; 
otherwise we say that {\it the map $f^m$ at the point $z$ (of period $m$) is non-reversing}\index{map $f^m$ at the point $z$ 
(of period $m$) is non-reversing}.

	Now let $f:I\rightarrow I$ be a continuous interval map. 
The map $f$ is said to have {\it the i-specification property}\index{i-specification property} 
or simply {\it the i-specification} if for any
$\varepsilon>0$ there exists an integer $M=M(\varepsilon)$ such that for any $k>1$, any $k$ points $x_1,x_2,\ldots,x_k\in I$, 
any semi-neighborhoods $U_i\ni x_i$ with $\lambda(U_i)=\varepsilon$,  
any integers $a_1\le b_1<a_2\le b_2<\ldots <a_k\le b_k$ with $a_i-b_{i-1} \ge M,\,2\le i\le k$ and any integer $p$ with 
$p\ge M+b_k-a_1$ there exists a point $x\in I$ of period $p$ such that $f^p$ is non-reversing at the point $x$ and, moreover,  
$d(f^nx,f^nx_i)\le \varepsilon$ for $a_i\le n \le b_i,1\le i\le k$
and $f^{a_i}z\in U_i,\;1\le i\le n$. The additional properties which are required by the i-specification comparing with the usual
specification  
give us the possibility to lift some properties of mixing interval maps (which as we are going to prove have the i-specification)  to
interval maps on basic sets. 

\noindent
{\bf Theorem 8.7[Bl7].}
{\it If a map $f:[0,1]\rightarrow [0,1]$ is mixing then it has the i-specification property.}

\noindent
{\bf Proof.} We will consider some cases depending on the structure of the set $A(f)$ (see Lemma 8.5).

	First we consider the case $A(f)=\emptyset$. Suppose that $\eta>0$. Choose $M=M(\eta)$ such that for any interval $U$ we have
$f^MU=[0,1]$ provided $\lambda(U)>\eta/2$ (which is possible by Lemma 8.5). Let us consider points $x_1,\ldots,x_n$ with
semi-neighborhoods $U_i\ni x_i$ of length $\eta$ and integers $a_1\le b_1<a_2\le b_2<\ldots<a_n\le b_n,\;p\;$ such that
$b_i-a_{i-1}\ge M\:(2\le i\le n),\;p\ge M+b_n-a_1$. From now on without loss of generality we will suppose that $a_1=0$.
We have to find a periodic point $z$ of period $p$ such that $f^p$ is non-reversing at $z$ and, moreover, 
$|f^tz-f^tx_i|\le \eta$ for $a_i\le t\le b_i$ and
$f^{a_i}z\in U_i\;(1\le i\le n)$.

	First let us find an interval $W$ with an orbit which approximates pieces of orbits $\{f^tx_i:a_i\le t\le b_i\}^n_{i=1}$
quite well; we show that one can find W in such a way that $f^pW=[0,1]$. Recall the following      

\noindent
{\bf Property C4}(see Section 2).
{\it Let $U$ be an interval, $x\in U$ be a point, $\lambda(U)\ge \eta>0,\, n>0$. Then there exists an interval $V$ such that
$x\in V\subset U,\,\lambda(f^iV)\le \eta\:(0\le \eta \le n)$ and $\lambda(f^jV)=\eta$ for some $j\le n$.}

	By Property C4 there exists an interval $V_1$ such that $x_1\in V_1\subset U_1,\;\lambda(f^iV_1)\le \eta\;(a_1\le i\le b_1)$
and $\lambda(f^{t_1}V_1)=\eta$ for some $t_1,\;0=a_1\le t_1\le b_1$. Clearly, $[0,1]=f^{a_2-b_1}(f^{b_1}V_1)=f^{a_2-t_1}(f^{t_1}V_1)$
since $a_2-t_1\ge a_2-b_1\ge M$. Then we can find an interval $W_1\subset V_1$ such that $f^{a_2}W_1=U_2$. Repeating this argument
we get an interval $W=[\alpha,\beta]$ such that for any $1\le i\le n$ and $a_i\le t\le b_i$ we have 
$f^tW\subset [f^tx_i-\eta,f^tx_i+\eta],\;f^{a_i}W\subset U_i$ and for some $a_n\le l\le b_n$ we have $\lambda(f^lW)=\eta$. Since 
$p\ge M+b_n-a_1=M+b_n$ we see that $f^pW=f^{p-l}(f^lW)=[0,1]$.

	It remains to show that there exists a periodic point $z\in W$ of period $p$ such that $f^p$ is non-reversing at $z$.  
Suppose that $f^p$ is reversing at all $p$-periodic points in $W$.  
Then it is easy to see that there is only one $p$-periodic point $z\in W$ and $z\in int\,W=(\alpha,\beta)$.
At the same time $\lambda(f^lW)\ge \eta$ and so we may assume that, say, $\lambda([f^lz,f^l\beta])\ge \eta/2$; by the choice of $M$
it implies that $f^p[z,\beta]=f^{p-l}(f^l[z,\beta]=[0,1]$ and hence there is another $p$-periodic point
in $(z,\beta]$ which is a contradiction. It completes the consideration of the case $A(f)=\emptyset$. 

	Consider the case $A(f)=\{0\},f(0)=0$; the other cases may be considered similarly. 
Again suppose that $\eta>0$. We will say that
a point $y$ {\it $\delta$-approximates }\index{$\delta$-approximation of a point} a point $x$  
if $|f^nx-f^ny|\le \delta\;(\forall n)$. Let us prove the following

\noindent
{\bf Assertion 1.}
{\it There exists a closed interval $I$ such that $I\cap A(f)=\emptyset$ and for any $x\in [0,1]$ there exists $y\in I$ which
$\eta/3$-approximates $x$; moreover, if $x\in I$ then we can set $y=x$.}

	Indeed, by Lemma 8.6 we can find two fixed points $0<e<d$ such that $d<\eta/3,f[0,e]\subset [0,d]$. Let us show that
$I=[e,1]$ has the required property. 

	We may assume that $x\in [0,e]$ (otherwise we can set $y=x$). If $orb\,x\subset [0,e]$ then set $y=e$. 
If $orb\,x\not \subset [0,e]$ then first let us choose the smallest $n$ such that $f^nx\not \in [0,e]$. Clearly, $f^nx\in (e,d]$.
Now it is easy to see that there exists $y\in (e,d]$ such that $f^iy\in (e,d]$ for $0\le i\le n-1$ and
$f^ny=f^nx$. Obviously $y$ is the required point which completes the proof of Assertion 1.

	Let $M=M(\eta)$ be an integer such that for any interval $U$ longer than $\eta/6$ we have $f^mU\supset I$ for any
$m\ge M$. To show that $f$ has the i-specification property let us consider points $x_1,\ldots,x_n$ with semi-neighborhoods
$U_i\ni x_i$ of length $\eta$ and integers $0=a_1\le b_1<a_2\le b_2<\ldots<a_n\le b_n,\;p$ such that 
$b_i-a_{i-1}\ge M\:(2\le i\le n),\;p\ge M+b_n-a_1$. We have to find a periodic point $z$ of period $p$ such that $f^p$ is non-reversing
at $z$ and, moreover,  
$|f^tz-f^tx_i|\le \eta$ for $a_i\le t\le b_i$ and
$f^{a_i}z\in U_i\;(1\le i\le n)$.

	First let us find points $y_i\in I$ which $\eta/3$-approximate points $x_i$ and belong to $U_i$ (it is possible by Assertion 1
and the fact that if $x_i\not \in I$ then the only semi-neighborhood of $x_i$ of length $\eta$ is $U_i=[x_i,x_i+\eta)$).
Then choose one-sided semi-neighborhoods $V_i$ of $y_i$ such that $V_i\subset U_i,\lambda(V_i)=\eta/3,V_i\subset I\;(1\le i\le n)$.
Now it is easy to see that one can replace $U_i$ by $V_i$, then repeat the arguments from the case $A(f)=\emptyset$ and get a point 
$z$ with the required properties. This completes the proof.$\Box$
\vspace{.2in}

\noindent
{\large \bf 9. Corollaries concerning periods of cycles}
\vspace{.2in}

	Let us pass to the corollaries concerning periods of cycles of continuous maps of the interval. Theorem Sh1 well-known
properties of the topological entropy imply that $h(f)=h(f|\overline{Per\,f})$. 
However, it is possible to get a set $D$ such that $h(f)=h(f|D)$ using
essentially fewer periodic points of $f$. Indeed, let $A\subset {\Bbb N},\;K_f(A)=\{y\in Per\,f$: minimal period of $y$ belongs to
$A\}$.

\noindent
{\bf Theorem 9.1[Bl4,Bl7].} {\it The following two properties of $A\subset \Bbb N$ are equivalent:

	1) $h(f)=h(f|\overline {K_f(A)})$ for any $f$;

	2) for any $k$ there exists $n\in A$ which is a multiple of $k$.}

\noindent
{\bf Proof.} First suppose that statement 2) holds and prove that it implies statement 1). By the Decomposition Theorem it is enough
to show that $\overline {\cup B_i}\subset \overline {K_f(A)}$ where $\cup B_i$ is the union of all basic sets of $f$. Fix a basic
set $B=B(orb\,I,f)$; then by Theorem 4.1.f) 
we see that there is an interval $J\subset I$, a number $m$ such that $f^mJ=J$, a set $\widetilde B=\overline {int\,J\cap B}$ and
a monotone map $\phi:J\rightarrow [0,1]$ such that $\bigcup^{m-1}_{i=0}f^i \widetilde B=B$ and 
$f^m|\widetilde B$ is almost conjugate by $\phi$ to a mixing map
$g:[0,1]\rightarrow [0,1]$. By Theorem 8.7 the map $g$ has the specification property. Now we need the following easy property
of maps with the specification.

\noindent
{\bf Property X.}
{\it If $\psi:X\rightarrow X$ is a map with the specification and $H\subset \Bbb N$ is infinite then $\overline {K_\psi (H)}=X$.}

	To prove Property X it is necessary to observe first that there exist at least two different $\psi$-periodic orbits.
Now we need to show that for any $z\in X$ there is a point from $K_\psi(H)$ in any open $U\ni z$. To this end we may 
apply the specification property and pick up a point $y\in U$ which first approximates the orbit of $z$ for a lot of time,  
then  approximates one of the previously chosen periodic orbits 
for only one iteration of $f$ and also has the property $\psi^Ny=y$ where $N\in H$ is a large number
(the periodic orbit we consider here should not contain $z$; that is why first needed to find two
distinct periodic orbits). Clearly, taking the appropriate constants and large enough number $N$ from $H$ 
we can see that the minimal period of $y$ is exactly $N$ which completes the proof of Property X.

	Let us return to the proof of Theorem 9.1. Consider the set $A'=\{n: mn\in A\}$. Then by statement 2)
from Theorem 9.1 we see that $A'$ is infinite; so by Property X we have that $K_g(A')$ is dense in $[0,1]$. Now by the properties of
almost conjugations we see that $\widetilde B \subset \overline {K_f(A)}$. It completes the sketch of the proof of 
the fact that statement 2) implies statement 1) of Theorem 9.1. 

	To show that statement 1) implies statement 2) suppose that $A\subset \Bbb N$ is a set such that for
some $k$ there are no multiples of $k$ in $A$. We need to construct a map $f$ such that $h(f)>h(f|\overline {K_f(A)}$.
To this end consider some pm-map $g$ with a periodic interval $I$ of period $k$. Let us construct a new map $f$
which coincides with $g$ on the set $[0,1]\setminus orb_gI$ and may be obtained by changing of the map $g$
only on the set $orb_gI$ in such a way that $orb_gI=orb_fI$ remains the cycle of intervals for the map $f$ as well as 
for the map $g$ and 

	$h(f|orb_f I)>h(f|\{x:f^nx\not\in orb_f I \:(\forall n)\})=h(f|\{x:f^nx\not\in orb_f I \:(\forall n)\})$. 

\noindent Clearly, it is possible and this way we will get a map $f$ such that $h(f)>h(f|\overline{K_f(A)}$. It completes the sketch of
the proof of Theorem 9.1.$\Box$

	Now we are going to study how the sets $\Omega(f),\Omega(f^2),\ldots$ vary for maps with a fixed set of periods of cycles.
In what follows by a period of a periodic point we always mean the minimal period of the point. 
	In [Sh1] A.N. Sharkovskii introduced the notion of L-scheme.

\noindent
{\bf L-scheme.} If there exist a fixed point $x$ and a point $y$ such that either $f^2\le x<y<fy$ or $fy<y<x\le f^2y$
then it is said that $f$ {\it has L-scheme} and points $x,y$ {\it form L-scheme}.

\noindent
{\bf Theorem Sh4[Sh1].} {\it If $f$ has L-scheme then $f$ has cycles of all periods.}

\noindent
{\bf Lemma 9.2.} {\it If $f$ has L-scheme then $h(f)\ge ln 2$.}

\noindent
{\bf Proof.} It follows from the well-known results on the connection between symbolic 
dynamics and one-dimensional dynamical systems (see, for example, [BGMY]).$\Box$

\noindent
{\bf Lemma 9.3[Bl2,Bl7].} {\it Let $f:[0,1]\rightarrow [0,1]$ be a transitive continuous map. Then:

	1) $f^2$ has L-scheme;

	2) $h(f)\ge 1/2\cdot{\ln 2}$;

	3) $f$ has cycles of all even periods.}

\noindent
{\bf Proof.} By Theorem Sh4 and Lemma 9.2 it is sufficient to prove statement 1). Consider some cases.

\noindent
{\bf Case 1.}{\it There exist $0\le a<b\le 1$ such that $fa=a, fb=b$.}

	Assuming $z<fz$ for $z\in (a,b)$ let us prove that $a\in f[b,1]$. Indeed, otherwise $I=[b,1]\cup f[b,1]\neq [0,1]$
is an $f$-invariant interval which contradicts the transitivity. Choose the smallest $c\in [b,1]$ such that $fc=a$; then
$b<c$. It is easy to see that there exists $d\in (a,c)$ such that $fd=c$ and points $a,d$ form L-scheme. 
In other words, we have shown that in this case the map $f$ itself has L-scheme; in particular, 
if $f(0)=0$ or $f(1)=1$ then $f$ has L-scheme.

\noindent
{\bf Case 2.} {\it There exists a fixed point $t\in (0,1)$ such that $fy>y$ for any $y\in [0,t)$ and $fy<y$ for any $y\in (t,1]$.}

	If $f[0,t]=[t,1],\;f[t,1]=[0,t]$ then by Case 1 we may conclude that $f^2$ has L-scheme. So by Lemma 8.3 we may assume
that $f$ is mixing which implies that $f^2$ is transitive and has an $f^2$-fixed point $y\neq t$ (by Lemma 8.2). 
Now Case 1 implies the conclusion.$\Box$

	Note that Lemma 9.3 implies statement d) of Theorem 4.1.

\noindent
{\bf Theorem 9.4[Bl4,Bl7,Bl8].} {\it Let $n\ge 0, k\ge 0$ be fixed, $f$ have no cycles  of period $2^n(2k+1)$. Then:

	1) if $B=B(orb\,I,f)$ is a basic set and $I$ has a period $m$ then $2^n(2k+1)\prec m\prec 2^{n-1}$;

	2) $\Omega(f)=\Omega(f^{2^n})$;

	3) if $f$ is of type $2^l,\:l\le \infty\:$ then $\Omega(f)=\Omega(f^r)\:(\forall r)$.}

\noindent
{\bf Proof.} 1) By the Sharkovskii theorem about the coexistence of periods of cycles for interval maps 
and Theorem 4.1 we have $2^n(2k+1)\prec m$. Suppose $m=2^i,\;i\le n-1$. 
Then by Lemma 9.3 and Theorem 4.1 $f$ has a cycle of period $2^i\cdot 2(2k+1)\prec 2^n(2k+1)$ which is a contradiction.

	2) It is sufficient to prove that if $x\in \Omega(f)\setminus \omega(f)$ then $x\in \Omega(f^{2^n})$; 
indeed, obviously $\omega(f)\in \Omega(f^r)$ for any $r$ 
and so $\omega(f)\subset \Omega(f^{2^n})$.  
By Theorem 3.1 if $x$ belongs to a solenoidal set then $orb\,x$ is infinite and so by Theorem CN $x\in \Omega(f^{2^n})$
(remind that Theorem CN was formulated in Section 4). 
Now let $x\in B'(orb\,I,f)$ where $I$ is chosen by Lemma 4.6. Suppose that $I$ has a period $m$; 
by statement 1) $m=2^nj,\;1\le j$. By Lemma 4.6 it implies that
$x\in \Omega(f^m)\subset \Omega(f^{2^n})$.

	3) Follows from statement 2) and Theorem CN.1).$\Box$
\vspace{.2in}

\noindent
{\large \bf 10. Invariant measures}
\vspace{.2in}

	It is well-known that the specification property has a lot of consequences concerning invariant measures (see, for
example, [DGS]). We summarized some of them in Theorem DGS in subsection 1.10 of Section 1.
In the rest of Section 10 we rely on the results of Sections 2-5 to make use of Theorem 8.7 and Theorem DGS. First we need the
following

\noindent
{\bf Lemma 10.1.} 
{\it Let $f:[0,1]\rightarrow [0,1]$ be continuous, $B=B([0,1],f)\neq \emptyset$ and $f|B$ be mixing. Let also $\eta>0$ and
$x_1,x_2,\ldots,x_m\in Per\,(f|B)$. Then one can find $M=M(\{x_i\}^m_{i=1},\eta)$ such that for any integers
$a_1\le b_1<a_2\le b_2<\ldots<a_m\le b_m,\;p$ with $a_{i+1}-b_i\ge M\;(1\le i\le m-1),\:p\ge M+b_m-a_1$ there
exists a periodic point $z\in B$ of period $p$ such that $f^p$ is non-reversing at $z$ and, moreover, $|f^nz-f^nx_i|\le \eta$ 
for $a_i\le n\le b_i\:(1\le i\le m)$.}

\noindent
{\bf Proof.} First consider the case when $m=2$; let $x_2=y$. For the sake of convenience
let us reformulate our lemma in this situation. Namely,  
$x,y\in Per\,(f|B)$ and we have to find $M=M(\{x,y\},\eta)$ such that for any $a_1\le b_1<a_2\le b_2,\;p$ 
with $a_2-b_1\ge M,p\ge M+b_2-a_1$ there
exists a periodic point $z\in B$ of period $p$ such that $f^p$ is non-reversing at $z$ and, moreover, $|f^nz-f^nx|\le \eta$ 
for $a_1\le n\le b_1$ and $|f^nz-f^ny|\le \eta$ for
$a_2\le n\le b_2$.

	Let us assume that $x$ and $y$ are fixed points; the result in the general situation 
may be deduced from this case or may be proved by the similar arguments. Choose a semi-neighborhood $V$ of $x$ in the following way.
First choose a side $T$ of $x$ such that $x$ is a limit point for $B$ from the side $T$. If $x$ is not an endpoint of some
interval complementary to $B$ then let $V=V_T(x)$ be a semi-neighborhood of $x$ of length smaller than $\eta$. If,
for example, $(x,\alpha)$ is an interval complementary to $B$ then let $V=V_T(x)$ have the properties $f \overline V \not \ni \alpha$ 
and $\lambda(V)<\eta$. Similarly we find a semi-neighborhood $W$ of $y$. We may assume $V\cap W=\emptyset$. 

	By Theorem 4.1 there exist a mixing map $g:[0,1]\rightarrow [0,1]$ and a non-strictly increasing map $\phi:[0,1]\rightarrow [0,1]$
such that $\phi$ almost conjugates $f$ to $g$. We may assume that $\phi(W)=W'$ and $\phi(V)=V'$ have the same length $\delta$ 
and, moreover, $W=\phi^{-1}(W')$ and $V=\phi^{-1}V'$; by the construction $V'$ and $W'$ are semi-neighborhoods of $\phi(x)=x'$ and 
$\phi(y)=y'$ respectively. Furthermore, we may assume that if $x$ is not an endpoint of an interval complementary to $B$ then
$[x'-\delta,x'+\delta]\subset int\,(\phi[x-\eta,x+\eta])$ and the similar property holds for $y$.

	By Theorem 8.7 there exists $M=M(\delta)$ corresponding to the constant $\delta$ in the i-specification property for $g$.
Again we may assume without loss of generality that $a_1=0$. 
Now let $0=a_1\le b_1< a_2\le b_2,\;p$ be integers  with the properties from Lemma 10.1 with this number $M$. Applying Theorem 8.7
to the points $x',y'$ with the semi-neighborhoods $V',W'$ and the integers $0=a_1\le b_1< a_2\le b_2,\;p$ we can find a 
periodic point $z'$ such that  $g^p$ is non-reversing at $z'$ and, moreover, 
$|g^nz'-g^nx'|\le \delta$ for $a_1\le n\le b_1$, $|g^nz'-g^ny'|\le \delta$ for
$a_2\le n\le b_2$ and $z'=g^{a_1}z'\in V',g^{a_2}z'\in W'$.

	Properties of $\phi$ imply that $\phi^{-1}(z')$ is either a point or a closure of an interval complementary to $B$. In the
first case set $z=\phi^{-1}(z')$. In the second case it is easy to see that since $z'$ is a $g$-periodic point of
period $p$ at which $g^p$ is non-reversing then there exists an endpoint $z$ of the interval $\phi^{-1}(z')$ such that $f^pz=z$. 
In any case we get a 
$f$-periodic point $z\in B$ of period $p$ such that $f^p$ is non-reversing at $z$ and $\phi(z)=z'$.

	Let us show that $z$ is the required point. Suppose that $x$ is not an endpoint of an interval complementary to $B$.
Then $|g^nz'-g^nx'|=|g^nz'-x'|\le \delta$ implies $|f^nz-f^nx|=|f^nz-x|<\eta$ by the choice of $\delta$. So we may assume that
$(x,\alpha)$ is an interval complementary to $B$. By the construction $z'=g^{a_1}z'\in V'$ and so $z=f^{a_1}z\in V$. 
Suppose that there exist numbers $r\le b_1$ such that $f^rz\not \in V$ and let $n$ be the smallest such number.  
If $f^nz$ lies to the left of $V$ then $|\phi(f^nz)-x'|=|g^nz'-x'|>\delta$ although by the i-specification property 
$|g^nz'-x'|\le \delta$ (since $n\le b_1$). Thus $f^nz$ lies to the right of $V$ which means that it lies to the right of $\alpha$. 
At the same time $f^{n-1}z\in V,\;fx=x$ and  by the choice of $V$ we have 
$f \overline V\not \ni \alpha$. Clearly, we get to the contradiction and so $f^rz\in V,a_1\le r\le b_1$. Applying the similar
arguments to the point $y$ we obtain the conclusion.

	The proof in case when $m>2$ is similar and left to the reader. $\Box$

\noindent
{\bf Corollary 10.2.}
{\it Let $d_1,\ldots,d_n$ be periodic points belonging to a basic set $B,\;l\in {\Bbb N}$ 
and $\mu=\sum^n_{i=1}\alpha_i\cdot\nu(d_i)$
be an invariant measure. Then $\mu$ can be approximated by CO-measures with supports in $B$ and minimal periods greater than $l$.}

\noindent
{\bf Proof.} We only outline here the proof which is very is similar to that of Proposition 21.8 from [DGS] 
(note that we are going to apply Lemma 10.1 instead of the specification property).

	Namely, suppose that a neighborhood of $\mu$ is given. 
We may assume that $n>1$ and orbits of $d_1,\ldots,d_n$ are pairwise distinct. Choose $\eta$ such
that $dist(orb\,d_i,\:orb\,d_j)>10\eta\:(i\neq j)$. Then approximate the measure $\mu$
by a measure of the same type, i.e. by a measure $\mu'= \sum^n_{i=1}\beta_i\cdot\nu(d_i)$, where $\beta_i$
are sufficiently chosen and very close to $\alpha_i$ rationals. The next step is to construct
a collection of integers $a_1=0<b_1<a_2<b_2<\ldots<a_n<b_n,\:p$ which are required in Lemma 10.1 in such a way that
for any $1\le i\le n$ we have $(b_i-a_i)/p=\beta_i$ and $b_i-a_i\gg M=M(\{d_i\}^n_{j=1},\eta)$; furthermore,
we may assume that $p\gg l$. Take the periodic point $z$ of period $p$ which exists for this collection of integers and
periodic points by Lemma 10.1 and
approximates pieces of orbits of $d_1,\ldots,d_n$. Then because of the choice of $\eta$ it is easy to see that $p$ is the minimal
period of $z$. At the same time  
similarly to the proof of Proposition 21.8 from [DGS] it is easy to see
that in fact the constants may be chosen in such a way that the point $z$ generates the required CO-measure $\nu(z)$; in
other words, we may assume that $\nu(z)$ approximates $\mu$, lying in the previously given neighborhood of $\mu$. 
It completes the proof. $\Box$
 
\noindent
{\bf Theorem 10.3} (cf. Theorem DGS). 
{\it Let $B$ be a basic set. Then the following statements are true.

	1) For any $l \in \Bbb N$ the set $\bigcup_{p\ge l}P_f(p)$ is dense in $M_{f|B}$.

	2) The set of ergodic non-atomic invariant measures $\mu$ with $supp\,\mu=B$ is residual in $M_{f|B}$.

	3) The set of all invariant measures which are not strongly mixing is a residual subset of $M_{f|B}$.

	4) Let $V\subset M_{f|B}$ be a non-empty closed connected set. Then the set of all points $x\in B$ such that $V_f(x)=V$ is
dense in $B$ (in particular, every measure $\mu\in M_{f|B}$ has generic points).

	5) The set of points with maximal oscillation for $f|B$ is residual in $B$.}

\noindent
{\bf Proof.} First observe that if $g$ is a transitive non-strictly periodic map then it is easy to see that Theorem 10.3
holds for $g$ by Theorem DGS, Theorem 8.7 and Lemma 8.3. Now let us pass to the proof of statement 1) assuming that $B$
is a Cantor set.

	Let $B=B(orb\,I,f),\;g$ be a transitive non-strictly periodic map and $\phi$ almost conjugate $f|orb\,I$ to $g$
(maps $\phi$ and $g$ exist by Theorem 4.1). Let $\mu\in M_{f|B}$ and $l\in \Bbb N$. We have to prove that $\mu$ belongs to the
closure of $\bigcup_{p\ge l}P_f(p)$ in $M_{f|B}$.

	The case when $\mu$ is non-atomic is quite clear and we leave it to the reader  
(indeed, it is enough to consider the measure $\mu'\in M_g$ which is the $\phi$-image of $\mu$, apply Theorem DGS
to the measure $\mu'$ and then lift the approximation we found for the measure $\mu'$ to the approximation 
of the measure $\mu$ which is possible since $\mu$ is non-atomic). 
On the other hand it is easy to see that any invariant measure from $M_{f|B}$ may be approximated by a measure $\mu$
of type $\mu=\alpha_0\cdot \tilde\mu+\sum^N_{i=1}\alpha_i\cdot\nu(e_i)$ where $\tilde\mu$
is non-atomic and $N<\infty$. By the non-atomic case we can approximate $\tilde\mu$ by a CO-measure  
$\nu(e_0)$. Applying Corollary 10.2 we can approximate the measure $\sum^N_{i=0}\alpha_i\cdot\nu(e_i)$ by a CO-measure $\nu(c)$
where $c$ is a periodic point with a minimal period $m\ge l$. This completes the proof of statement 1).

	Looking through the proofs of Propositions 21.9-21.21 from [DGS, Section 21] which correspond to statements 2)-5) of
Theorem DGS one can check that they are based on statement 1) of Theorem DGS and the property of invariant measures which is proved
in Corollary 10.2. Repeating the arguments from [DGS, Section 21] one can prove statements 2)-5) of Theorem 10.3.$\Box$

	Property 5) from Theorem 10.3 shows that if $f$ is a transitive interval map then points with maximal oscillation form
a residual subset of the interval. Applying this result we can easily specify Theorem 6.2 as it was 
explained in the proof of this theorem. Namely, in the notation from Theorem 6.2 we need to choose 
a residual subset $\Pi_{orb\,I}$ of any cycle of intervals $orb\,I$ where $f|orb\,I$ is transitive;
to specify Theorem 6.2 one can now choose the set of points with maximal oscillation as the set  $\Pi_{orb\,I}$. It will lead us to the
following 

\noindent  
{\bf Theorem 6.2'}(cf.[Bl1],[Bl8]).
{\it Let $f:[0,1]\rightarrow [0,1]$ be a continuous map without wandering intervals. Then there exists a residual subset
$G\subset [0,1]$ such that for any $x\in G$ one of the following possibilities holds:

	1)$\omega(x)$ is a cycle;

	2) $\omega(x)$ is a solenoid;

	3) $\omega(x)=orb\,I$ is a cycle of intervals and $V_f(x)=M_{f|orb\,I}$.}

\noindent
{\bf Theorem 10.4.} 
{\it Let $\mu$ be an invariant measure. Then the following properties of $\mu$ are equivalent.

	1) There exists $x\in [0,1]$ such that $supp\,\mu\subset\omega(x)$.

	2) The measure $\mu$ has a generic point.

	3) The measure $\mu$ can be approximated by CO-measures.}

\noindent
{\bf Remark.} For non-atomic measures Theorem 10.4 was proved in [Bl4,BL7].

\noindent
{\bf Proof.} Clearly, 2)$\Rightarrow$1). If $\omega(x)$ is a cycle then the implications 1)$\Rightarrow$2) and 1)$\Rightarrow$3)
are trivial. If $\omega(x)$ is a basic set then   the implications 1)$\Rightarrow$2) and 1)$\Rightarrow$3) follow from Theorem 10.3.
The case when $\omega(x)$ is a solenoidal set may be easily deduced from Theorem 3.1;this case is left to the reader.

	It remains to prove that 3)$\Rightarrow$1). Let $\{e_i\}$ be a sequence of periodic points such that 
$\nu(e_i)\rightarrow \mu$. Set $L\equiv\{z:$ for any open $U\ni z$ there exists a sequence $n_k\rightarrow \infty$ such that
$orb\,e_{n_k}\cap U\neq\emptyset\:(\forall k)\}$. Obviously, $L$ is compact, $supp\,\mu\subset L,fL=L$. We may assume that
$e_i\searrow e$. Consider the set $P^R(e)=P^R$; then $L\subset P^R$. By Lemma 2.2 there are the following possibilities
for $P^R$.

	1)$P^R$ {\it is a cycle}. This case is trivial. 

	2) $P^R$ {\it is a solenoidal set}. Then by Theorem 3.1 the fact that $supp\,\mu\subset L\subset P^R$ implies that 
$supp\,\mu=S$ where $S$ is the unique minimal subset if $P^R$. This completes the consideration of the case 2).

	3) $\{P^R\}$ {\it is a cycle of intervals}. Consider two subcases.

	3a) $e$ {\it is the right endpoint of a component $[d,e]$ of $P^R$.} Then $orb\,e_i\cap P^R=\emptyset$ and hence 
$L\subset \partial(P^R)$. Surjectivity of $f|L$ implies that $e\in Per\,f$ and we may assume that $fe=e$.
Clearly, it implies that $\{L\}=\{e\}$ and completes the consideration of the subcase 3a).

	3b) $e\in [z,y)$ {\it where $[z,y]$ is a component of $P^R$.} Then it is easy to see that $L\subset E(P^R,f)$
(the definition of the set $E(orb\,I,f)$ for cycle of intervals $orb\,I$ may be found in Section 4 before Lemma 4.5).
Indeed, we may assume that $orb\,e_i\subset P^R$. Let $\zeta\in L$ and $T$ is a side of $\zeta$ from which points of $orb\,e_{n_k}$
approach the point $\zeta$. Then $T$ is a side of $\zeta$ in the corresponding component of
$P^R$. Consider $P^T(\zeta)$; clearly, $P^T(\zeta)\subset P^R$. At the same time it is easy to see that for any semi-neighborhood
$W_T(\zeta)$ and any $n$ the set $\overline {\bigcup_{i>n} f^iW_T(\zeta)}$ contains some right semi-neighborhood of $e$ which
implies that $P^T(\zeta)\supset P^R$. Finally $P^T(\zeta)=P^R$ and so $L\subset E(P^R,f)$ by the definition. 
Hence by Theorem 4.1 and Lemma 4.5 either $L$ is a cycle or $L\subset B(P^R,f)$. This completes the proof of Theorem 10.4.$\Box$

\noindent
{\bf Corollary 10.5[Bl4,Bl7].} 
{\it CO-measures are dense in all ergodic measures of an interval map}.

\noindent
{\bf Remark.} In [Bl4,Bl7] Corollary 10.5 was deduced from the version of Theorem 10.4 for non-atomic measures proved in [Bl4,Bl7].

\noindent
{\bf Proof.} Left to the reader.$\Box$
\vspace{.2in}

\pagebreak

\noindent
{\large \bf 11. Discussion of some recent results of Block and Coven and Xiong Jincheng}
\vspace{.2in}

	There are some recent papers ([BC], [X]) in which the authors investigate the sets $\omega(f)\setminus \overline {Per\,f}$
and $\Omega(f)\setminus \overline{Per\,f}$. Let us discuss some of their results.

	First observe that by the Decomposition Theorem if $x\in  \omega(f)\setminus \overline {Per\,f}$ then $x\in S_\omega$
for some solenoidal set $S_\omega$ and thus by Theorem 3.1 $\omega(x)=S$ is a minimal solenoidal set. It implies the following
theorem proved in [BC].

\noindent
{\bf Theorem BC.} {\it If $x\in  \omega(f)\setminus \overline {Per\,f}$ then $\omega(x)$ is an infinite minimal set.}

	In [X] some new notions were introduced. Let us recall them. For a set $Y\subset [0,1]$ by $\Lambda(Y)$ we denote
the set $\bigcup_{x\in Y}\omega(x)$; let $\Lambda^1=\Lambda([0,1])=\omega(f),\;\Lambda^2=\Lambda(\Lambda^1)$ etc. 
Obviously $\Lambda^1\supset \Lambda^2\supset\ldots$; let $\Lambda^\infty\equiv\bigcap^\infty_{n=1}\Lambda^n$.

	By $\alpha(x)$\index{$\alpha(x)$} we denote the set of all $\alpha$-limit points of $x$; in other words, $y\in \alpha(x)$ 
if and only if there exist sequences $x_{-i}\rightarrow y$ and $n_i\rightarrow\infty$ such that $f^{n_i}x_{-i}=x$ for any $i$.
A point $y$ is called a $\gamma$-limit point of $x$\index{$\gamma$-limit point} if $y\in \omega(x)\cap \alpha(x)$.
Let $\gamma(x)\equiv \omega(x)\cap \alpha(x)$\index{$\gamma(x)$} and 
$\Gamma(f)\equiv\Gamma\equiv\bigcup_{x\in [0,1]}\gamma(x)$\index{$\Gamma(f)$}.

	In the following lemma we use the notation from the Decomposition Theorem.

\noindent
{\bf Lemma 11.1.} $\Gamma=(\bigcup_i B_i)\cup (\bigcup_{\beta\in \cal A}S^{(\beta)})\cup X_f$.  

\noindent
{\bf Proof.} First let us prove that $\Gamma \supset(\bigcup_i B_i)\cup (\bigcup_{\beta\in \cal A}S^{(\beta)})\cup X_f$.  
Clearly, $X_f\cup(\bigcup_{\beta\in \cal A}S^{(\beta)})\subset \Gamma$
(for $S^{(\beta)}$ it follows for example from the fact that $f|S^{(\beta)}$ is minimal by Theorem 3.1).  
By Theorem 4.1 to prove that $B_i\subset \Gamma\;(\forall i)$ 
it is sufficient to show that $\Gamma(g)=[0,1]$ provided $g:[0,1]\rightarrow [0,1]$ is a transitive map.
Indeed, if $x\in (0,1)$ then by Lemma 8.3 $\alpha(x)=[0,1]$. Thus if $x\in (0,1)$ has a dense orbit in $[0,1]$ then
$\gamma (x)=[0,1]$ and so $\Gamma(g)=[0,1]$. Finally, we may conclude that
$\Gamma \supset (\bigcup_i B_i)\cup (\bigcup_{\beta\in \cal A}S^{(\beta)})\cup X_f$.

	Now let us prove that $\Gamma \subset (\bigcup_i B_i)\cup (\bigcup_{\beta\in \cal A}S^{(\beta)})\cup X_f$
Indeed, 
$\Gamma \subset \omega(f)=(\bigcup_i B_i)\cup (\bigcup_{\beta\in \cal A}S^{(\beta)}_\omega)\cup X_f$ by the definition of $\Gamma$.
So to prove Lemma 11.1  
it remains to show that if $x\in S^{(\beta)}_\omega\setminus S^{(\beta)}$ then $x\not \in \Gamma$ (here $\beta \in \cal A$).
Suppose there exists $z$ such that $x\in \omega(z)\cap \alpha(z)$. Then the fact that
$x\not \in S^{(\beta)}$ implies that $z\not \in Q^{(\beta)}$ because otherwise
$x\in \omega(z)=S^{(\beta)}$ by Theorem 3.1. Hence $\alpha(z)\cap Q^{(\beta)}=\emptyset$
which is a contradiction with the fact that $x\in \alpha(z)\cap Q^{(\beta)}$. It completes the proof of Lemma 11.1.$\Box$

	Let us show how to deduce some of the results of [X] from our results.

\noindent
{\bf Theorem X1[X].}
{\it 1) $\Omega(f)\setminus \Gamma$ is at most countable.

	2) $\Lambda^1\setminus \Gamma$ is either empty or countable.

	3) $\overline {Per\,f}\setminus \Gamma$ is either empty or countable.}

\noindent
{\bf Proof.} 1) By the Decomposition Theorem $\Omega(f)\setminus \overline {Per\,f}$ is at most countable. By Theorem 3.1
$S^{(\beta)}_p\neq S^{(\beta)}$ for at most countable family of solenoidal sets and $S^{(\beta)}_p\setminus S^{(\beta)}$
is at most countable. This implies statement 1).

	2) First recall that $\Lambda^1=\omega(f)$. If $\Lambda^1\setminus \Gamma\neq\emptyset$ then by  the Decomposition Theorem and
Lemma 11.1 there exist a solenoidal set $\omega(z)$ and a point $x\in \omega(z)\setminus S$ where $S$ is the unique minimal set
belonging to $\omega(z)$ (see Theorem 3.1); actually $\omega(z)\setminus S\subset \Lambda^1\setminus \Gamma$. 
Now the fact that $f|\omega(z)$ is surjective implies that $\omega(z)\setminus S$ is
countable and the inclusion $\omega(z)\setminus S\subset \Lambda^1\setminus \Gamma$ implies the conclusion.

	3) Consider the case when $\overline {Per\,f}\setminus \Gamma\neq \emptyset$. Similarly to the proof of statement 2)
we see that then there exists a solenoidal set $Q=\bigcap_i orb\,J_i$ such that $(\overline {Per\,f}\cap Q)\setminus S\neq\emptyset$ 
where $S$ is the unique minimal set belonging to $Q$. Denote $(\overline {Per\,f}\cap Q)$ by $R$.
  
	We are going to prove the fact that $R\setminus S$ is a countable set by repeating the arguments from the proof of statement 2) 
replacing $\omega(z)$ by $R$. However, to this end we need to show that $f|R$ is surjective.
Consider a point $y\in R$ and show that it has $f$-preimages in $R$. 
The fact that
$y\in R\subset \overline {Per\,f}$ implies that
there is a point $z\in \overline {Per\,f}$ such that $fz=y$. Let us prove that $z\in Q^{(\beta)}$.
Suppose that $z\not \in Q^{(\beta)}$. Then the fact that $fz=y$ and Theorem 3.1 imply that there exist an open $U\ni z$ and a number
$N$ such that $U\cap Q^{(\beta)}=\emptyset$ and for any $n>N$ we have $f^nU\subset Q^{(\beta)}$. Clearly, it contradicts the fact
that $z\in \overline {Per\,f}$ and shows that actually $z\in Q^{(\beta)}$; hence $z\in Q^{(\beta)}\cap \overline {Per\,f}=R$
and so $f|R$ is surjective. Now the fact that $R\setminus S$ is a countable set may be proved similarly to 
statement 2). $\Box$ 

\noindent
{\bf Theorem X2[X].} $\Lambda^\infty=\ldots=\Lambda^3=\Lambda^2=\Lambda(\overline {Per\,f})=\Lambda(\Omega(f))=\Gamma$.

\noindent
{\bf Proof.} By Lemma 11.1 and properties of basic sets (Theorem 4.1), solenoidal sets (Theorem 3.1)
and cycles we have $\Lambda(\Gamma)=\Gamma$ and so $\Lambda(\Omega(f))\supset \Gamma$
because $\Gamma \subset \Omega (f)$. 
On the other hand the
Decomposition Theorem and the definition of $\Gamma$ imply that $\Lambda(\Omega(f))\subset \Gamma$; indeed, 
in the notation from the Decomposition Theorem we have $\Lambda(B'_i)\subset B_i$ for all $i$, $\Lambda(X_f)\subset X_f$
and $\Lambda(Q^{(\alpha)})\subset S^{(\alpha)}$ for any $\alpha$ (the last assertion follows from Theorem 3.1). 
So $\Lambda(\Omega(f))=\Gamma=\Lambda(\Gamma)$ which completes the proof.$\Box$.

\noindent
{\bf Theorem X3[X].} {\it The following properties of a map $f$ are equivalent:

	1) the type of $f$ is $2^i,\;i\le \infty$;

	2) every $\gamma$-limit point of $f$ is recurrent.}

\noindent
{\bf Proof}. As we have shown in Section 1 the fact that $f$ has type $2^i,\;i\le \infty$ is equivalent to the absence of basic sets  
(see the part of Section 1 where we discuss the connection between the Misiurewicz theorem on maps with zero entropy
and the ``spectral'' decomposition). So in this case by Theorem X1 we see that 
$\Gamma=(\bigcup_{\beta\in \cal A}S^{(\beta)})\cup X_f$. But by Theorem 3.1 every
point of $S^{(\beta)}$ is recurrent and ,clearly, every point of of $X_f$ is recurrent.
Hence if the type of $f$ is $2^i, i\le \infty$ then every $\gamma$-limit point of $f$ is recurrent. 

	On the other hand if there is a basic set $B$ of $f$ then it is easy to find a non-recurrent point $z\in B$
(it follows, for example, from Theorem 4.1 and Lemma 8.3). Now by Lemma 11.1 $B\subset \Gamma$ which 
shows that there exist non-recurrent points in $\Gamma$ and completes the proof.
$\Box$

\pagebreak

\begin{center} 
{\bf REFERENCES}
\end{center}

\noindent
[AKMcA] R.L. Adler, A.G. Konheim, M.H.McAndrew. {\em Topological entropy,} 
Trans. Amer. Math. Soc., {\bf 114 } (1965), 309-319.
\vspace{.1in}

\noindent
[ALM ] L. Alsed\'a, J. Llibre, M. Misiurewicz. {\em Periodic orbits of maps of Y,} Trans. Amer. Math. Soc., {\bf 313 } (1989 ), 475-538.
\vspace{.1in}

\noindent
[Ba] S. Baldwin. {\em An extension of \u Sarkovski\u i's theorem to the $n$-od,} Ergod. Th. and Dynam. Syst., {\bf 11}, no.2 (1991), 
p.249-271
\vspace{.1in} 

\noindent
[BC ] L. Block, E.M. Coven. {\em $\omega$-limit sets for maps of the interval,} Erg. Th. and Dyn. Syst., {\bf 6 } (1986 ), 335-344.
\vspace{.1in}

\noindent
[BGMY] L. Block, J. Guckenheimer, M. Misiurewicz, L.-S. Young. {\em Periodic orbits and topological entropy of one-dimensional
maps. } In: Global Theory of Dynamical Systems, Lecture Notes in Mathematics, {\bf 819 }, Springer: Berlin( 1980), 18-34.
\vspace{.1in}

\noindent
[Bl1] A.M. Blokh. {\em On the limit behavior of one-dimensional dynamical systems,} Russ. Math. Surv., {\bf 37}, no.1 (1982), 157-158.
\vspace{.1in}

\noindent
[Bl2] A.M. Blokh. {\em On sensitive mappings of the interval,} Russ. Math. Surv., {\bf 37}, no.2 (1982), 203-204.
\vspace{.1in}

\noindent
[Bl3] A.M. Blokh. {\em On the ``spectral'' decomposition for piecewise monotone maps of segment,} Russ. Math. Surv., {\bf 37 }, 
(1982), 198-199.
\vspace{.1in}

\noindent
[Bl4] A.M. Blokh {\em Decomposition of dynamical systems on an interval,} Russ. Math. Surv., {\bf 38}, no. 5 (1983), 133-134.
\vspace{.1in}

\noindent
[Bl5] A.M. Blokh. {\em On the connection between entropy and transitivity for one-\linebreak 
dimensional mappings,} Russ. Math. Surv., 
{\bf 42}, no.5 (1987), 165-166.
\vspace{.1in}

\noindent
[Bl6] A.M. Blokh. {\em A letter to editors,} Russ. Math. Surv., {\bf 42}, no.6 (1987).
\vspace{.1in}

\noindent
[Bl7] A.M. Blokh. {\em On the limit behavior of one-dimensional dynamical systems.1 }(in Russian), Preprint no.1156-82,
{\bf VINITI}, Moscow (1982).
\vspace{.1in}

\noindent
[Bl8] A.M. Blokh. {\em On the limit behavior of one-dimensional dynamical systems.2 }(in Russian), Preprint no.2704-82, 
{\bf VINITI}, Moscow (1982).
\vspace{.1in}

\noindent
[Bl9] A.M. Blokh. {\em On transitive maps of one-dimensional branched manifolds }(in Russian). In: Differential-difference Equations 
and Problems of Mathematical Physics, Kiev (1984), 3-9. 
\vspace{.1in}

\noindent
[Bl10] A.M. Blokh. {\em On some properties of maps of the interval with constant slope }(in Russian). In: Mathematical Physics and 
Functional Analysis, Kiev (1986), 127-136.
\vspace{.1in}

\noindent
[Bl11] A.M. Blokh. {\em On dynamical systems on one-dimensional branched manifolds. 1, 2, 3 }(in Russian):
1, Theory of Functions, Functional Analysis and Applications, Kharkov, {\bf 46} (1986), 8-18; 2, Theory of Functions, Functional Analysis
and Applications, Kharkov, {\bf 47} (1987), 67-77; 3, Theory of Functions, Functional Analysis and Applications, Kharkov, 
{\bf 48} (1987), 32-46.
\vspace{.1in}

\noindent
[Bl12] A.M. Blokh. {\em On $C^0$-continuity of entropy }(in Russian). Preprint (1989).
\vspace{.1in}

\noindent
[Bl13] A.M. Blokh. {\em The spectral decomposition, periods of cycles and Misiurewicz conjecture for graph maps } (1990), submitted
to ``Proceedings of the Conference on Dynamical Systems in G\"ustrow'' (to appear in Lecture Notes in Mathematics). 
\vspace{.1in}

\noindent
[Bl14] A.M. Blokh. {\em  On some properties of graph maps: spectral decomposition, Misiurewicz conjecture and abstract sets of periods,} 
Max-Planck-Institut f\"ur Mathematik, Preprint no.35 (June, 1991).
\vspace{.1in}

\noindent
[Bl15] A.M. Blokh. {\em Periods implying almost all periods, trees with snowflakes, and zero entropy maps,} 
SUNY, Institute for Mathematical Sciences, Preprint no.13 (August, 1991).
\vspace{.1in}

\noindent
[BlL] A.M. Blokh, M.Yu. Lyubich. {\em Non-existence of wandering intervals and structure of topological attractors of one-dimensional
dynamical systems. 2. The smooth case.,} Erg. Th. and Dyn. Syst, {\bf 9 }, no.4 (1989), 751-758. 
\vspace{.1in}

\noindent
[B1] R. Bowen. {\em Entropy for group endomorphisms and homogeneous spaces,} Trans. Amer. Math. Soc., {\bf 153} (1971), 401-413.
\vspace{.1in}

\noindent
[B2] R. Bowen. {\em Periodic points and measures for axiom A-diffeomorphisms,} Trans. Amer. Math. Soc., {\bf 154} (1971), 377-397.
\vspace{.1in}

\noindent
[B3] R. Bowen. {\em Topological entropy for noncompact sets,} Trans. Amer. Math. Soc., {\bf 184} (1973), 125-136.
\vspace{.1in}

\noindent
[B4] R. Bowen. {\em Equilibrium states and the ergodic theory of Anosov diffeomorphisms.} Lecture Notes in Mathematics, 
{\bf 470}, Springer: Berlin (1975).
\vspace{.1in}

\noindent
[BF] R. Bowen, J. Franks. {\em The periodic points of maps of the disk and the interval,} Topology, {\bf 15} (1976), 337-342.
\vspace{.1in}

\noindent
[CE] P. Collet, J.-P. Eckmann. {\em Iterated maps on the interval as dynamical systems. } Progress in Physics, {\bf 1 },  
Birkh\"auser: Boston (1980).
\vspace{.1in}

\noindent
[CN] E.M. Coven, Z. Nitecki. {\em Non-wandering sets of the powers of maps of the interval,} Erg. Th. and Dyn. Syst., {\bf 1} 
(1981), 9-31
\vspace{.1in}

\noindent
[DGS] M. Denker, C. Grillenberger, K. Sigmund. {\em Ergodic theory on compact spaces. } Lecture Notes in Mathematics, {\bf 527 ,} 
Springer: Berlin (1976).
\vspace{.1in}

\noindent
[D] A. Denjoy. {\em Sur les courbes definies par les \'equations differentielles \`a la surface du tore,}
 J. Math. Pures et Appl., {\bf 11} (1932), 333-375.
\vspace{.1in}

\noindent
[Di] E.I. Dinaburg. {\em The relation between topological entropy and metric entropy,} Soviet Math. Dokl., {\bf 11}, no.1 (1970), 13-16.
\vspace{.1in}

\noindent
[F] M. Feigenbaum. {\em Quantitative universality for a class of nonlinear transformations,} J. Stat. Phys., {\bf 19} (1978), 25-52.
\vspace{.1in}

\noindent
[Gu] J. Guckenheimer. {\em Sensitive dependence to initial conditions for one dimensional maps,} Comm. Math. Phys., 
{\bf 70} (1979), 133-160.
\vspace{.1in}

\noindent
[Go] T.N.T. Goodman. {\em Relating topological entropy with measure theoretic entropy,} Bull. Lond. Math. Soc., {\bf 3} (1971), 176-180.
\vspace{.1in}

\noindent
[GZ] B.M. Gurevich, A.S. Zargaryan. {\em A continuous one-dimensional map without maximal measure,} Funct. Anal. and its
Appl., {\bf 20}, no.2 (1986), 60-61.
\vspace{.1in}

\noindent
[H1] F. Hofbauer. {\em The structure of piecewise monotone transformations,} Erg. Th. and Dyn. Syst., {\bf 1}( 1981), 135-143.
\vspace{.1in}

\noindent
[H2] F. Hofbauer. {\em Piecewise invertible dynamical systems,} Probab. Th. Rel. Fields, {\bf 72} (1986), 359-386.
\vspace{.1in}

\noindent
[H3] F. Hofbauer. {\em Generic properties of invariant measures for continuous piecewise monotonic transformations,}
Monat. Math., {\bf 106} (1988), 301-312.
\vspace{.1in}

\noindent
[HR] F. Hofbauer, R. Raith. {\em Topologically transitive subsets of piecewise monotonic maps, which contain no periodic 
points,} Monat. Math., {\bf 107} (1989), 217-240.
\vspace{.1in}

\noindent
[JR1] L. Jonker, D. Rand. {\em Bifurcations in one dimension.1: The non-wandering set,} Inv. Math., {\bf 62} (1981), 347-365.
\vspace{.1in}

\noindent
[JR2] L. Jonker, D. Rand. {\em Bifurcations in one dimension.2: A versal model for bifurcations,} Inv. Math., {\bf 63} (1981), 1-15.
\vspace{.1in}

\noindent
[K] A.N. Kolmogorov. {\em A new metric invariant of transitive dynamical systems and automorphisms of Lebesgue spaces,} 
Dokl. Acad. Nauk SSSR, {\bf 119} (1958), 861-864.
\vspace{.1in}

\noindent
[L] M.Yu. Lyubich. {\em Non-existence of wandering intervals and structure of topological attractors of one-
dimensional dynamical systems. 1. The case of negative Schwarzian derivative,} Erg. Th. and Dyn. Syst., {\bf 9}, no.4 (1989), 737-750.
\vspace{.1in}

\noindent
[Ma] R. Ma\~{n}\'{e}. {\em Ergodic Theory and Differentiable Dynamics. } A Series of Modern Surveys in Mathematics, {\bf 8}, 
Springer:Berlin (1987).
\vspace{.1in} 

\noindent
[MMSt] M. Martens, W. de Melo, S.J. van Strien. {\em Julia-Fatou-Sullivan theory for real one-dimensional dynamics.}  
Preprint (1988).
\vspace{.1in}

\noindent
[MSt] W. de Melo, S.J. van Strien. {\em A structure theorem in one dimensional dynamics,} Ann. of Math., {\bf 129} (1989), 519-546.
\vspace{.1in}

\noindent
[MilT] J. Milnor, W. Thurston. {\em On iterated maps of the interval. }In: Dynamical systems, Lecture Notes in Mathematics,
{\bf 1342}, Springer: Berlin (1988), 465-564.
\vspace{.1in}

\noindent
[Mi1] M. Misiurewicz. {\em Structure of mappings of the interval with zero entropy. } \linebreak 
Preprint IHES (1978).
\vspace{.1in}

\noindent
[Mi2] M. Misiurewicz. {\em Horseshoes for mappings of the interval,} Bull. Acad. Pol. Sci., ser. sci. math., {\bf 27}, no. 2 (1979),
167-169.
\vspace{.1in}

\noindent
[Mi3] M. Misiurewicz. {\em Invariant measures for continuous transformations of $[0,1]$ with zero topological entropy.}
In: Ergodic theory, Lecture Notes in Mathematics, {\bf 729}, Springer: Berlin (1979), 144-152.
\vspace{.1in}

\noindent
[Mi4] M. Misiurewicz. {\em Periodic points of maps of degree one of a circle,} Erg. Th. and Dyn. Syst., {\bf 2} (1982), 221-227.
\vspace{.1in}

\noindent
[Mi5] M. Misiurewicz. {\em Jumps of entropy in one dimension,} Fund. Math., {\bf 132} (1989), 215-226.
\vspace{.1in}

\noindent
[MiS] M. Misiurewicz, W. Szlenk. {\em Entropy of piecewise monotone mappings,} Studia Mathematica, {\bf 67}, no.1 (1980), 45-53.
\vspace{.1in}

\noindent
[Mi\'Sl] M. Misiurewicz, S.V. \'{S}la\'{c}kov. {\em Entropy of piecewise continuous interval maps. } Preprint (1988).
\vspace{.1in}

\noindent
[N1] Z. Nitecki. {\em Periodic  and limit orbits and the depth of the center for piecewise monotone interval maps,} 
Proc. Amer. Math. Soc., {\bf 80} (1980), 511-514.
\vspace{.1in}

\noindent
[N2] Z.Nitecki. {\em Topological dynamics on the interval. } In: Ergodic Theory and dynamical systems, 2.
Progress in Math., {\bf 21}, Birkh\"auser: Boston (1982), 1-73.
\vspace{.1in}

\noindent
[P1] C. Preston. {\em Iterates of maps on an interval. } Lecture Notes in Mathematics, {\bf 999}, Springer: Berlin (1983).
\vspace{.1in}

\noindent
[P2] C. Preston. {\em Iterates of piecewise monotone mappings on an interval. } Lecture Notes in Mathematics, {\bf 1347}, 
Springer, Berlin (1988).
\vspace{.1in}

\noindent
[Sh1] A.N. Sharkovskii. {\em Co-existence of cycles of continuous maps of the line into itself }(in Russian), 
Ukr. Math. J., {\bf 16} (1964), 61-71.
\vspace{.1in}

\noindent
[Sh2] A.N. Sharkovskii. {\em Non-wandering points and the center of a continuous map of the line into itself }(in Ukrainian), 
Dop. Acad. Nauk Ukr. RSR Ser. A (1964), 865-868.
\vspace{.1in}

\noindent
[Sh3] A.N. Sharkovskii. {\em The behavior of a map in a neighborhood of an attracting set }(in Russian), Ukr. Math. J., {\bf 18} (1966),
60-83.
\vspace{.1in}

\noindent
[Sh4] A.N. Sharkovskii. {\em The partially ordered system of attracting sets,} Soviet Math. Dokl., {\bf 7} (1966), 1384-1386.
\vspace{.1in}

\noindent
[Sh5] A.N. Sharkovskii. {\em On a theorem of G.D. Birkhoff }(in Russian), Dop. Acad. Nauk Ukr. RSR Ser. A (1967), 429-432.
\vspace{.1in}

\noindent
[Sh6] A.N. Sharkovskii. {\em Attracting sets containing no cycles }(in Russian), Ukr. Math. J., {\bf 20} (1968), 136-142.
\vspace{.1in}

\noindent
[Li] Shihai Li. {\em Chain recurrent set and turning points,} to appear in Proc. Amer. Math. Soc.
\vspace{.1in}

\noindent
[Si1] K. Sigmund.  {\em Generic properties of invariant measures for axiom-A-diffeo\-mor\-phisms,} Inv. Math., {\bf 11} (1970), 99-109.
\vspace{.1in}

\noindent
[Si2] K. Sigmund. {\em On dynamical systems with the specification property,} Trans. Amer. Math. Soc., {\bf 190} (1974), 285-299.
\vspace{.1in}

\noindent
[S] S. Smale. {\em Differentiable dynamical systems,} Bull. Amer. Math. Soc., {\bf 73} (1967), 747-817.
\vspace{.1in}

\noindent
[Str] S.J. van Strien. {\em On the bifurcation creating horseshoes. }In: Lecture Notes in Mathematics, {\bf 898},  
Springer: Berlin (1980), 316-351.
\vspace{.1in}

\noindent
[W] J. Willms. {\em Asymptotic behavior of iterated piecewise monotone maps,} Erg. Th. and Dyn. Syst., {\bf 8} (1988), 111-131.
\vspace{.1in}

\noindent
[X] J.-C. Xiong. {\em The attracting center of a continuous self-map of the interval,} Erg. Th. and Dyn. Syst., {\bf 8} (1988), 205-213.
\vspace{.1in}

\noindent
[Y] J.C. Yoccoz. {\em Il n'y a pas de contre-exemple de Denjoy analitiques,} C.R. Acad. Sci. Paris, ser. Math., {\bf 298 }, no. 7 (1984), 
141-144.
\vspace{.1in}

\noindent
[Yo] L.-S. Young. {\em A closing lemma on the interval ,} Inv. Math., {\bf 54} (1979), 179-184.

\pagebreak

\begin{theindex}

  \item $(n,\varepsilon)$-separated set, 6
  \item $(x,T)$, 34
  \item $(x,T)_U$, 34
  \item $2$-adic solenoidal set, $2$-adic solenoid, 7
  \item $B'(M,f)$, 9
  \item $B(M,f)$, 8
  \item $C(T)$, 2
  \item $G_\delta$-set, 13
  \item $H(D)$, 7
  \item $K_f(A)$, 17
  \item $M(X)$, 18
  \item $M_T(X) \equiv M_T$, 18
  \item $NS$, 30
  \item $NS_n$, 30
  \item $P^{\cal U}(x)$, 36
  \item $P^{\cal U}_M(x,f)\equiv P^{\cal U}_M$, 36
  \item $P_M(x)$, 36
  \item $P_T(p)$, 18
  \item $Q(\{I_j\}^\infty_{j=0}) \equiv Q$, 7
  \item $R(T)$, 2
  \item $S_\Omega(Q) \equiv S_\Omega$, 7
  \item $S_\omega(Q) \equiv S_\omega$, 7
  \item $S_p(Q) \equiv S_p$, 7
  \item $Sf$, 30
  \item $Si(x)$, 34
  \item $Si_U(x)$, 34
  \item $Sm$, 30
  \item $Sm_n$, 30
  \item $U$-pair, pair in $U$, 34
  \item $V_T(x)$, 18
  \item $X_f$, 11
  \item $\Gamma(f)$, 75
  \item $\Omega$-basic set, 9
  \item $\Omega(T)$, 1
  \item $\alpha(x)$, 75
  \item $\delta$-approximation of a point, 64
  \item $\delta$-measure, 18
  \item $\delta_x$, 18
  \item $\gamma$-limit point, 75
  \item $\gamma(x)$, 75
  \item $\nu(a)$, 18
  \item $\omega (x)$, 1
  \item $\omega$-limit set of a point, 1
  \item $\omega(T)$, 2
  \item $\overline Z$, 2
  \item $\partial Z$, 8
  \item $\succ$-ordering among basic sets, 23
  \item $\widehat U$, 34
  \item $int\, Z$, 8
  \item $orb\,J$, 5
  \item $orb\,x$, 1
  \item $p(B)$, 23
  \item $supp\,\mu$, 18

  \indexspace

  \item almost conjugation, 8

  \indexspace

  \item basic set, 8

  \indexspace

  \item C(f), 20
  \item center of a map, 2
  \item CO-measure, 18
  \item conjugation, 5
  \item cycle, 2
  \item cycle of intervals, 5
  \item cyclical submanifold, 32

  \indexspace

  \item ergodic measure, 19

  \indexspace

  \item generating intervals, 6
  \item generic point for a measure, 19

  \indexspace

  \item homterval, 27

  \indexspace

  \item i-specification property, 62
  \item interval map of type $2^\infty$, 4
  \item interval map of type $m$, 4
  \item invariant measure, 18

  \indexspace

  \item limit point of an orbit from the side, 36
  \item limit side of a point, 36

  \indexspace

  \item map $f^m$ at the point $z$ (of period $m$) is non-reversing, 62
  \item map $f^m$ at the point $z$ (of period $m$) is reversing, 62
  \item map of a constant slope, 28
  \item map which is not degenerate on a side of a point, 35
  \item map with negative Schwarzian, 30
  \item maximal oscillation, 19
  \item minimal, 6

  \indexspace

  \item non-atomic measure, 18
  \item non-strictly periodic map, 8
  \item non-wandering point, 1
  \item non-wandering set, 1

  \indexspace

  \item o-extremum, 52
  \item orbit of a point, 1

  \indexspace

  \item pair, 34
  \item pair belonging to the image of another pair, 35
  \item Per\,T, 2
  \item period of a basic set, 23
  \item period of a periodic point, 1
  \item periodic interval, 5
  \item periodic point, 1
  \item piecewise-monotone map, pm-map, 20

  \indexspace

  \item recurrent point, 2
  \item residual set, 13

  \indexspace

  \item Schwarzian derivative, 30
  \item semiconjugation, 5
  \item set of genus 0, 3
  \item set of genus 1, 3
  \item set of genus 2, 3
  \item Sharkovskii ordering, 4
  \item side of a point in an interval, 35
  \item smooth map of an interval, 30
  \item solenoid, 6
  \item solenoidal set, 6
  \item source side, 39
  \item specification property, 16
  \item strongly mixing measure, 19
  \item support of $\mu$, 18

  \indexspace

  \item topological entropy, 6
  \item topologically expanding, expanding, 27
  \item topologically generic, 13
  \item topologically mixing, mixing map, 6
  \item transitive map, 6
  \item turning point, 20

  \indexspace

  \item unimodal map, 27

  \indexspace

  \item wandering interval, 12
  \item weak cycle of intervals, 5
  \item weakly periodic interval, 5

\end{theindex}

\end{document}